\documentclass[12pt]{article}
\usepackage{amsmath,amssymb,amsthm}
\author{I.\,A.\,Dynnikov\\
\parbox{0.7\textwidth}{\normalsize
\begin{center}Dept.\ of Mech.\ \& Math.\\
Moscow State University\\
Moscow 119992 GSP-2, Russia\\
e-mail: dynnikov@mech.math.msu.su\end{center}}}
\title{Arc-presentations of links.\\
Monotonic simplification}
\date{}

\def\integer{{\mathbb Z}}\def\real{{\mathbb R}}\def\disk{{\cal D}}
\newtheorem{theorem}{Theorem}
\newtheorem{lemma}{Lemma}
\newtheorem{prop}{Proposition}
\newtheorem{coro}{Corollary}
\theoremstyle{remark}
\newtheorem{rem}{Remark}

\begin{document}
\maketitle
\abstract{In the beginning of 90's J.\,Birman and W.\,Menasco
worked out a nice technique for studying links
presented in the form of a closed braid.
The technique is based on certain foliated surfaces
and uses tricks similar to those that were
introduced earlier by D.\,Bennequin.
A few years later P.\,Cromwell adapted Birman--Menasco's
method for studying so-called arc-presentations of links and established
some of their basic properties.
Here we exhibit a further development of that technique
and of the theory of arc-presentations, and
prove that any arc-presentation of the unknot
admits a (non-strictly) monotonic simplification by
elementary moves; this yields a simple algorithm
for recognizing the unknot.
We show also that the problem of recognizing split links
and that of factorizing a composite link
can be solved in a similar manner. We also define two easily checked
sufficient conditions for knottedness.
}

\tableofcontents

\section*{Introduction}
The most important problem of knot theory is the
classification of knots and links.
To construct an algorithm that, for any two given links, decides
whether they are isotopic or not is known to be a very hard problem.
A solution published in~\cite{matvbook} uses contributions of
many mathematicians including W.\,Haken, F.\,Waldhausen, W.\,Jaco,
P.\,Shalen, K.\,Johannson, G.\,Hemion,
W.\,Thurston, and S.\,Matveev. The previously announced solution
by G.\,Hemion was shown to be incomplete by S.\,Matveev.
It is not difficult to construct
an algorithm that enumerates all diagrams of links, which,
in conjunction with the previously mentioned algorithm, gives
the possibility of generating a list of diagrams in which
each isotopy class of links is presented exactly once and,
for any other diagram, to locate the corresponding isotopy
class in the list.

The algorithm for recognizing links mentioned above is extremely complex
and gives only a theoretical solution of the problem.
It actually deals with three-manifolds which are link complements
provided with a meridian-longitude pair in
each connected component of the boundary.
Up to now, this algorithm cannot be used in practice because of the huge
number of operations needed for its implementation,
even for links with a small number of crossings.

At the same time, a theoretical solution of the problem will be
given, if, for some computable function $f(n)$ (say, $\exp(\exp(\exp(n!)))$),
one shows that any two diagrams of the same link type
that have $\leqslant n$ crossings can be obtained from each other
by $\leqslant f(n)$ Reidemeister moves. However, an
explicit formula for such
a function is still unknown.

It seems to be widely understood that the following
problems, which are important ingredients of
the general classification problem, are simpler
from the practical point of view:
\begin{enumerate}
\item recognizing the trivial knot;
\item recognizing a split link and presenting it as the distant union
of two non-empty links;
\item recognizing a composite non-split link and presenting
it as the connected sum of two non-trivial links.
\end{enumerate}
We shall call the aggregate of these three problem the \emph{decomposition
problem}, since a solution of all three allows one to
express any given link in terms of prime non-split
links and the unknot by using the connected sum and distant union operations.
Such a decomposition is known to be essentially unique, but to solve the
general classification problem, one must also construct an algorithm
for comparing isotopy classes of two prime non-split
links that are presented by their diagrams.

The decomposition problem was solved by W.\,Haken~\cite{haken} and
H.\,Schubert~\cite{shubert} by using Haken's method
of normal surfaces, which plays
a very important r\^ole in modern three-dimensional topology.
(In particular, normal surfaces are used in the
algorithm for recognizing Haken manifolds mentioned above
and in the Rubinstein--Thompson algorithm for recognizing
the three-sphere).

Finding normal surfaces seems to be an exponentially hard
problem. At least, the known realizations are that hard, both theoretically
and in practice. So, any algorithm that uses Haken's technique
has almost no chance to be implementable in a reasonable time.

In recent years, a few attempts have been made to find an algorithm
for recognizing the trivial knot by using some \emph{monotonic
simplification}. By the latter we mean the following. One chooses a
way of presenting knots by diagrams of certain type and
introduces the notion of \emph{complexity} $c(D)$
of a diagram. For the chosen type of presentation, one also fixes
a set of \emph{elementary moves} that do not alter the topological
type of a knot. Then, for a given diagram $D$ of a knot, one
searches for a sequence of diagrams
\begin{equation}\label{diagrams}
D_0=D,D_1,D_2,\ldots
\end{equation}
such that, for
each $i$, the diagram $D_{i+1}$ is obtained from $D_i$ by an
elementary move and we have $c(D_{i+1})\leqslant c(D_i)$
(or $c(D_{i+1})<c(D_i)$, in which case we shall speak of a
\emph{strictly monotonic simplification}). If there
exists an algorithm that, for any $D$, produces such a sequence
$D_0=D,D_1,\dots,D_N=D'$ in which the last diagram $D'$
cannot be simplified anymore,
we shall say that $D$ is monotonically simplifiable to $D'$.
Such an algorithm certainly exists if, for any $n$,
there are only finitely many diagrams $D$ such that $c(D)<n$
(which will be the case for arc-presentations) because,
in this case, we can find \emph{all} monotonic simplification
sequences in finite time by using an exhaustive search.
If, for any diagram $D$ of the unknot, the diagram $D'$
obtained from $D$ by a simplification algorithm is always
the \emph{trivial diagram}, we shall speak about the
recognition of the unknot by monotonic simplification.

It is well known that, in the case of ordinary planar diagrams
with Reidemeister moves as elementary moves and
the crossing number as the measure of complexity,
recognizing the unknot by monotonic simplification is not
possible because there exist diagrams of the unknot that
cannot be simplified to the trivial circle by using only
Reidemeister moves not increasing the crossing number.

According to A.\,E.\,Hatcher's solution of the Smale conjecture~\cite{hat},
there is no topological obstruction to the existence of
a flow $\omega$ on the space of knots such that
any unknotted curve will evolve under $\omega$
to a round circle. Such a flow
may be the gradient flow of some energy function, which
plays the r\^ole of the measure of complexity. Appropriate
discretization of the flow may lead to a strict monotonic
simplification algorithm.
Some investigation of the gradient flows on the knot space and
numerical experiments have been
done in the literature (see, {\it e.g.},~\cite{oha}, \cite{fhw},
\cite{kusner}).

A monotonic simplification procedure for spines of three-manifolds
has been worked out by S.\,V.\,Matveev in~\cite{matv}
for recognizing three-manifolds. It also can be used
for recognizing the unknot: for a given knot one first
constructs a spine of the knot complement and then applies
the simplification algorithm. Numerical experiments
of this kind has been made by E.\,Fominykh.

In \cite{dyn}, the author introduced three-page presentations
of knots. The complexity of a three-page knot
is defined to be the number of vertices at
the binding line of the three-page book. A simplification
procedure for three-page knots has been
tested on a series of examples.

In all three mentioned situations, it has not been proved that
a monotonic simplification to a trivial circle
is always possible for any presentation of the unknot,
but no counterexamples have been found either,
though attempts have been made.

The first successful attempt to find a measure
of complexity with respect to which any presentation
of the unknot admits a monotonic simplification
was made by J.\,Birman and W.\,Menasco in~\cite{bm5}.
The authors studied presentations of a knot as the closure
of a braid and showed that the number of strands in a non-trivial
braid presenting the unknot can always be reduced
by certain type moves. The moves include the Markov moves
that do not increase the number of strands (\emph{i.e.}, braid
conjugations and ``destabilizations'', which are the Markov moves
that reduce the number of strands)
and so-called exchange moves, which preserve the number of strands.
The ideas of~\cite{bm5} and the known solution of the conjugacy
problem for braid groups were used by J.\,Birman and M.\,Hirsch in~\cite{bh}
to construct a new algorithm for recognizing the unknot.
This algorithm, however, does not have the form a monotonic simplification
algorithm in our sense, and it is not clear whether
the straightforward idea of monotonic simplification can be used here,
since no algorithm is known to decide whether a given braid conjugacy
class admits a destabilization after finitely may exchanges.

The problem of recognizing the unknot happens to be very closely
related to two other problems: recognizing split links and
recognizing composite links. In many cases,
one and the same technique allows to solve all three problems.
The distant union and connected sum operations
can be defined at the level of diagrams.
But a split (or composite) link can usually be presented
by a non-split (respectively, prime) diagrams.
So, the general question is this:
for a given diagram $D$, can one apply finitely many
elementary moves without increasing the complexity
so that the final diagram is obtained
by the distant union and connected sum operations
from trivial diagrams and diagrams of prime non-split links?
J.\,Birman and W.\,Menasco have shown in~\cite{bm4} that
the answer is positive in the case of closed braids
(there is a gap in the proof of the result on composite links,
but the assertion is true).

At first look, arc-presentations, which we study in this paper,
seem to have nothing in common with
braids. However, the braid foliation technique that
J.\,Birman and W.\,Menasco developed in a series
of papers including~\cite{bm4} and \cite{bm5}
turned out to be very well adopted for
studying arc-presentations of knots.
The fact that the Birman--Menasco foliated surface technique
can be extended to arc-presentations
was discovered by P.~Cromwell in paper~\cite{cromwell},
where he proved the additivity of
the arc index under the connected sum operation. The possibility of
a monotonic simplification and other algorithmic questions are
not discussed in~\cite{cromwell}, but a sequence
of moves preserving the complexity appears in the proof.
The generalized type~I move of~\cite{cromwell} can be easily
decomposed into elementary moves
(Proposition~\ref{genexchange} of this paper).
In conjunction with this remark, the arguments of~\cite{cromwell}
would suffice for a proof of recognizibility of split
links by monotonic simplification. As for
recognizing composite links, there is a gap in
the final part of the proof of the main result of~\cite{cromwell},
which is covered in this paper. We also provide some
technical details that are not mentioned in~\cite{cromwell}.

In addition to this, we extend the foliated surface technique
to spanning discs of arc-presentations of the unknot.
Compared to the case of a splitting or factorizing sphere,
spanning discs require delicate care of their behaviour
near the boundary.

In the earlier version of the present manuscript, the treatment
of the boundary was not careful enough, and a gap occurred
in the proof of the main result. The mistake was discovered
by W.\,Menasco and A.\,Sikora
who also suggested an idea for filling the gap.
We use their suggestion in modified form.

The main result of this paper is the following.

\begin{theorem}\label{th1}
The decomposition problem of arc-presentations
is solvable by monotonic simplification.
\end{theorem}

\begin{figure}[ht]\caption{A rectangular diagram}\label{example}\vskip20pt
\centerline{\begin{picture}(150,150)
\put(80,0){\circle{3}}\put(82,0){\line(1,0){66}}\put(150,0){\circle{3}}
\put(150,2){\line(0,1){131}}\put(150,135){\circle{3}}
\put(148,135){\line(-1,0){20}}\put(122,135){\line(-1,0){10}}
\put(110,135){\circle{3}}\put(110,133){\line(0,-1){16}}
\put(110,115){\circle{3}}\put(108,115){\line(-1,0){66}}
\put(40,115){\circle{3}}\put(40,117){\line(0,1){31}}
\put(40,150){\circle{3}}\put(38,150){\line(-1,0){36}}
\put(0,150){\circle{3}}\put(0,148){\line(0,-1){126}}\put(0,20){\circle{3}}
\put(2,20){\line(1,0){75}}\put(83,20){\line(1,0){40}}\put(125,20){\circle{3}}
\put(125,22){\line(0,1){121}}\put(125,145){\circle{3}}
\put(123,145){\line(-1,0){71}}\put(50,145){\circle{3}}
\put(50,143){\line(0,-1){11}}\put(50,130){\circle{3}}
\put(48,130){\line(-1,0){5}}\put(37,130){\line(-1,0){10}}
\put(25,130){\circle{3}}\put(25,128){\line(0,-1){86}}\put(25,40){\circle{3}}
\put(27,40){\line(1,0){51}}\put(80,40){\circle{3}}\put(80,38){\line(0,-1){36}}
\end{picture}}
\end{figure}
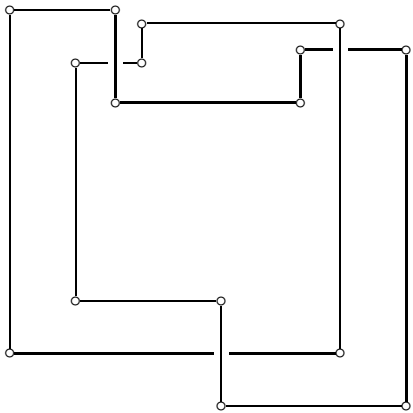

In order to make this claim less abstract in this superficial
introduction, we explain it in a very elementary language.
Let us call a \emph{rectangular diagram} an ordinary
planar diagram $D$ of a link satisfying the following
restrictions:
\begin{enumerate}
\item $D$ consists only of vertical and horizontal straight
line segments, which we call \emph{edges};
\item at each crossing of $D$, the vertical arc is
overcrossing and the horizontal one undercrossing;
\item no two edges are collinear (see Fig.~\ref{example}).
\end{enumerate}
Two such diagrams are called \emph{combinatorially equivalent} if they
are isotopic in the plane via an ambient isotopy
$h$ that has the form $h(x,y)=(f(x),g(y))$.

Since the type of a crossing is determined by its
position, there is no need to indicate in figures
which arc is overcrossing and which one undercrossing.

Clearly, any planar diagram of a link is isotopic to a
rectangular diagram (see below).

The following transforms of rectangular diagrams
are \emph{elementary moves}:
\begin{enumerate}\setcounter{enumi}{-1}
\item cyclic permutation of horizontal (vertical)
edges (see Fig.~\ref{cyclicv});
\begin{figure}[ht]\caption{Cyclic permutation of vertical edges}
\label{cyclicv}\vskip20pt
\centerline{\begin{picture}(300,100)
\multiput(10,0)(8,0){13}{\line(1,0){4}}
\multiput(10,100)(8,0){13}{\line(1,0){4}}
\multiput(10,0)(0,8){13}{\line(0,1){4}}
\multiput(110,0)(0,8){13}{\line(0,1){4}}
\put(30,70){\circle{3}}\put(28,70){\line(-1,0){26}}\put(0,70){\circle{3}}
\put(0,68){\line(0,-1){46}}\put(0,20){\circle{3}}\put(2,20){\line(1,0){96}}
\put(100,20){\circle{3}}\multiput(190,0)(8,0){13}{\line(1,0){4}}
\multiput(190,100)(8,0){13}{\line(1,0){4}}
\multiput(190,0)(0,8){13}{\line(0,1){4}}
\multiput(290,0)(0,8){13}{\line(0,1){4}}
\put(210,70){\circle{3}}\put(212,70){\line(1,0){86}}\put(300,70){\circle{3}}
\put(300,68){\line(0,-1){46}}\put(300,20){\circle{3}}
\put(298,20){\line(-1,0){16}}\put(280,20){\circle{3}}
\put(120,48){\hbox to 60pt{$\leftarrow$\hss$-$\hss$-$\hss
$-$\hss$-$\hss$-$\hss$-$\hss$\rightarrow$}}
\end{picture}}
\end{figure}
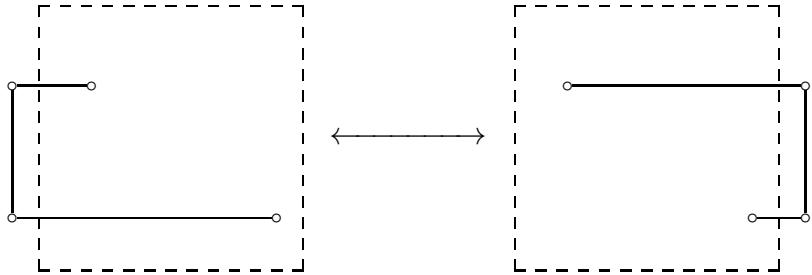
\item stabilization and destabilization (Fig.~\ref{(de)stab});
\begin{figure}[ht]\caption{Stabilization and destabilization moves}
\label{(de)stab}\vskip20pt
\centerline{\begin{picture}(340,60)
\put(0,50){\line(1,0){58}}\put(60,50){\circle{3}}\put(60,52){\line(0,1){6}}
\put(60,60){\circle{3}}\put(58,60){\line(-1,0){6}}\put(50,60){\circle{3}}
\put(50,58){\line(0,-1){58}}
\put(70,33){\hbox to 60pt{$\leftarrow$\hss$-$\hss$-$\hss$-$\hss
$-$\hss$-$\hss$-$\hss$-$}}
\put(70,40){\hbox to 60pt{\hss\scriptsize stabilization\hss}}
\put(70,23){\hbox to 60pt{$-$\hss$-$\hss
$-$\hss$-$\hss$-$\hss$-$\hss$-$\hss$\rightarrow$}}
\put(70,17){\hbox to 60pt{\hss\scriptsize destabilization\hss}}
\put(140,50){\line(1,0){48}}\put(190,50){\circle{3}}
\put(190,48){\line(0,-1){48}}
\put(210,33){\hbox to 60pt{$-$\hss$-$\hss$-$\hss
$-$\hss$-$\hss$-$\hss$-$\hss$\rightarrow$}}
\put(210,40){\hbox to 60pt{\hss\scriptsize stabilization\hss}}
\put(210,23){\hbox to 60pt{$\leftarrow$\hss$-$\hss$-$\hss
$-$\hss$-$\hss$-$\hss$-$\hss$-$}}
\put(210,17){\hbox to 60pt{\hss\scriptsize destabilization\hss}}
\put(280,50){\line(1,0){38}}\put(320,50){\circle{3}}
\put(320,48){\line(0,-1){6}}\put(320,40){\circle{3}}
\put(322,40){\line(1,0){6}}\put(330,40){\circle{3}}
\put(330,38){\line(0,-1){38}}
\end{picture}}
\vskip20pt
\centerline{\begin{picture}(340,60)
\put(0,50){\line(1,0){38}}\put(40,50){\circle{3}}\put(40,52){\line(0,1){6}}
\put(40,60){\circle{3}}\put(42,60){\line(1,0){6}}\put(50,60){\circle{3}}
\put(50,58){\line(0,-1){58}}
\put(70,33){\hbox to 60pt{$\leftarrow$\hss$-$\hss$-$\hss$-$\hss
$-$\hss$-$\hss$-$\hss$-$}}
\put(70,40){\hbox to 60pt{\hss\scriptsize stabilization\hss}}
\put(70,23){\hbox to 60pt{$-$\hss$-$\hss
$-$\hss$-$\hss$-$\hss$-$\hss$-$\hss$\rightarrow$}}
\put(70,17){\hbox to 60pt{\hss\scriptsize destabilization\hss}}
\put(140,50){\line(1,0){48}}\put(190,50){\circle{3}}
\put(190,48){\line(0,-1){48}}
\put(210,33){\hbox to 60pt{$-$\hss$-$\hss$-$\hss
$-$\hss$-$\hss$-$\hss$-$\hss$\rightarrow$}}
\put(210,40){\hbox to 60pt{\hss\scriptsize stabilization\hss}}
\put(210,23){\hbox to 60pt{$\leftarrow$\hss$-$\hss$-$\hss
$-$\hss$-$\hss$-$\hss$-$\hss$-$}}
\put(210,17){\hbox to 60pt{\hss\scriptsize destabilization\hss}}
\put(280,50){\line(1,0){58}}\put(340,50){\circle{3}}
\put(340,48){\line(0,-1){6}}\put(340,40){\circle{3}}
\put(338,40){\line(-1,0){6}}\put(330,40){\circle{3}}
\put(330,38){\line(0,-1){38}}
\end{picture}}
\end{figure}
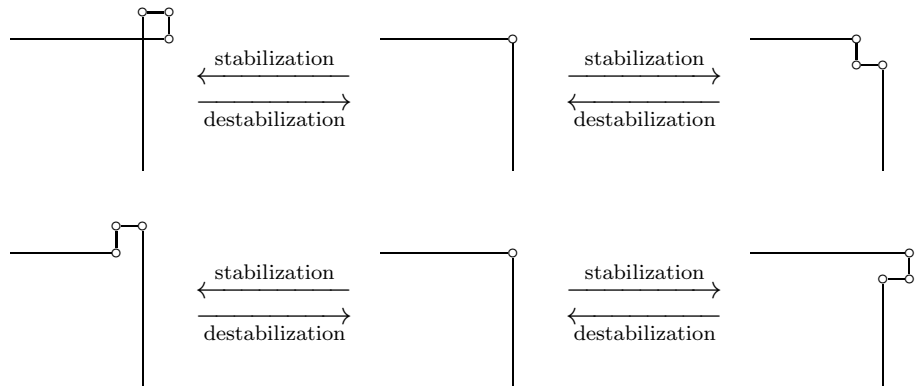
\item interchanging neighbouring edges if their pairs of
endpoints do not interleave (Fig.~\ref{inthor}).
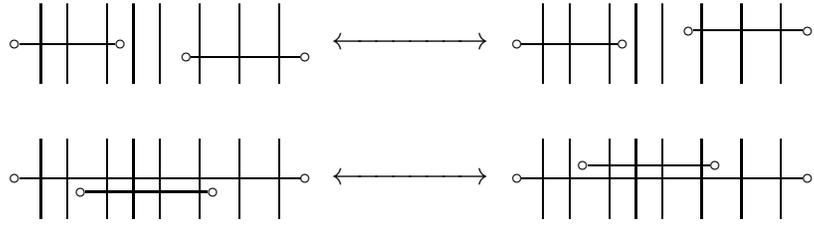
\begin{figure}[ht]\caption{Interchanging horizontal edges}
\label{inthor}\vskip20pt
\centerline{\begin{picture}(300,30)
\put(0,15){\circle{3}}\put(2,15){\line(1,0){36}}\put(40,15){\circle{3}}
\put(65,10){\circle{3}}\put(67,10){\line(1,0){41}}\put(110,10){\circle{3}}
\put(10,0){\line(0,1){30}}\put(20,0){\line(0,1){30}}
\put(35,0){\line(0,1){30}}\put(45,0){\line(0,1){30}}
\put(55,0){\line(0,1){30}}\put(70,0){\line(0,1){30}}
\put(85,0){\line(0,1){30}}\put(100,0){\line(0,1){30}}
\put(120,13){\hbox to 60pt{$\leftarrow$\hss$-$\hss$-$\hss
$-$\hss$-$\hss$-$\hss$-$\hss$\rightarrow$}}
\put(190,15){\circle{3}}\put(192,15){\line(1,0){36}}\put(230,15){\circle{3}}
\put(255,20){\circle{3}}\put(257,20){\line(1,0){41}}\put(300,20){\circle{3}}
\put(200,0){\line(0,1){30}}\put(210,0){\line(0,1){30}}
\put(225,0){\line(0,1){30}}\put(235,0){\line(0,1){30}}
\put(245,0){\line(0,1){30}}\put(260,0){\line(0,1){30}}
\put(275,0){\line(0,1){30}}\put(290,0){\line(0,1){30}}
\end{picture}}
\vskip20pt
\centerline{\begin{picture}(300,30)
\put(0,15){\circle{3}}\put(2,15){\line(1,0){106}}\put(110,15){\circle{3}}
\put(25,10){\circle{3}}\put(27,10){\line(1,0){46}}\put(75,10){\circle{3}}
\put(10,0){\line(0,1){30}}\put(20,0){\line(0,1){30}}
\put(35,0){\line(0,1){30}}\put(45,0){\line(0,1){30}}
\put(55,0){\line(0,1){30}}\put(70,0){\line(0,1){30}}
\put(85,0){\line(0,1){30}}\put(100,0){\line(0,1){30}}
\put(120,13){\hbox to 60pt{$\leftarrow$\hss$-$\hss$-$\hss
$-$\hss$-$\hss$-$\hss$-$\hss$\rightarrow$}}
\put(190,15){\circle{3}}\put(192,15){\line(1,0){106}}\put(300,15){\circle{3}}
\put(215,20){\circle{3}}\put(217,20){\line(1,0){46}}\put(265,20){\circle{3}}
\put(200,0){\line(0,1){30}}\put(210,0){\line(0,1){30}}
\put(225,0){\line(0,1){30}}\put(235,0){\line(0,1){30}}
\put(245,0){\line(0,1){30}}\put(260,0){\line(0,1){30}}
\put(275,0){\line(0,1){30}}\put(290,0){\line(0,1){30}}
\end{picture}}
\end{figure}
\end{enumerate}
Figures~\ref{cyclicv}--\ref{inthor} illustrate particular cases
of the moves. To obtain all the other cases, one should apply
all possible rotations by $\pi k/2$, $k\in\integer$ and reflections
in horizontal and vertical lines if necessary.

The \emph{complexity} $c(D)$ of a rectangular diagram
$D$ is the number of
vertical edges in $D$. The only operations changing the
complexity are stabilization (which increases it by $1$)
and destabilization (which is inverse to stabilization).
We would like to stress that the number of crossings is \emph{not}
taken into account in the definition of complexity.
The number of crossings of a rectangular diagram of complexity $n$
is bounded from above by $(n-1)^2/2$.
The simplest possible rectangular diagram, the \emph{trivial
diagram}, is an ordinary rectangle. Its complexity is $2$.

A rectangular diagram $D$ is a \emph{distant union} of diagrams
$D_1$, $D_2$ if there is a vertical line $l$ that does not
intersect $D$ and splits it into two non-trivial parts which
are equivalent to $D_1,D_2$. We write $D=D_1\sqcup D_2$ in this case.

A rectangular diagram $D$ is a \emph{connected sum} of diagrams $D_1$ and
$D_2$ if there exists a rectangular disk $\mathcal R$ with vertical
and horizontal sides whose boundary
$\partial\mathcal R$ meets $D$ twice so that the rectangular
diagrams $D_1$ and $D_2$ (in either order) are obtained from
$\mathcal R\cap D$ and $(\real^2\setminus\mathcal R)
\cap D$ by connecting the two points $(\partial\mathcal R)\cap D$
by a simple path lying outside or inside $\mathcal R$, respectively.
In this case, we say that the diagram $D$ is \emph{composite} and
write $D=D_1\#D_2$. In particular, a rectangular diagram $D$
is composite if there is a vertical or horizontal line
intersecting $D$ exactly twice such that on both its sides
there are at least four vertices of $D$.

As we shall see in~\ref{arcandrect},
rectangular diagrams are just another way of
thinking about arc-presentations.
Thus, the assertion of Theorem~\ref{th1} is actually this.
For any given rectangular diagram $D$,
one can find a sequence of elementary moves
$D\mapsto D_1\mapsto D_2\mapsto\ldots\mapsto
D_N$, not including stabilization, such that the
final diagram $D_N$ is obtained by the connected sum and
distant union operations from diagrams of prime non-trivial
non-split links and trivial diagrams.

The idea of specifying arc-presentations by drawing
the corresponding rectangular diagrams
is not new. It is already present in P.\,Cromwell's work~\cite{cromwell},
where he calls them ``loops and lines'' diagrams and uses them to show
that any link has an arc-presentation. A similar result is
proved in the manuscript~\cite{phd}
for graphs embedded in $\real^3$ by using
an appropriate generalization of the notion of rectangular diagram.

\section{Preliminaries}
\label{prelim}
\subsection{Definition of arc-presentations}
We shall regard a three-dimensional sphere $S^3$ as the join
of two circles: $S^3=S^1*S^1$, and use the coordinate system
$(\varphi,\tau,\theta)$ on $S^3$, where $\varphi$ and $\theta$
are coordinates on the circles, $\varphi,\theta\in\real/(2\pi\integer)$,
and $\tau$ takes values in the interval $[0,1]$.
We have $(\varphi,0,\theta_1)\sim(\varphi,0,\theta_2)$
and $(\varphi_1,1,\theta)\sim(\varphi_2,1,\theta)$
for all $\varphi,\varphi_1,\varphi_2,\theta,\theta_1,\theta_2\in S^1$.
In order to distinguish between the circles $\tau=0$ and $\tau=1$, we
shall denote them by $S^1_\varphi$ and $S^1_\theta$, respectively.
The circle $S^1_\varphi$ will be called the \emph{binding circle}.
We shall denote by $\disk_t$ the open disk defined by $\theta=t$, $\tau>0$.
Such a disk will be called a \emph{page}.
We have $\partial\disk_t=S^1_\varphi$.

We shall regard $\varphi$ and $\theta$ as functions on $S^3$.
As a rule, we use the notation $s$, $s'$, $s_1$, $s_2,\dots$
for points from $S^1_\varphi$ and $t$, $t'$, $t_1$, $t_2,\dots$
for points from $S^1_\theta$.

Let $L$ be a link in $S^3$.
By an \emph{arc presentation} of $L$ we shall mean a link $L'$
in $S^3$ isotopic to $L$ such that the
set $L'\cap S^1_\varphi$, whose elements
are called \emph{vertices}, is finite,
and, for any $t\in S^1$, the intersection of $L'$ with any page
$\disk_t$ is either the empty set or an open
arc approaching two distinct vertices.

\begin{rem}
The term `arc-presentation' was introduced by J.\,Birman and
W.\,Menasco, though the object itself appeared already in~\cite{brunn}.
However, it has received almost no attention before P.\,Cromwell's
paper~\cite{cromwell}.
\end{rem}

\subsection{Duality of arc-presentations}

Two arc-presentations, say $L_1$ and $L_2$, are not distinguished
if they have the same set of vertices, and, for any $t\in S^1_\theta$,
we have either $L_1\cap\disk_t=L_2\cap\disk_t=\varnothing$ or
$L_1\cap\disk_t$ and $L_2\cap\disk_t$ are arcs with the same endpoints.
Having agreed about this, we may assume without loss of
generality that any arc $L\cap\disk_t$ of an arc-presentation $L$
consists of two radii of the disk $\disk_t$. This means that
$L$ is a link consisting of segments of the form
$P*Q\subset S^1_\varphi*S^1_\theta$, where $P\in S^1_\varphi$,
$Q\in S^1_\theta$. Notice that, in this definition, the r\^oles of
the circles $S^1_\varphi$ and $S^1_\theta$ are the same,
and $L$ can be considered as an arc-presentation with respect
to both circles. Vertices of one of the presentations
are the centers of arcs of the other.

So, there is a duality operation on the set of arc-presentations,
which is defined by the mapping $\xi:S^3\rightarrow S^3$
written in our coordinate system as $\xi(\varphi,\tau,\theta)=
(\theta,1-\tau,\varphi)$.
The r\^oles of arcs and vertices are interchanged under this duality.

\begin{rem}
There is a triangulation $T(L)$ of $S^3$ associated naturally with
any arc-presentation $L$ having $n\geqslant2$ vertices.
$T(L)$ has $n^2$ simplices that have the form $I'*I''$, where
$I'\subset S^1_\varphi$, (respectively, $I''\subset S^1_\theta$)
is an interval between neighbouring vertices of the
arc-presentation $L$ (respectively, $\xi(L)$). In
this construction, $L$ is a subset of the 1-skeleton of $T(L)$
containing all the vertices. The main result of this
paper can be proved by using
normal surfaces
if one chooses $T(L)$ as the main triangulation. In this case,
the reasoning becomes ``self-dual'' with respect to $\xi$.
However, it seems to us that Birman--Menasco's approach,
which ``breaks the symmetry'', is easier to use in our case.
\end{rem}

\subsection{Arc-presentations and rectangular diagrams}
\label{arcandrect}

\begin{prop}
Let $L$ be an arc-presentation of a link, $D$ a rectangular
diagram\/ {\rm(}see the Introduction\/{\rm)} with vertices at all points
$(x,y)\in[0,2\pi)\times[0,2\pi)$
such that $L\cap\disk_y$ is an arc one of whose endpoints is
$(x,0,0)\in S^1_\varphi$. Then $D$ presents the same link as $L$ does.

The correspondence $L\mapsto D$ between arc-presentations
and rectangular diagrams in $[0,2\pi)\times[0,2\pi)$ is
one-to-one.
\end{prop}

\begin{proof}
The definition of $D$ makes sense, since each arc has
two endpoints, and two arcs are attached to each vertex.
Thus, in any vertical or horizontal straight line in the $(x,y)$-plane,
we will have either none or two vertices. In the latter
case, we connect them by a straight line segment,
obtaining the rectangular diagram $D$.

It is easy to see that one can obtain $D$ from $L$ by
cutting $S^3$ along disks $\disk_0$ and $\xi(\disk_0)$ and
then projecting the result appropriately
to the plane $\tau=1/2$.
See more explanations in~\cite{cromwell}.

The last assertion of the proposition is obvious.
\end{proof}

We shall regard rectangular diagrams as a convenient
way for depicting arc-presentations. Vertical edges
of a rectangular diagram correspond to vertices, horizontal
edges to arcs of an arc-presentation. We shall often
draw appropriate rectangular diagrams for illustrating
properties or transforms of arc-presentations.

\subsection{Arc-presentations and closed braids}\label{arcbraid}

There is also an easy way to convert a rectangular
diagram (an arc-presentation) to a closed braid.
An arc presentation or a rectangular diagram can be
endowed with an orientation as an ordinary link
or an ordinary planar link diagram, respectively.
For an oriented rectangular
diagram $D$, we shall call a horizontal edge of $D$
\emph{negative} if its orientation is opposite
to the orientation of the $x$-axis.

\begin{prop}\label{BD}
Let $D$ be an oriented rectangular diagram in $[0,2\pi)\times[0,2\pi)$.
Let us replace each negative horizontal edge in $D$ of the form
$[x_1,x_2]\times y$, $x_1<x_2$, by the union of two straight line segments
$([-\varepsilon,x_1]\cup[x_2,2\pi])\times y$, where $\varepsilon>0$.
As before, at each crossing, we interpret the vertical arc
as overcrossing and the horizontal one as undercrossing.
After a small deformation of the obtained picture,
we can get a planar diagram of a braid, say $b_D$.
Then the closure of $b_D$ is a link equivalent to the one
defined by $D$.

Any diagram of a braid can be obtained from some rectangular
diagram by the method just described.
\end{prop}

\begin{figure}[ht]\caption{Converting a rectangular diagram
to a braid}\label{braid}\vskip20pt
\centerline{\begin{picture}(365,75)
\put(0,75){\circle{3}}\put(2,75){\line(1,0){56}}\put(2,75){\vector(1,0){30}}
\put(60,75){\circle{3}}\put(60,73){\line(0,-1){26}}\put(60,45){\circle{3}}
\put(58,45){\line(-1,0){10}}\put(42,45){\line(-1,0){25}}
\put(15,45){\circle{3}}\put(15,43){\line(0,-1){26}}\put(15,15){\circle{3}}
\put(17,15){\line(1,0){10}}\put(33,15){\line(1,0){10}}\put(45,15){\circle{3}}
\put(45,17){\line(0,1){41}}\put(45,60){\circle{3}}\put(47,60){\line(1,0){10}}
\put(63,60){\line(1,0){10}}\put(75,60){\circle{3}}
\put(75,58){\line(0,-1){56}}\put(75,0){\circle{3}}\put(73,0){\line(-1,0){41}}
\put(30,0){\circle{3}}\put(30,2){\line(0,1){26}}\put(30,30){\circle{3}}
\put(28,30){\line(-1,0){10}}\put(12,30){\line(-1,0){10}}
\put(0,30){\circle{3}}\put(0,32){\line(0,1){41}}
\put(95,36){$\rightarrow$}
\put(135,75){\circle{3}}\put(137,75){\line(1,0){56}}\put(195,75){\circle{3}}
\put(195,73){\line(0,-1){26}}\put(195,45){\circle{3}}
\put(197,45){\line(1,0){10}}\put(213,45){\line(1,0){7}}
\put(148,45){\line(-1,0){10}}\put(132,45){\line(-1,0){7}}
\put(150,45){\circle{3}}\put(150,43){\line(0,-1){26}}\put(150,15){\circle{3}}
\put(152,15){\line(1,0){10}}\put(168,15){\line(1,0){10}}
\put(180,15){\circle{3}}\put(180,17){\line(0,1){41}}\put(180,60){\circle{3}}
\put(182,60){\line(1,0){10}}\put(198,60){\line(1,0){10}}
\put(210,60){\circle{3}}\put(210,58){\line(0,-1){56}}\put(210,0){\circle{3}}
\put(212,0){\line(1,0){8}}\put(163,0){\line(-1,0){38}}\put(165,0){\circle{3}}
\put(165,2){\line(0,1){26}}\put(165,30){\circle{3}}
\put(167,30){\line(1,0){10}}\put(183,30){\line(1,0){24}}
\put(213,30){\line(1,0){7}}\put(133,30){\line(-1,0){8}}
\put(135,30){\circle{3}}\put(135,32){\line(0,1){41}}
\put(240,36){$\rightarrow$}
\bezier200(270,30)(279,30)(280,45)
\bezier200(280,45)(282,75)(300,75)
\put(300,75){\line(1,0){30}}
\bezier200(330,75)(339,75)(340,60)
\bezier200(340,60)(341,50)(349,50)
\put(355,50){\line(1,0){10}}
\put(270,0){\line(1,0){30}}
\bezier200(300,0)(309,0)(310,15)
\bezier200(310,15)(311,30)(322,30)
\put(328,30){\line(1,0){24}}
\put(358,30){\line(1,0){7}}
\put(270,45){\line(1,0){7}}
\bezier200(283,45)(294,45)(295,30)
\bezier200(295,30)(296,15)(307,15)
\bezier200(313,15)(324,15)(325,30)
\bezier200(325,30)(327,60)(337,60)
\bezier200(343,60)(353,60)(355,30)
\bezier200(355,30)(357,0)(365,0)
\end{picture}}
\end{figure}
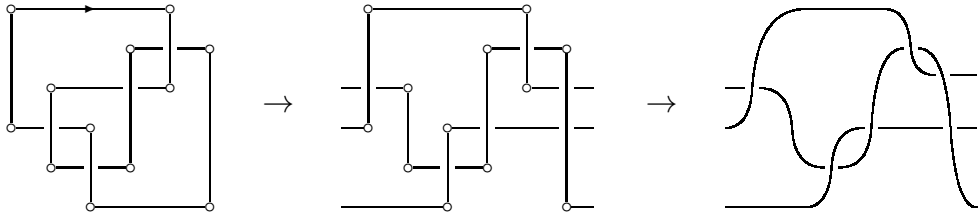

The idea is illustrated in Fig.~\ref{braid}. We skip the easy proof.
See also~\cite{cromwell}.

\subsection{Elementary moves}

Actually, elementary moves of arc-presentations have already been
introduced in the language of rectangular diagrams (see the Introduction).
Here we describe them directly.

Throughout this section, we assume that
each arc of any arc-presentation $L$
consists of two radii of a disk $\disk_t$. Thus, we have
$t\in S^1_\theta\cap L$ if and only if the page $\disk_t$ contains an
arc of $L$.

Suppose that we have a continuous
family $L(u)$ of arc presentations, {\it i.e.}, such
that the sets $L(u)\cap S^1_\varphi$ and $L(u)\cap S^1_\theta$
depend on $u$ continuously. Then all the links $L(u)$ are isotopic
to each other via isotopies of the form $(\varphi,\tau,\theta)\mapsto
(f(\varphi),\tau,g(\theta))$. Let us see what happens
to the corresponding rectangular diagrams $D(u)$ when $u$ changes.
If no point in $L(u)\cap S^1_\varphi$ and $L(u)\cap S^1_\theta$
passes through the origin $0\in S^1$, the diagram $D(u)$
also changes continuously, and the combinatorial type
of $D(u)$ stays unchanged. If one of the points in $L(u)\cap S^1_\varphi$
(respectively, $L(u)\cap S^1_\theta$) passes through the origin,
the corresponding rectangular diagram is changed by
a cyclic permutation of vertical (respectively, horizontal) edges.
Two arc-presentations that can be included into such a continuous
family are said to be \emph{combinatorially equivalent}.

The following assertion, although quite obvious, is very important
from the algorithmic point of view.

\begin{prop}
For any $n$, there are only finitely many classes
of combinatorially equivalent arc-presentations
of complexity $n$.
\end{prop}

\begin{rem}
One can show that the number $N(n)$ of pairwise distinct classes
of combinatorially equivalent arc-presentations
of complexity $n$ satisfies the inequalities
$((n-1)!)^2/(2n)<N(n)\leqslant((n-1)!)^2$.
\end{rem}

For three distinct points $x_1,x_2,x_3$ of $S^1$,
we shall write $x_2\in(x_1,x_3)$
if, when going in the positive direction of $S^1$
we meet the points $x_i$ in this order:
$\ldots,x_1,x_2,x_3,x_1$, $x_2,\ldots$.
We shall write $x_2\in[x_1,x_3]$ if we have either $x_2\in(x_1,x_3)$
or $x_2\in\{x_1,x_3\}$.
By $\gamma_{t,s,s'}$ we shall denote an arc in $\disk_t$
such that $\partial\gamma_{t,s,s'}=\{s,s'\}\subset S^1_\varphi$.
When using this notation,
we shall not distinguish between two arcs lying in the same
page and having the same endpoints. Thus, we write $\gamma_{t,s,s'}
\subset L$ whenever we want to say that $L\cap\disk_t$
contains an arc with endpoints $s,s'$.

Let $\alpha=\gamma_{t,s_1,s_2}$ be an arc of an arc-presentation $L$.
For sufficiently small
$\varepsilon>0$, there is no arc of $L$ in $\disk_{t'}$
if $t'\in(t-\varepsilon,t+\varepsilon)$, except for $\alpha$,
and there is no vertex of $L$ in
$(s_1-\varepsilon,s_1+\varepsilon)\subset S^1_\varphi$,
except for $s_1$. Replacing the arc $\alpha$ by the
arcs $\gamma_{t,s_1+\varepsilon_1,s_2}$ and $\gamma_{t+\varepsilon_2,
s_1,s_1+\varepsilon_1}$, where
$\varepsilon_{1,2}\in\{\varepsilon,-\varepsilon\}$,
will be called a \emph{stabilization move}.
The inverse operation will be called a \emph{destabilization move}.
In the language of rectangular diagrams, these operations
coincide with stabilization and destabilization moves of
rectangular diagrams, respectively,
provided that we have $0\notin(s_1-\varepsilon,s_1+\varepsilon)$
and $0\notin(t-\varepsilon,t+\varepsilon)$, which can
always be achieved by a small perturbation of $L$ and taking $\varepsilon$
sufficiently small.

Notice that the (de)stabilization move is self-dual with respect
to $\xi$.

Let $\alpha_1,\alpha_2$ be two arcs of
an arc-presentation $L$, $\alpha_i\subset\disk_{t_i}$,
and there is no arc of $L$ in
$\cup_{t\in(t_1,t_2)}\disk_t$ or in $\cup_{t\in(t_2,t_1)}\disk_t$.
In this case, we shall say that the arcs $\alpha_1$ and $\alpha_2$
are \emph{neighbouring}.

\begin{figure}[ht]\caption{Interleaved and non-interleaved arcs}
\label{interleave}
$$\begin{array}{cc}
\begin{picture}(80,80)
\put(40,40){\circle{40}}\put(60,40){\circle*{3}}\put(20,40){\circle*{3}}
\put(40,20){\circle*{3}}\put(40,60){\circle*{3}}\put(20,40){\line(1,0){40}}
\put(40,20){\line(0,1){40}}\put(7,40){$s_1$}\put(65,40){$s_2$}
\put(37,67){$s_3$}\put(37,10){$s_4$}
\end{picture}
&\begin{picture}(80,80)
\put(40,40){\circle{40}}\put(60,40){\circle*{3}}\put(20,40){\circle*{3}}
\put(40,20){\circle*{3}}\put(40,60){\circle*{3}}\put(20,40){\line(1,1){20}}
\put(40,20){\line(1,1){20}}\put(7,40){$s_1$}\put(65,40){$s_3$}
\put(37,67){$s_2$}\put(37,10){$s_4$}
\end{picture}\\
\mbox{interleaved arcs}&\mbox{non-interleaved arcs}
\end{array}$$
\end{figure}
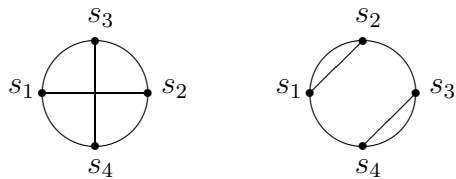

Let $\alpha_1$ and $\alpha_2$ be neighbouring arcs,
$s_1,s_2\in S^1_\varphi$ the endpoints
of $\alpha_1$, and $s_3,s_4$ be the endpoints of
$\alpha_2$, $\partial\alpha_1,\cap\partial\alpha_2=\varnothing$.
We say that the arcs $\alpha_1,\alpha_2$
are \emph{non-interleaved} if we have either
$s_1,s_2\in(s_3,s_4)$ or $s_1,s_2\in(s_4,s_3)$.
Otherwise, they are said to be \emph{interleaved}.
In other words, arcs are interleaved if the corresponding
chords of $S^1$ intersect each other (see Fig.~\ref{interleave}).

If $\alpha_1\subset\disk_{t_1}$,
$\alpha_2\subset\disk_{t_2}$ are non-interleaved neighbouring
arcs, we define an \emph{arc exchange move} on $L$
as replacing $\alpha_1,\alpha_2$ by arcs $\alpha_1'\in\disk_{t_2}$,
$\alpha_2'\in\disk_{t_1}$ such that $\partial\alpha_i'=
\partial\alpha_i$. This operation corresponds to
interchanging horizontal edges of a rectangular diagram, probably
combined with a cyclic permutation of horizontal edges.

We define a \emph{vertex exchange move} to be
the operation dual via $\xi$ to
an arc exchange move. It corresponds to
interchanging vertical edges of a rectangular diagram,
probably combined with a cyclic permutation of vertical edges.
We shall use the term `exchange move' for both operations,
an arc exchange move and a vertex exchange move.

In the sequel, stabilization, destabilization, and exchange moves will be
referred to as \emph{elementary} moves.
When appropriate, we shall not distinguish between combinatorially
equivalent arc-presentations and think about elementary
operations as being performed on the corresponding
combinatorial classes.

\begin{prop}\label{arc-pres}
Any link in $S^3$ has an arc-presentation.

Two arc-presentations $L_1,L_2$ of the same link
can be obtained from each other by finitely many
elementary moves.
\end{prop}

\begin{figure}[ht]\caption{A realization of the
Reidemeister moves}\label{reid}\vskip10pt
\centerline{\begin{picture}(80,60)
\put(0,10){\line(0,1){20}}\put(0,30){\line(1,0){20}}
\put(24,28){$\leftrightarrow$}
\put(60,10){\line(0,1){40}}\put(60,50){\line(-1,0){20}}
\put(40,50){\line(0,-1){20}}\put(40,30){\line(1,0){17}}
\put(63,30){\line(1,0){17}}
\end{picture}
\qquad
\begin{picture}(60,60)
\put(0,20){\line(1,0){10}}\put(10,20){\line(0,1){20}}
\put(10,40){\line(-1,0){10}}\put(20,10){\line(0,1){40}}
\put(24,28){$\leftrightarrow$}\put(40,20){\line(1,0){7}}
\put(53,20){\line(1,0){7}}\put(60,20){\line(0,1){20}}
\put(40,40){\line(1,0){7}}\put(53,40){\line(1,0){7}}
\put(50,10){\line(0,1){40}}
\end{picture}
\qquad
\begin{picture}(110,60)
\put(15,60){\line(0,-1){60}}\put(0,45){\line(1,0){12}}
\put(18,45){\line(1,0){12}}\put(30,45){\line(0,-1){30}}
\put(30,15){\line(1,0){15}}\put(0,30){\line(1,0){12}}
\put(18,30){\line(1,0){9}}\put(33,30){\line(1,0){12}}
\put(49,28){$\leftrightarrow$}\put(95,60){\line(0,-1){60}}
\put(65,45){\line(1,0){15}}\put(80,45){\line(0,-1){30}}
\put(80,15){\line(1,0){12}}\put(98,15){\line(1,0){12}}
\put(65,30){\line(1,0){12}}\put(83,30){\line(1,0){9}}
\put(98,30){\line(1,0){12}}
\end{picture}}
\end{figure}
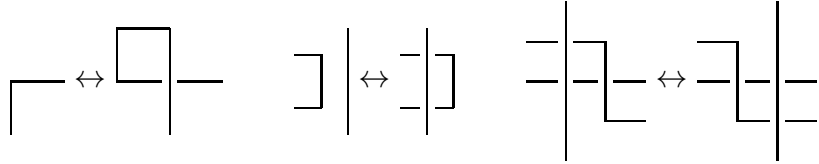

\begin{proof}
We provide here only a sketch, the details are very
elementary and standard. (In~\cite{cromwell}, one can find a proof
that uses Markov's theorem. The latter seems to us to be
much less elementary than the assertion to be proved.
We would prefer to use Reidemeister's theorem.)

As we have already mentioned, any planar diagram of a link
is isotopic to a rectangular diagram. Indeed,
for a given planar diagram, one can first deform it near crossings
so that overcrossing arcs become vertical and undercrossing arcs
horizontal. Then one can approximate the rest of the diagram
by a step-line consisting of vertical and horizontal straight
line segments. In the generic case, there will be no
collinear segments. In order to deform such a diagram into an isotopic one,
we might need to add more edges or contract some, which is achieved
by (de)stabilization moves, and pass through codimension one
degenerations, when two non-overlapping edges become collinear.
The latter is achieved by exchange moves.

In order to prove the second assertion, one shows first
that rotating a crossing by $\pi$ as shown in Fig.~\ref{rotate}
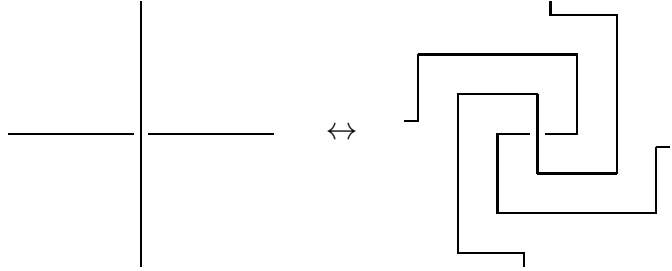
\begin{figure}[ht]\caption{Rotating a crossing}\label{rotate}\vskip20pt
\centerline{\begin{picture}(250,100)
\put(50,0){\line(0,1){100}}\put(0,50){\line(1,0){47}}
\put(100,50){\line(-1,0){47}}\put(120,48){$\leftrightarrow$}
\put(200,35){\line(0,1){30}}\put(185,50){\line(1,0){12}}
\put(215,50){\line(-1,0){12}}\put(200,35){\line(1,0){30}}
\put(230,35){\line(0,1){60}}\put(230,95){\line(-1,0){25}}
\put(205,95){\line(0,1){5}}\put(200,65){\line(-1,0){30}}
\put(170,65){\line(0,-1){60}}\put(170,5){\line(1,0){25}}
\put(195,5){\line(0,-1){5}}\put(185,50){\line(0,-1){30}}
\put(185,20){\line(1,0){60}}\put(245,20){\line(0,1){25}}
\put(245,45){\line(1,0){5}}\put(215,50){\line(0,1){30}}
\put(215,80){\line(-1,0){60}}\put(155,80){\line(0,-1){25}}
\put(155,55){\line(-1,0){5}}
\end{picture}}
\end{figure}
can be performed by applying finitely many elementary moves.
Then one do the same with Reidemeister moves, assuming
that the involved crossings are already positioned as one wishes
(see Fig.~\ref{reid}).\end{proof}

\begin{rem} In~\cite{cromwell}, the set of allowed (de)stabilization-type
moves in the formulations
of an analogue of Proposition~\ref{arc-pres} is larger. We shall now
see that those, more general, moves can be expressed in terms of
our elementary moves.
\end{rem}

\subsection{Generalized moves}\label{generalized}
Let $\alpha\in\disk_t$ be an arc of an arc-presentation $L$ and
$s\in S^1_\varphi$ be a point distinct from vertices of $L$.
Let $s_1,s_2\in S^1_\varphi$ be the endpoints of $\alpha$.
Replacing $\alpha$ by the arcs
$\gamma_{t,s_1,s}$ and $\gamma_{t+\varepsilon,s,s_2}$,
where $\varepsilon>0$ is sufficiently small, will be
called a \emph{generalized stabilization move}.
(In~\cite{cromwell}, it is called a type~IV move. The dual
operation, which we shall not need, is called a type~III move
in~\cite{cromwell}.) The operation inverse to a generalized
stabilization move is called a \emph{generalized destabilization move}.

Let $L$ be an arc-presentation, $s_1,s_2,s_3$ three points
of $S^1_\varphi$ distinct from vertices of $L$ such that $s_2\in(s_1,s_3)$.
Suppose that for some $t_1,t_2\in S^1_\theta$,
the link $L$ does not intersect any of the arcs
$\gamma_{t,s_1,s_2}$ if $t\in[t_1,t_2]$ and
$\gamma_{t,s_2,s_3}$ if $t\in[t_2,t_1]$.
In this case, we shall call the operation consisting in
interchanging the intervals $(s_1,s_2)$ $(s_2,s_3)$
a \emph{generalized} (\emph{vertex}) \emph{exchange move}.
This operation can be described more formally as follows.
Define a function $f:S^1_\varphi\rightarrow S^1_\varphi$ by
$$s\mapsto\left\{\begin{array}{ll}
s+s_3-s_2&\mbox{if }s\in(s_1,s_2),\\
s+s_1-s_2&\mbox{if }s\in(s_2,s_3),\\
s&\mbox{otherwise}.
\end{array}\right.$$
Then the generalized exchange move just defined consists in
replacing each arc of the form $\gamma_{t,s,s'}$ by the arc
$\gamma_{t,f(s),f(s')}$. (In~\cite{cromwell}, this operation
is called generalized type~I move.)
The dual operation, generalized arc exchange move, will
not be needed here.

\begin{prop}\label{genexchange}
Any generalized\/ {\rm(}de\/{\rm)}stabilization move can be presented
as the composition of one ordinary {\rm(}de\/{\rm)}stabilization move
and a few exchange moves.
Any generalized exchange move can be presented as the composition of
ordinary exchange moves.
\end{prop}

\begin{proof}
We illustrate the idea of the proof in Figures~\ref{genstab}
and \ref{genexch} by using corresponding rectangular
diagrams. The details are easy.
\end{proof}

\begin{figure}[ht]\caption{Decomposition of a generalized
stabilization move}\label{genstab}
\vskip20pt
\centerline{\begin{picture}(370,50)
\put(20,20){\line(0,1){8}}
\put(20,30){\circle{3}}
\put(22,30){\line(1,0){66}}
\put(90,30){\circle{3}}
\put(90,32){\line(0,1){18}}
\put(0,30){${}_t$}\multiput(4,30)(4,0){4}{\line(1,0){2}}
\put(18,4){${}_{s_1}$}
\multiput(20,8)(0,4){3}{\line(0,1){2}}
\put(88,4){${}_{s_2}$}
\multiput(90,8)(0,4){5}{\line(0,1){2}}
\put(40,20){\line(0,1){30}}
\put(50,20){\line(0,1){30}}
\put(70,20){\line(0,1){30}}
\put(160,20){\line(0,1){8}}
\put(160,30){\circle{3}}
\put(162,30){\line(1,0){6}}
\put(170,30){\circle{3}}
\put(170,32){\line(0,1){6}}
\put(170,40){\circle{3}}
\put(172,40){\line(1,0){56}}
\put(230,40){\circle{3}}
\put(230,42){\line(0,1){8}}
\put(140,30){${}_t$}
\multiput(144,30)(4,0){4}{\line(1,0){2}}
\put(140,40){${}_{t+\varepsilon}$}
\multiput(154,40)(4,0){4}{\line(1,0){2}}
\put(158,4){${}_{s_1}$}
\multiput(160,8)(0,4){3}{\line(0,1){2}}
\put(162,10){${}_{s_1+\varepsilon}$}
\multiput(170,14)(0,4){4}{\line(0,1){2}}
\put(228,4){${}_{s_2}$}
\multiput(230,8)(0,4){8}{\line(0,1){2}}
\put(180,20){\line(0,1){30}}
\put(190,20){\line(0,1){30}}
\put(210,20){\line(0,1){30}}
\put(300,20){\line(0,1){8}}
\put(300,30){\circle{3}}
\put(302,30){\line(1,0){36}}
\put(340,30){\circle{3}}
\put(340,32){\line(0,1){6}}
\put(340,40){\circle{3}}
\put(342,40){\line(1,0){26}}
\put(370,40){\circle{3}}
\put(370,42){\line(0,1){8}}
\put(280,30){${}_t$}
\multiput(284,30)(4,0){4}{\line(1,0){2}}
\put(280,40){${}_{t+\varepsilon}$}
\multiput(294,40)(4,0){11}{\line(1,0){2}}
\put(298,4){${}_{s_1}$}
\multiput(300,8)(0,4){3}{\line(0,1){2}}
\put(338,4){${}_{s}$}
\multiput(340,8)(0,4){5}{\line(0,1){2}}
\put(368,4){${}_{s_2}$}
\multiput(370,8)(0,4){8}{\line(0,1){2}}
\put(320,20){\line(0,1){30}}
\put(330,20){\line(0,1){30}}
\put(350,20){\line(0,1){30}}
\put(95,30){\hbox to 40pt{$-$\hss$-$\hss$-$\hss$-$\hss$\rightarrow$}}
\put(100,36){\scriptsize stabili-}
\put(100,26){\scriptsize zation}
\put(235,30){\hbox to 40pt{$-$\hss$-$\hss$-$\hss$-$\hss$\rightarrow$}}
\put(242,36){\scriptsize vertex}
\put(237,26){\scriptsize exchange}
\end{picture}}
\end{figure}
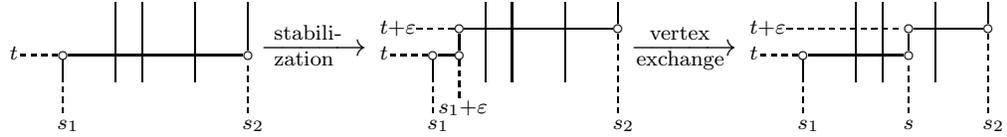
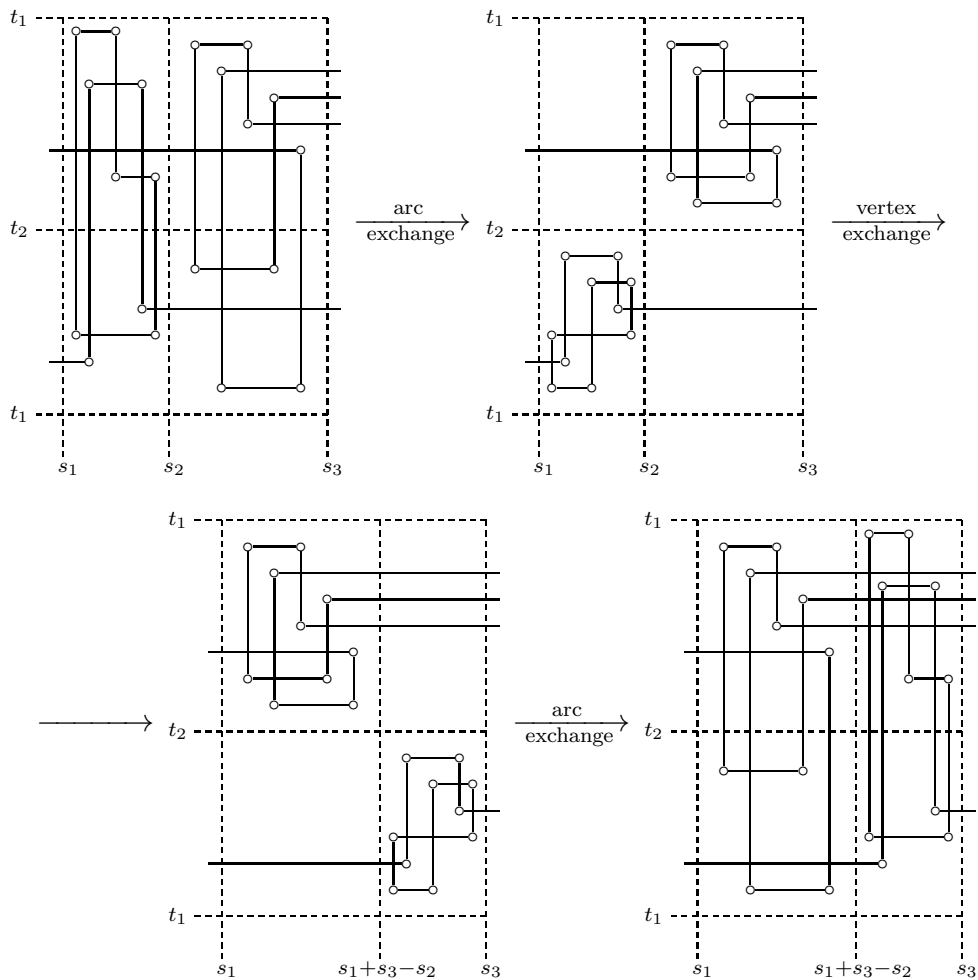
\begin{figure}[ht]\caption{Decomposition of a generalized
exchange move}\label{genexch}
\vskip20pt
\centerline{\begin{picture}(360,360)
\put(0,190){\begin{picture}(120,170)
\put(0,20){${}_{t_1}$}
\multiput(10,20)(4,0){28}{\line(1,0){2}}
\put(0,90){${}_{t_2}$}
\multiput(10,90)(4,0){28}{\line(1,0){2}}
\put(0,170){${}_{t_1}$}
\multiput(10,170)(4,0){28}{\line(1,0){2}}
\put(18,0){${}_{s_1}$}
\multiput(20,4)(0,4){42}{\line(0,1){2}}
\put(58,0){${}_{s_2}$}
\multiput(60,4)(0,4){42}{\line(0,1){2}}
\put(118,0){${}_{s_3}$}
\multiput(120,4)(0,4){42}{\line(0,1){2}}
\put(15,40){\line(1,0){13}}
\put(30,40){\circle{3}}
\put(30,42){\line(0,1){101}}
\put(30,145){\circle{3}}
\put(32,145){\line(1,0){16}}
\put(50,145){\circle{3}}
\put(50,143){\line(0,-1){81}}
\put(50,60){\circle{3}}
\put(52,60){\line(1,0){73}}
\put(25,165){\circle{3}}
\put(25,163){\line(0,-1){111}}
\put(25,50){\circle{3}}
\put(27,50){\line(1,0){26}}
\put(55,50){\circle{3}}
\put(55,52){\line(0,1){56}}
\put(55,110){\circle{3}}
\put(53,110){\line(-1,0){11}}
\put(40,110){\circle{3}}
\put(40,112){\line(0,1){51}}
\put(40,165){\circle{3}}
\put(38,165){\line(-1,0){11}}
\put(15,120){\line(1,0){93}}
\put(110,120){\circle{3}}
\put(110,118){\line(0,-1){86}}
\put(110,30){\circle{3}}
\put(108,30){\line(-1,0){26}}
\put(80,30){\circle{3}}
\put(80,32){\line(0,1){116}}
\put(80,150){\circle{3}}
\put(82,150){\line(1,0){43}}
\put(125,140){\line(-1,0){23}}
\put(100,140){\circle{3}}
\put(100,138){\line(0,-1){61}}
\put(100,75){\circle{3}}
\put(98,75){\line(-1,0){26}}
\put(70,75){\circle{3}}
\put(70,77){\line(0,1){81}}
\put(70,160){\circle{3}}
\put(72,160){\line(1,0){16}}
\put(90,160){\circle{3}}
\put(90,158){\line(0,-1){26}}
\put(90,130){\circle{3}}
\put(92,130){\line(1,0){33}}
\end{picture}}
\put(130,280){\hbox to45pt{$-$\hss$-$\hss$-$\hss$-$\hss$-$\hss$\rightarrow$}}
\put(145,286){\scriptsize arc}
\put(135,276){\scriptsize exchange}
\put(180,190){\begin{picture}(120,170)
\put(0,20){${}_{t_1}$}
\multiput(10,20)(4,0){28}{\line(1,0){2}}
\put(0,90){${}_{t_2}$}
\multiput(10,90)(4,0){28}{\line(1,0){2}}
\put(0,170){${}_{t_1}$}
\multiput(10,170)(4,0){28}{\line(1,0){2}}
\put(18,0){${}_{s_1}$}
\multiput(20,4)(0,4){42}{\line(0,1){2}}
\put(58,0){${}_{s_2}$}
\multiput(60,4)(0,4){42}{\line(0,1){2}}
\put(118,0){${}_{s_3}$}
\multiput(120,4)(0,4){42}{\line(0,1){2}}
\put(15,40){\line(1,0){13}}
\put(30,40){\circle{3}}
\put(30,42){\line(0,1){36}}
\put(30,80){\circle{3}}
\put(32,80){\line(1,0){16}}
\put(50,80){\circle{3}}
\put(50,78){\line(0,-1){16}}
\put(50,60){\circle{3}}
\put(52,60){\line(1,0){73}}
\put(25,30){\circle{3}}
\put(25,32){\line(0,1){16}}
\put(25,50){\circle{3}}
\put(27,50){\line(1,0){26}}
\put(55,50){\circle{3}}
\put(55,52){\line(0,1){16}}
\put(55,70){\circle{3}}
\put(53,70){\line(-1,0){11}}
\put(40,70){\circle{3}}
\put(40,68){\line(0,-1){36}}
\put(40,30){\circle{3}}
\put(38,30){\line(-1,0){11}}
\put(15,120){\line(1,0){93}}
\put(110,120){\circle{3}}
\put(110,118){\line(0,-1){16}}
\put(110,100){\circle{3}}
\put(108,100){\line(-1,0){26}}
\put(80,100){\circle{3}}
\put(80,102){\line(0,1){46}}
\put(80,150){\circle{3}}
\put(82,150){\line(1,0){43}}
\put(125,140){\line(-1,0){23}}
\put(100,140){\circle{3}}
\put(100,138){\line(0,-1){26}}
\put(100,110){\circle{3}}
\put(98,110){\line(-1,0){26}}
\put(70,110){\circle{3}}
\put(70,112){\line(0,1){46}}
\put(70,160){\circle{3}}
\put(72,160){\line(1,0){16}}
\put(90,160){\circle{3}}
\put(90,158){\line(0,-1){26}}
\put(90,130){\circle{3}}
\put(92,130){\line(1,0){33}}
\end{picture}}
\put(310,280){\hbox to45pt{$-$\hss$-$\hss$-$\hss$-$\hss$-$\hss$\rightarrow$}}
\put(320,286){\scriptsize vertex}
\put(315,276){\scriptsize exchange}
\put(10,90){\hbox to45pt{$-$\hss$-$\hss$-$\hss$-$\hss$-$\hss$\rightarrow$}}
\put(60,0){\begin{picture}(120,170)
\put(0,20){${}_{t_1}$}
\multiput(10,20)(4,0){28}{\line(1,0){2}}
\put(0,90){${}_{t_2}$}
\multiput(10,90)(4,0){28}{\line(1,0){2}}
\put(0,170){${}_{t_1}$}
\multiput(10,170)(4,0){28}{\line(1,0){2}}
\put(18,0){${}_{s_1}$}
\multiput(20,4)(0,4){42}{\line(0,1){2}}
\put(64,0){${}_{s_1+s_3-s_2}$}
\multiput(80,4)(0,4){42}{\line(0,1){2}}
\put(118,0){${}_{s_3}$}
\multiput(120,4)(0,4){42}{\line(0,1){2}}
\put(15,40){\line(1,0){73}}
\put(90,40){\circle{3}}
\put(90,42){\line(0,1){36}}
\put(90,80){\circle{3}}
\put(92,80){\line(1,0){16}}
\put(110,80){\circle{3}}
\put(110,78){\line(0,-1){16}}
\put(110,60){\circle{3}}
\put(112,60){\line(1,0){13}}
\put(85,30){\circle{3}}
\put(85,32){\line(0,1){16}}
\put(85,50){\circle{3}}
\put(87,50){\line(1,0){26}}
\put(115,50){\circle{3}}
\put(115,52){\line(0,1){16}}
\put(115,70){\circle{3}}
\put(113,70){\line(-1,0){11}}
\put(100,70){\circle{3}}
\put(100,68){\line(0,-1){36}}
\put(100,30){\circle{3}}
\put(98,30){\line(-1,0){11}}
\put(15,120){\line(1,0){53}}
\put(70,120){\circle{3}}
\put(70,118){\line(0,-1){16}}
\put(70,100){\circle{3}}
\put(68,100){\line(-1,0){26}}
\put(40,100){\circle{3}}
\put(40,102){\line(0,1){46}}
\put(40,150){\circle{3}}
\put(42,150){\line(1,0){83}}
\put(125,140){\line(-1,0){63}}
\put(60,140){\circle{3}}
\put(60,138){\line(0,-1){26}}
\put(60,110){\circle{3}}
\put(58,110){\line(-1,0){26}}
\put(30,110){\circle{3}}
\put(30,112){\line(0,1){46}}
\put(30,160){\circle{3}}
\put(32,160){\line(1,0){16}}
\put(50,160){\circle{3}}
\put(50,158){\line(0,-1){26}}
\put(50,130){\circle{3}}
\put(52,130){\line(1,0){73}}
\end{picture}}
\put(190,90){\hbox to45pt{$-$\hss$-$\hss$-$\hss$-$\hss$-$\hss$\rightarrow$}}
\put(205,96){\scriptsize arc}
\put(195,86){\scriptsize exchange}
\put(240,0){\begin{picture}(120,170)
\put(0,20){${}_{t_1}$}
\multiput(10,20)(4,0){28}{\line(1,0){2}}
\put(0,90){${}_{t_2}$}
\multiput(10,90)(4,0){28}{\line(1,0){2}}
\put(0,170){${}_{t_1}$}
\multiput(10,170)(4,0){28}{\line(1,0){2}}
\put(18,0){${}_{s_1}$}
\multiput(20,4)(0,4){42}{\line(0,1){2}}
\put(64,0){${}_{s_1+s_3-s_2}$}
\multiput(80,4)(0,4){42}{\line(0,1){2}}
\put(118,0){${}_{s_3}$}
\multiput(120,4)(0,4){42}{\line(0,1){2}}
\put(15,40){\line(1,0){73}}
\put(90,40){\circle{3}}
\put(90,42){\line(0,1){101}}
\put(90,145){\circle{3}}
\put(92,145){\line(1,0){16}}
\put(110,145){\circle{3}}
\put(110,143){\line(0,-1){81}}
\put(110,60){\circle{3}}
\put(112,60){\line(1,0){13}}
\put(85,165){\circle{3}}
\put(85,163){\line(0,-1){111}}
\put(85,50){\circle{3}}
\put(87,50){\line(1,0){26}}
\put(115,50){\circle{3}}
\put(115,52){\line(0,1){56}}
\put(115,110){\circle{3}}
\put(113,110){\line(-1,0){11}}
\put(100,110){\circle{3}}
\put(100,112){\line(0,1){51}}
\put(100,165){\circle{3}}
\put(98,165){\line(-1,0){11}}
\put(15,120){\line(1,0){53}}
\put(70,120){\circle{3}}
\put(70,118){\line(0,-1){86}}
\put(70,30){\circle{3}}
\put(68,30){\line(-1,0){26}}
\put(40,30){\circle{3}}
\put(40,32){\line(0,1){116}}
\put(40,150){\circle{3}}
\put(42,150){\line(1,0){83}}
\put(125,140){\line(-1,0){63}}
\put(60,140){\circle{3}}
\put(60,138){\line(0,-1){61}}
\put(60,75){\circle{3}}
\put(58,75){\line(-1,0){26}}
\put(30,75){\circle{3}}
\put(30,77){\line(0,1){81}}
\put(30,160){\circle{3}}
\put(32,160){\line(1,0){16}}
\put(50,160){\circle{3}}
\put(50,158){\line(0,-1){26}}
\put(50,130){\circle{3}}
\put(52,130){\line(1,0){73}}
\end{picture}}
\end{picture}}
\end{figure}

\section{Proof of the main result}

Recall, that a rectangular diagram $D$ is said to be \emph{trivial}
if it consists of four edges. It is \emph{split} if there
is a vertical line $l$ non-intersecting the diagram
such that $D$ has edges on both sides of $l$.

In this section, a rectangular diagram $D$ is said to be \emph{composite}
if there is a vertical line $l$ intersecting $D$
in two points such that $D$ has horizontal edges on both
sides of $l$. This is a particular case of the definition
given in the Introduction.

An arc-presentations whose rectangular diagram is trivial,
split, or composite is called \emph{trivial}, \emph{split}, or
\emph{composite}, respectively.

The assertion of Theorem~\ref{th1} will be a consequence of
the following three Propositions.

\begin{prop}\label{triv}
If $L$ is an arc-presentation of the unknot, then there
exists a finite sequence of exchange and destabilization moves
$L\mapsto L_1\mapsto L_2\mapsto\ldots\mapsto L_n$
such that the arc-presentation $L_n$ is trivial.
\end{prop}

\begin{prop}\label{split}
If $L$ is an arc-presentation of a split link, then there
exists a finite sequence of exchange and destabilization moves
$L\mapsto L_1\mapsto L_2\mapsto\ldots\mapsto L_n$
such that the arc-presentation $L_n$ is split.
\end{prop}

\begin{prop}\label{comp}
If $L$ is an arc-presentation of a non-split composite link, then there
exists a finite sequence of exchange and destabilization moves
$L\mapsto L_1\mapsto L_2\mapsto\ldots\mapsto L_n$
such that the arc-presentation $L_n$ is composite.
\end{prop}

The proof of these three Propositions follows the same scheme, only
technical details are different. So, we shall be proving
all three statements
simultaneously.
The strategy and methods of the proof are very close to those in
Birman--Menasco works~\cite{bm4},
\cite{bm5}. However, we need to concentrate on many technical
details, which are specific in our case. This
forces us to write a complete proof instead of just using
existing results.

To begin with, we describe the general outline
of the proof.

\begin{enumerate}
\item The assumptions of Propositions~\ref{triv}--\ref{comp}
imply the existence of a certain surface $M$: a disk whose boundary
is $L$, a two-sphere that does not intersect $L$ and is such that
there are non-trivial parts of $L$ on both sides of $M$,
or a two-sphere that intersects $L$ twice and cuts $L$ in two
 non-trivial tangles. We show that such a surface $M$,
which we call a \emph{characteristic surface},
can be isotoped to satisfy certain restrictions.
Compared to the previous works~\cite{bm4,bm5,cromwell},
a new thing here is only a restriction on the behaviour
of the spanning disk near the boundary in the case of the unknot.

\item We consider the foliation $\cal F$ on $M$ defined
by the equation $d\theta=0$. We define the \emph{complexity}
$c(M)\in\integer$ of the characteristic surface and show that, if
$c(M)>0$, then the foliation $\cal F$ contains
certain patterns. In the case when $M$ is a two-sphere
we just follow~\cite{bm4,bm5}, and extend the argument to the case of
the disk.

\item We show that, if $\cal F$ contains a pattern defined at the
previous step, then there exists an arc-presentation $L'$ and a surface
$M'$ such that $L'$ is obtained from $L$ by
finitely many exchange and destabilization moves, the pair $(L',M')$
is isotopic to $(L, M)$, the surface $M'$
satisfies the restrictions introduced at Step~1, and we have
either $c(L')<c(L)$ or $c(L')=c(L)$ and $c(M')<c(M)$.
Similarly to~\cite{bm4,bm5,cromwell}, this is done by using
two tricks that are closely related to those introduced by
D.\,Bennequin in~\cite{ben}.
In the unknot case we also need one more, easier,
trick that allows to simplify the boundary of a spanning disk.
\item We notice that, if $c(F)=0$, then we are done, {\it i.e.}, the
arc-presentation $L$ is trivial, split, or composite, respectively.
\end{enumerate}

\subsection{Characteristic surfaces}
Let $L$ be an arc-presentation satisfying the assumptions of one
of the Propositions~\ref{triv}--\ref{comp}. From now on, we shall assume
that each arc of the form $L\cap\disk_t$ is smooth, which
was not the case under the agreement of Section~\ref{prelim}.
Let $M\subset S^3$
be
\begin{itemize}
\item an embedded two-dimensional disk whose boundary is $L$ if $L$
is an arc-pre\-sen\-ta\-tion of the unknot;
\item an embedded two-sphere splitting $L$
into two non-empty links if $L$ is a split link;
\item a factorizing sphere, {\it i.e.},
an embedded sphere that meets $L$ in two points, cutting $L$
into two non-trivial tangles, if $L$ is a composite non-split link.
\end{itemize}
In this case, we call $M$ a \emph{characteristic surface} for $L$.

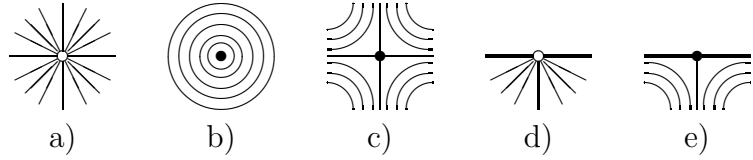
\begin{figure}[ht]\caption{Vertices and singularities of
$\cal F$}\label{singularities}
$$\begin{array}{ccccccccc}
\begin{picture}(40,40)
\put(20,20){\circle{4}}
\put(20,0){\line(0,1){18}}
\put(0,20){\line(1,0){18}}
\put(20,40){\line(0,-1){18}}
\put(40,20){\line(-1,0){18}}
\put(6,6){\line(1,1){12}}
\put(2,11){\line(2,1){16}}
\put(11,2){\line(1,2){8}}
\put(34,6){\line(-1,1){12}}
\put(38,11){\line(-2,1){16}}
\put(29,2){\line(-1,2){8}}
\put(6,34){\line(1,-1){12}}
\put(2,29){\line(2,-1){16}}
\put(11,38){\line(1,-2){8}}
\put(34,34){\line(-1,-1){12}}
\put(38,29){\line(-2,-1){16}}
\put(29,38){\line(-1,-2){8}}
\end{picture}&&
\begin{picture}(40,40)
\put(20,20){\circle*{4}}
\put(20,20){\circle{10}}
\put(20,20){\circle{16}}
\put(20,20){\circle{25}}
\put(20,20){\circle{30}}
\put(20,20){\circle{40}}
\end{picture}&&
\begin{picture}(40,40)
\put(20,20){\circle*{4}}
\put(20,0){\line(0,1){40}}
\put(0,20){\line(1,0){40}}
\put(0,0){\oval(20,20)[tr]}
\put(0,0){\oval(27,27)[tr]}
\put(0,0){\oval(35,35)[tr]}
\put(40,0){\oval(20,20)[tl]}
\put(40,0){\oval(27,27)[tl]}
\put(40,0){\oval(35,35)[tl]}
\put(0,40){\oval(20,20)[br]}
\put(0,40){\oval(27,27)[br]}
\put(0,40){\oval(35,35)[br]}
\put(40,40){\oval(20,20)[bl]}
\put(40,40){\oval(27,27)[bl]}
\put(40,40){\oval(35,35)[bl]}
\end{picture}&&
\begin{picture}(40,40)
\put(20,20){\circle{4}}
{\linethickness{1pt}\put(0,20){\line(1,0){18}}
\put(40,20){\line(-1,0){18}}}
\put(20,0){\line(0,1){18}}
\put(6,6){\line(1,1){12}}
\put(2,11){\line(2,1){16}}
\put(11,2){\line(1,2){8}}
\put(34,6){\line(-1,1){12}}
\put(38,11){\line(-2,1){16}}
\put(29,2){\line(-1,2){8}}
\end{picture}&&
\begin{picture}(40,40)
\put(20,20){\circle*{4}}
{\linethickness{1pt}\put(0,20){\line(1,0){18}}
\put(40,20){\line(-1,0){18}}}
\put(20,0){\line(0,1){18}}
\put(0,0){\oval(20,20)[tr]}
\put(0,0){\oval(27,27)[tr]}
\put(0,0){\oval(35,35)[tr]}
\put(40,0){\oval(20,20)[tl]}
\put(40,0){\oval(27,27)[tl]}
\put(40,0){\oval(35,35)[tl]}
\end{picture}
\\\mbox{a)}&&\mbox{b)}&&\mbox{c)}&&\mbox{d)}&&\mbox{e)}
\end{array}$$
\end{figure}

We say that a characteristic surface $M$ is \emph{admissible} if
\begin{itemize}
\item the surface $M$ is smooth everywhere up to boundary, except
at $(\partial M)\cap S^1_\varphi$;
\item $M\setminus\partial M$ intersects
the binding circle $S^1_\varphi$
transversely at finitely many points;
\item the foliation $\cal F$ on $M\setminus
S^1_\varphi$ defined by $d\theta=0$ has only finitely
many singularities, which are points of tangency
of $M$ with pages $\disk_t$;
\item all singularities of $\cal F$ are of Morse type,
{\it i.e.}, local extrema (Fig.~\ref{singularities}~b) or saddle
critical points (Fig.~\ref{singularities}~c) of
the (multivalued) function $\theta|_{M\setminus S^1_\varphi}$;
\item near any point $p\in(\partial M)\cap S^1_\varphi$,
the foliation $\cal F$ is radial (see Fig.~\ref{singularities}~d);
\item there is at most one point $p\in(\partial M)\cap S^1_\varphi$
at which $|\int_\gamma d\theta|>2\pi$, where $\gamma\subset M$
is a proper arc in a small neighbourhood of $p$ such that the endpoints
of $\gamma$ lie in $\partial M$ on different sides of $p$.
Such a point is called a \emph{winding vertex};
\item there is at most one point
$p\in(\partial M)\setminus S^1_\varphi$ at which
the surface $M$ is not transversal to the
corresponding page $\disk_{\theta(p)}$;
at the exceptional point, the foliation $\cal F$ must
have a saddle critical point (see Fig.~\ref{singularities}~e).
If such a saddle and a winding vertex are both present, then
the winding vertex is an endpoint of the edge containing the saddle;
\item each page $\disk_t$ contains no more than one of the following:
\begin{itemize}
\item an arc of $L$;
\item a singularity of ${\cal F}|_{M\setminus\partial M}$;
\end{itemize}
\item if $L$ is a non-split composite link, the two intersection
points $L\cap M$ are not vertices of $L$.
\end{itemize}

The behaviour of $\cal F$ near a point from $(M\setminus
\partial M)\cap S^1_\varphi$ is shown in Fig.~\ref{singularities}~a).
The points of intersection of $M$ with the binding circle,
where the foliation $\cal F$ is not defined, will be called
\emph{vertices} of $\cal F$. Sometimes we will distinguish
between \emph{interior} and \emph{boundary} vertices
that lie in $M\setminus\partial M$
and $\partial M$, respectively. The singularities
of $\cal F$ shown in Fig.~\ref{singularities}~b,c,e
will be called a \emph{pole}, a(n interior) \emph{saddle},
and a \emph{boundary
saddle}, respectively.

\begin{lemma}
For any arc-presentation satisfying the assumptions
of one of the Propositions\/~{\rm\ref{triv}--\ref{comp}},
there exists an admissible characteristic surface.
\end{lemma}

\begin{proof}
In the case of a split link or a composite link, there
is almost nothing to prove: a characteristic surface
is admissible if it is in general position with respect
to the binding circle, to the link, and to the foliation of
$S^3\setminus S^1_\varphi$ by disks $\disk_t$.

In the case of unknot, the admissibility of $M$
means something more than just general position.
By a small perturbation of an arbitrary characteristic
surface $M$, we can achieve that $M$
satisfy almost all the restrictions.
The only thing that we may not avoid in this way is that $\cal F$
have \emph{many} boundary saddles, whereas
we allow it to have just one.

Notice, that if the number of vertices of an arc-presentation of
the unknot is odd, then the foliation on the spanning disk is forced
to have a boundary saddle, since the parity of the number of
boundary saddles should be equal to that of the number of boundary vertices.
Indeed, the coorientation of $\partial M$ defined by
$\qopname\relax o{grad}\theta$
flips at each vertex and at each boundary saddle, so, in
total there should be an even number of them.

In order to avoid the occurrence of many boundary
saddles, we shall construct the disk $M$,
starting from a small neighbourhood of the boundary.
We construct initially a narrow ribbon $R\cong [0,1]\times S^1$ such that
$\partial R=L\cup L'$, where $L'$ is a circle unlinked
with $L$, and the foliation on $R$ defined by $d\theta=0$
has the desired behaviour near $L$. Then we attach a
two-dimensional disk along $L'$, obtaining a characteristic surface
$M$. Finally, we deform $M$ slightly if necessary,
keeping it fixed in a small neighbourhood of $L=\partial M$,
so that, after deformation, $M\setminus\partial M$ will
be in general position.
So, the only thing to explain is how to construct the ribbon $R$.

Notice that constructing the ribbon $R$
is \emph{not} equivalent to specifying
a trivial framing of $L$. It would have been so if $L$ were
a smooth curve.

We start by enumerating the vertices $s_1,s_2,\dots,s_n\in S^1_\varphi$
and arcs $\alpha_1,\dots,\alpha_n$ of $L$ in the order
that they follow in $L$:
$$\partial\alpha_i=\{s_i,s_{i+1}\},\quad\mbox{if}\quad i=1,\dots,n-1,
\quad\mbox{and}\quad\partial\alpha_n=\{s_n,s_1\}.$$
Denote by $r_u$ the rotation of $S^3$ about the binding circle
by the angle $u$: $r_u(\varphi,\tau,\theta)=(\varphi,\tau,\theta+u)$.
For an arc $\alpha\subset\disk_t$, we denote
by $d_\varepsilon(\alpha)$ the disk
$$\bigcup\limits_{u\in[0,\varepsilon]}r_u(\alpha),
\quad\mbox{if }\varepsilon>0$$
and
$$\bigcup\limits_{u\in[\varepsilon,0]}r_u(\alpha),
\quad\mbox{if }\varepsilon<0.$$

Take a point $p$ in the arc $\alpha_n$.
Let $\alpha_n',\alpha_n''$ be the two parts of the arc $\alpha$
cut by $p$: $\partial\alpha_n'=\{s_n,p\}$, $\partial\alpha_n''=
\{p,s_1\}$.
Pick an $\varepsilon>0$ which is smaller than
the $\theta$-distance between any two neighbouring arcs.
Consider the union $R_1$ of the following disks:
$$d_\varepsilon(\alpha_1),d_{-\varepsilon}(\alpha_2),
d_\varepsilon(\alpha_3),\dots,d_{(-1)^n\varepsilon}(\alpha_{n-1}),
d_{(-1)^{n+1}\varepsilon}(\alpha_n'),d_{-\varepsilon}(\alpha_n'').$$
If $n$ is even, the last two disks form the disk
$d_{-\varepsilon}(\alpha_n)$. If $n$ is odd, then the intersection
of $R_1$ with a neighbourhood of the point $p$ will look as
shown in Fig.~\ref{r1} on the left.
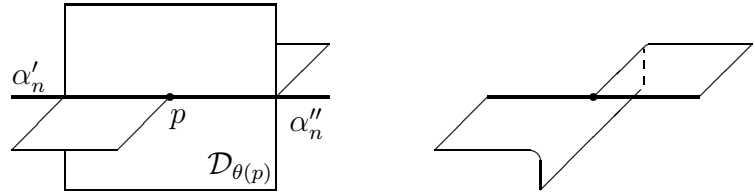
\begin{figure}[ht]\caption{The ribbon $R$ near a boundary saddle}\label{r1}
\vskip20pt
\centerline{\begin{picture}(120,70)
\put(60,35){\circle*{3}}
{\linethickness{1pt}\put(0,35){\line(1,0){120}}}
\put(0,15){\line(1,1){20}}
\put(40,15){\line(1,1){20}}
\put(100,35){\line(1,1){20}}
\put(0,15){\line(1,0){40}}
\put(100,55){\line(1,0){20}}
\put(60,25){$p$}
\put(0,40){$\alpha_n'$}
\put(105,22){$\alpha_n''$}
\put(20,0){\line(0,1){15}}
\put(20,35){\line(0,1){35}}
\put(100,0){\line(0,1){70}}
\put(20,70){\line(1,0){80}}
\put(20,0){\line(1,0){80}}
\put(74,6){$\disk_{\theta(p)}$}
\end{picture}\hskip40pt
\begin{picture}(120,70)
\put(60,35){\circle*{3}}
{\linethickness{1pt}\put(20,35){\line(1,0){80}}}
\put(0,15){\line(1,1){20}}
\put(100,35){\line(1,1){20}}
\put(0,15){\line(1,0){35}}
\put(120,55){\line(-1,0){38}}
\put(82,52){\oval(6,6)[tl]}
\multiput(79,49)(0,-6){2}{\line(0,-1){3}}
\put(40,10){\line(0,-1){10}}
\put(35,10){\oval(10,10)[tr]}
\put(40,0){\line(1,1){38}}
\put(60,35){\line(1,1){19}}
\end{picture}}
\end{figure}
We attach a small disk ``perpendicular'' to $\alpha_n$
along the segment $\{r_u(p)\}_{u\in[-\varepsilon,\varepsilon]}$
and then smooth the obtained surface
as shown in Fig.~\ref{r1} on the right, still denoting the
result by $R_1$. The point $p$ will be the only boundary
saddle of $\cal F$.

The union $R_1$ of disks is only a part of the ribbon under
construction. Near any vertex, $R_1$ looks as shown
in Fig.~\ref{r2} on the left.
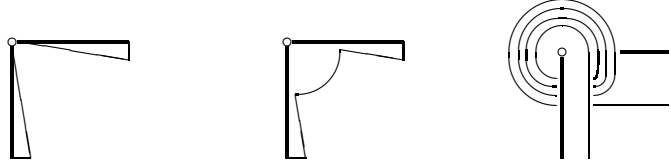
\begin{figure}[ht]\caption{The ribbon $R$ near a vertex of $L$}
\label{r2}
\centerline{
\begin{picture}(44,60)
\put(0,44){\circle{3}}
{\linethickness{1pt}\put(2,44){\line(1,0){42}}
\put(0,42){\line(0,-1){42}}}
\put(0,42){\line(1,-6){7}}
\put(0,0){\line(1,0){7}}
\put(2,44){\line(6,-1){42}}
\put(44,44){\line(0,-1){7}}
\end{picture}\hskip60pt
\begin{picture}(44,60)
\put(0,44){\circle{3}}
{\linethickness{1pt}\put(2,44){\line(1,0){42}}
\put(0,42){\line(0,-1){42}}}
\put(3,24){\line(1,-6){4}}
\put(0,0){\line(1,0){7}}
\put(20,41){\line(6,-1){24}}
\put(44,44){\line(0,-1){7}}
\put(3,41){\oval(34,34)[br]}
\end{picture}\hskip40pt
\begin{picture}(60,60)
\put(20,40){\circle{3}}
{\linethickness{1pt}\put(20,0){\line(0,1){38}}}
\put(30,0){\line(0,1){40}}
\put(20,40){\oval(20,20)[t]}
\put(18,40){\oval(16,20)[bl]}
\bezier50(32,30)(33.5,30)(33.5,40)
\put(20,40){\oval(27,26)[t]}
\put(18,40){\oval(23,26)[bl]}
\put(32,40){\oval(9,26)[br]}
\put(20,40){\oval(33,33)[t]}
\put(18,40){\oval(29,33)[bl]}
\put(32,40){\oval(16,33)[br]}
\put(20,40){\oval(40,40)[t]}
\put(18,40){\oval(36,40)[bl]}
\put(32,20){\line(1,0){28}}
{\linethickness{1pt}\put(42,40){\line(1,0){18}}}
\end{picture}}
\end{figure}
We transform it into a ribbon as shown in Fig.~\ref{r2} in
the middle, attaching a sector of a small disk
transversal to the binding circle. If we do this
near each vertex, then, possibly, the obtained
ribbon will be twisted, that is,
the connected components of $\partial R$ will be linked.
We can compensate this by twisting the ribbon
around the binding circle ``in the opposite direction''
at an arbitrarily chosen vertex as shown
in Fig.~\ref{r2} on the right.
This will create a winding vertex, so, if a boundary saddle
is present, we do this at one of the boundary vertices
adjacent to the saddle.
\end{proof}

\subsection{Moves of the characteristic surface}
We will think about an admissible characteristic surface $M$
of an arc-presentation $L$ as a periodic flow
of planar pictures that are intersections
of $M$ with the pages $\disk_t$, $t\in S^1$.
The pages $\disk_t$ containing singularities
of $\cal F$ and arcs of $L$ will be called \emph{singular},
the corresponding values of $t$ \emph{critical}.
All the others are \emph{regular}.

For simplicity, we assume in this subsection
that $\cal F$ has no poles and no closed
regular fibres. The notions and the claims of this subsection
will be used only when poles and closed fibres are absent.

If $\disk_t$ is a regular page, then the intersection $\disk_t\cap(L\cup M)$
is the union of pairwise non-intersecting arcs having
endpoints on the binding circle $S^1_\varphi=\partial\disk_t$.
Notice that two arcs in $\disk_t\cap M$ may have
a common endpoint if $L$ is an arc-presentation of the unknot
(refer to the right picture in Fig.~\ref{r2}), but
no arc forms a loop, which follow from an orientation argument.
Each section $\disk_t\cap(L\cup M)$ will be regarded
up to a homeomorphism of $\overline{\disk_t}$ fixed at $\partial\disk_t$,
which we shall express by saying that we consider the \emph{combinatorial
type} of the section $\disk_t\cap(L\cup M)$. This combinatorial
type can be specified by listing pairs of endpoints of all
arcs in $\disk_t$, like: $s_{i_1}s_{i_2}$, $s_{i_3}s_{i_4}$,
$s_{i_5}s_{i_6}$, \ldots.
The notation $s_is_j$ will be used for any arc with
endpoints $s_i,s_j\in S^1_\varphi$.
Notice: $s_is_j=s_js_i$, but $(s_i,s_j)\cap(s_j,s_i)=\varnothing$.
We shall also refer to the combinatorial type
of a section $\disk_t\cap(L\cup M)$ as a \emph{state}
of $(L,M)$.

We leave the following statement, which is standard, without proof.

\begin{prop}\label{combdescr}
Let $M$ and $M'$ be admissible characteristic surfaces for
arc-pre\-sen\-ta\-tions $L$ and $L'$, respectively, such that
we have $L\cap S^1_\varphi=L'\cap S^1_\varphi$,
$M\cap S^1_\varphi=M'\cap S^1_\varphi$, and
for all $t\in S^1_\theta$, the sections
$\disk_t\cap(L\cup M)$ and $\disk_t\cap(L'\cup M')$
have the same combinatorial type.
Then there is an isotopy $f$ fixed at $S^1_\varphi$,
preserving each page $\disk_t$, and taking
$(L,M)$ to $(L',M')$.
\end{prop}

If $\disk_t$ is a regular page, then the combinatorial type
of the section $\disk_t\cap(L\cup M)$ does not change if
$t$ varies slightly.
In addition to simple arcs, a singular page may contain
one of the following:
\begin{enumerate}
\item a T-joint, which is a boundary saddle of $\cal F$
with three arcs attached;
\item an X-joint, which is either an interior saddle of $\cal F$
with four arcs attached or a regular fibre of $\cal F$
crossed by an arc of $L$ (in the case of a composite link).
\end{enumerate}

When $t$ passes a critical value, the state of $(L,M)$ changes.
We shall call such a change
an \emph{event}. Below is the list of all possible events.

\medskip
\centerline{
\begin{tabular}{|p{4cm}|p{4cm}|p{4cm}|}
\hline
Event&Singularity&Case\\\hline
a single arc appears and immediately disappears&
an arc of $L$&
split link, composite link\\\hline
a single arc appears&
an arc of $L$&
unknot\\\hline
a single arc disappears&
an arc of $L$&
unknot\\\hline
an arc $s_1s_2$ is replaced by some arc
$s_2s_3$&
an arc $s_1s_3$ of $L$ containing a boundary saddle&
unknot\\\hline
two arcs $s_1s_2$, $s_3s_4$
are replaced by $s_1s_4$, $s_2s_3$&
a saddle of $\cal F$&
any\\\hline
\end{tabular}}

\medskip

If an arc with endpoints $s,s'\in S^1_\varphi$ is contained
in $\disk_t\cap(L\cup M)$, where $\disk_t$ is a singular page,
but there is no arc with endpoints $s,s'$ ``just before''
or ``just after'' the corresponding events, {\it i.e.},~in
the section $\disk_{t-\varepsilon}\cap(L\cup M)$ or
$\disk_{t+\varepsilon}\cap(L\cup M)$,
where $\varepsilon$ is sufficiently
small, we shall say that the arc $ss'$ \emph{participates} in the
corresponding event.

Clearly, if we know the state of $(L,M)$ just
before an event and have a description of
the event in terms of the participating arcs,
then we can recover the state of $(L,M)$ right after
the event. So, in order to specify the combinatorial type
of an arc-presentation
$L$ with a characteristic surface $M$,
we must describe the state of $(L,M)$
at just one non-critical moment
and provide the ordered list of all the events that occur when $t$
runs over the circle. Such a state with a list of events
will be called \emph{combinatorial description} of $(L,M)$.

We shall use combinatorial descriptions whenever
we will need to change the surface $M$ or an arc-presentation $L$.
Of course, an arbitrary initial state equipped with an
ordered list of events does not necessarily define
an arc-presentation and a characteristic surface;
the combinatorial data must be \emph{consistent}. It is easy, but
tedious, to list all the requirements for the data to be consistent.
Here are some of them: the final state (to which we come
after all the events) must coincide with the initial one,
an arc that is supposed to disappear in a forthcoming
event must be present just before the event, etc.
In each case when we change the combinatorial description
of the pair $(L,M)$, it will be clear that the
combinatorial data remain consistent.

\begin{lemma}\label{intevents}
Let $E_1,E_2$ be successive events such that
no arc participating in $E_1$ coincides or
interleaves with an arc participating in $E_2$.
If both events correspond to arcs of $L$, we
assume additionally that those arcs do
not have a common endpoint.
Then the events can be interchanged,
which will have no effect on the foliation $\cal F$
and the isotopy class of $(L,M)$.
If both singular pages corresponding to the events $E_1,E_2$
contain an arc of $L$, then the exchange of events
will result in one arc exchange move on $L$. Otherwise, there will be
no effect on the combinatorial type of $L$.
\end{lemma}

\begin{proof}
More precisely, the first assertion means the following.
Let $t_1=\theta(E_1)$ and $t_2=\theta(E_2)$ be the time
of the events. Then there exists an arc-presentation $L'$ and
an admissible characteristic surface $M'$ for $L'$ such that,
for all $t\in(t_2,t_1)$ the states
of $(L,M)$ and $(L',M')$ at the moment $t$ coincide.
At the moments $t_1,t_2$,
the events $E_2,E_1$, respectively, happen with
$\disk_t\cap(L'\cup M')$, and there is no event in between.

By saying that a transform has no effect on the foliation
$\cal F$, we mean that the new foliation $\cal F'$ is equivalent
to $\cal F$ via a homeomorphism $M\rightarrow M'$.

The assumptions of the lemma imply that there exists
a continuous family of arcs $\gamma_{t,s_1,s_2}\subset\disk_t$,
$t\in(t_1-\varepsilon,t_2+\varepsilon)$,
$\partial\gamma_{t,s_1,s_2}=\{s_1,s_2\}\subset S^1_\varphi$,
forming an open disk $d$ that does not intersect the surface $M$ and
the link $L$, such that the event $E_1$ occurs ``on the right'' of
the disk $d$ and the event $E_2$ occurs ``on the left'' of $d$.
The latter means, in particular, that all arcs participating
in $E_1$ (respectively, $E_2$)
have both endpoints on the interval $[s_1,s_2]$ (respectively, in
$[s_2,s_1]$). It is not forbidden that an arc participating
in the events have one or both endpoints on $\{s_1,s_2\}$.
If $s_1$ and $s_2$ are the endpoints of an arc participating
in $E_1$ or $E_2$, the interior of the arc must be on the appropriate
side of $\gamma_{t,s_1,s_2}$, where the corresponding event occurs.

Now we can define a self-homeomorphism $f$ of $S^3\setminus d$
preserving the foliation $d\theta=0$
by making the time $\theta$ ``go faster''
on the left of $d$ and ``go slower'' on the right of $d$,
which will result in $E_2$ occurring
before $E_1$. In the pages $\disk_t$ with $t\in[t_2+\varepsilon,
t_1-\varepsilon]$ the mapping $f$ is assumed to be identical.
The mapping $f$ will send the surface $M$
to another one, preserving the foliation. Clearly, the mapping
$f$ can be adjusted in a small neighbourhood of $d$
to become a self-homeomorphism of the whole $S^3$.

The last two assertions of the lemma are obvious.
\end{proof}

\begin{rem} We do \emph{not} exclude the case when an arc participating
in $E_1$ has a common endpoint with one participating in $E_2$.
This may occur when $L$ is an arc-presentation
of the unknot, $\partial M=L$,
and the common end of two arcs is a winding vertex of $\cal F$
(refer again to Fig.~\ref{r2}).
This is why we wrote $[s_1,s_2]$ and $[s_2,s_1]$ instead
of $(s_1,s_2)$, $(s_2,s_1)$ in the above proof.

\end{rem}

Whereas the transform described in
Lemma~\ref{intevents} generalizes an arc exchange move
to the case of an arc-presentation endowed with an
admissible characteristic surface, the following
construction is an analogue of a generalized
exchange move.

\begin{lemma}\label{bigexchange}
Let $s_1,s_2,s_3$ be three points of $S^1_\varphi$ disjoint
from $L$ and $M$ such that $s_2\in(s_1,s_3)$.
Suppose that there exist moments $t_1,t_2\in
S^1_\theta$ such that no arc in $\disk_t\cap(L\cup M)$
interleaves with $s_1s_2$ if $t\in[t_1,t_2]$
and with $s_2s_3$ if $t\in[t_2,t_1]$.
Then one can interchange all the vertices of $L$ and $\cal F$ lying
in $(s_1,s_2)$ with those in $(s_2,s_3)$, keeping
the relative order of vertices in each of the intervals fixed,
without change in the foliation $\cal F$ and of the isotopy
class of $(L, M)$.
\end{lemma}

Interchanging the vertices should be understood in the same
way as described in~\ref{generalized}.

\begin{proof}
Denote by $d_{t_1,t_2,s_1,s_2,s_3}$
a closed two-dimensional disk whose interior consists
of the vertex $s_2$ and the following four disks $d_1,d_2,d_3,d_4$
adjacent to $s_2$:
$$d_1=\cup_{t\in(t_1,t_2)}\gamma_{t,s_1,s_2},\qquad d_2=
\cup_{t\in(t_2,t_1)}\gamma_{t,s_2,s_3},$$
$d_3$ $d_4$ are the two ``triangles'' bounded by
$\gamma_{t_i,s_1,s_2}$, $\gamma_{t_i,s_2,s_3}$,
$\gamma_{t_i,s_1,s_3}$, where $i=1,2$.

The assumptions of the Lemma imply that the disk
$d^{(1)}=d_{t_1,t_2,s_1,s_2,s_3}$ is disjoint from $L$ and $M$.
Let $L'$ and $M'$ be the arc-presentation
and the surface obtained after the transform.
Then the disk $d^{(2)}=d_{t_2,t_1,s_1,s_1+s_3-s_2,s_3}$
is disjoint from $L'$ and $M'$.
It is easy to see that there exists a homeomorphism
$f:S^3\setminus d^{(1)}\rightarrow S^3\setminus d^{(2)}$
preserving each page $\disk_t$ and interchanging
the intervals $(s_1,s_2)$ and $(s_2,s_3)$ of the
binding circle. We have $f(L)=(L')$ and $f(M)= M'$.
There is a small neighbourhood $U$ of $d^{(1)}$ that
does not meet $L$ and $M$. By adjusting $f$ in $U$
we can make $f$ a self-homeomorphism of the whole sphere $S^3$.\end{proof}

\subsection{Patterns of $\cal F$ that are always present}

Let $L$ be an arc-presentation, $M$ an admissible characteristic
surface of $L$, and $\cal F$ the foliation on $M$ defined by
$d\theta=0$.  We define the \emph{complexity} $c(M)$
of $M$ as the total number of singularities of $\cal F$.

\emph{Fibres} of $\cal F$ are of the following types:
\begin{itemize}
\item a closed circle;
\item an open arc connecting two vertices;
\item an open arc connecting a vertex to a saddle
or a saddle to itself (we call such
a fibre a \emph{separatrix});
\item a singular point of $\cal F$.
\end{itemize}

Two regular fibres connecting the same pair of vertices
will called \emph{parallel} if they are included into
a continuous family of such fibres (\emph{i.e.} enclose
a disk with no singularity inside). Any other two fibres
connecting vertices will be called \emph{non-parallel}.

For a vertex $v$ of $\cal F$, we call the closure of the union
of all fibres of $\cal F$ approaching $v$
the \emph{star} of $v$. By the \emph{valence} of
$v$ we mean the number of separatrices in the star of $v$,
which coincides with the number of saddles in the star.

An interior vertex $v$ of $\cal F$ is said to be \emph{bad}
in three cases: 1) the star of $v$ contains at least two non-parallel
fibres each connecting $v$ to a boundary vertex;
2) the star of $v$ is pierced by $L$; 3) the star of $v$ contains a winding
vertex.
If a vertex is not bad, it is said to be \emph{good}.

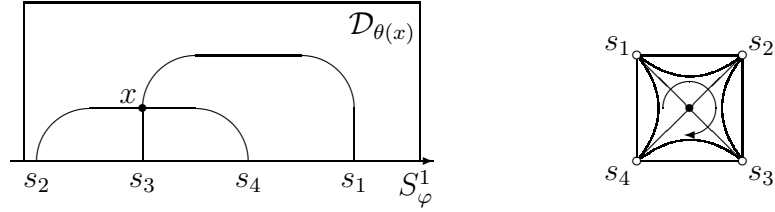
\begin{figure}[ht]
\caption{Orientation of a saddle}\label{orient}
\vskip2em
\centerline{
\begin{picture}(160,80)
\put(0,20){\vector(1,0){160}}
\put(160,7){\hbox to0pt{\hss$S^1_\varphi$}}
\put(5,20){\line(0,1){60}}
\put(5,80){\line(1,0){150}}
\put(155,20){\line(0,1){60}}
\put(128,68){$\disk_{\theta(x)}$}
\put(50,20){\oval(80,40)[t]}
\put(90,20){\oval(80,80)[t]}
\put(50,40){\circle*3}
\put(10,10){\hbox to0pt{\hss$s_2$\hss}}
\put(50,10){\hbox to0pt{\hss$s_3$\hss}}
\put(90,10){\hbox to0pt{\hss$s_4$\hss}}
\put(130,10){\hbox to0pt{\hss$s_1$\hss}}
\put(48,42){\hbox to0pt{\hss$x$}}
\end{picture}\hskip2cm
\begin{picture}(80,80)
\put(20,20){\circle3}
\put(20,60){\circle3}
\put(60,20){\circle3}
\put(60,60){\circle3}
\put(40,40){\oval(20,20)[t]}
\put(40,40){\oval(20,20)[br]}
\put(38,30){\vector(-1,0){0}}
\put(40,40){\circle*3}
\put(18,18){\vbox to0pt{\hbox to0pt{\hss$s_4$}\vss}}
\put(62,18){\vbox to0pt{\hbox to0pt{$s_3$\hss}\vss}}
\put(18,62){\hbox to0pt{\hss$s_1$}}
\put(62,62){\hbox to0pt{$s_2$\hss}}
\put(21.5,21.5){\line(1,1){37}}
\put(58.5,21,5){\line(-1,1){37}}
\bezier150(21.5,21)(40,35)(58.5,21)
\bezier150(21.5,59)(40,45)(58.5,59)
\bezier150(21,21.5)(35,40)(21,58.5)
\bezier150(59,21.5)(45,40)(59,58.5)
\put(21.5,20){\line(1,0){37}}
\put(21.5,60){\line(1,0){37}}
\put(20,21.5){\line(0,1){37}}
\put(60,21.5){\line(0,1){37}}
\end{picture}
}
\end{figure}

Each saddle $x$ will be endowed with an \emph{orientation}, which
is the circular order of the separatrices approaching $x$
induced by the order of the corresponding vertices on the axis $S^1_\varphi$,
see Fig.~\ref{orient}.

Now we consider the star of a winding vertex in more detail.

\begin{rem}
It was brought to the author's attention by Bill Menasco and Adam Sikora
that the presence of a winding vertex in the star of a good two-valent
vertex may produce an obstruction to applying the corresponding
generalized exchange move in the proof of Lemma~\ref{l5}.
The definitions of the characteristic surface and a bad vertex
have been changed,
the new Lemma~\ref{newlemma} below has been added, Lemma~\ref{cases} and
the proof of Lemma~\ref{l5} have been modified by introducing new
cases 6)--8) in order correct the mistake in the earlier version of the paper.
\end{rem}

Suppose that $M$ has a winding vertex. Denote it by $s^{\mathrm w}$ and
its star by $\sigma^{\mathrm w}$.
Let $s^{\mathrm w}_1,\dots,s^{\mathrm w}_q$ be
the vertices, and $x^{\mathrm w}_1,\dots,x^{\mathrm w}_{q-1}$
the interior saddles in $\partial\sigma^{\mathrm w}\setminus s^{\mathrm w}$
numbered
so that $x^{\mathrm w}_{i}$ lies between $s^{\mathrm w}_i$ and
$s^{\mathrm w}_{i+1}$, $i=1,\dots,q-1$. If there is a boundary saddle, we
assume that it is adjacent to $s^{\mathrm w}_q$ and denote it by
$x^{\mathrm w}_q$. Thus, the valence of $s^{\mathrm w}$ is equal to
$$q'=\left\{\begin{aligned}q-1&\quad\text{if there is no boundary saddle},\\
q&\quad\text{otherwise}.\end{aligned}\right.$$

Note that if a vertex from $\sigma^{\mathrm w}$
is connected to $s^{\mathrm w}$ by two or more non-parallel
fibres, it will appear in the sequence
$s^{\mathrm w}_1,\dots,s^{\mathrm w}_q$ more than once. So,
some of the $s^{\mathrm w}_i$ may coincide.

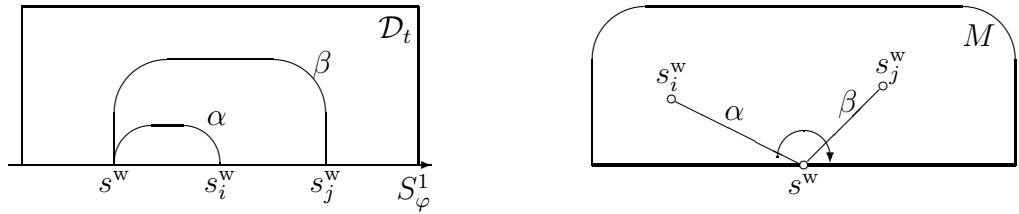
\begin{figure}[ht]
\caption{Orientation of the winding vertex}\label{windorient}
\vskip2em
\centerline{
\begin{picture}(160,80)
\put(0,20){\vector(1,0){160}}
\put(160,7){\hbox to0pt{\hss$S^1_\varphi$}}
\put(5,20){\line(0,1){60}}
\put(5,80){\line(1,0){150}}
\put(155,20){\line(0,1){60}}
\put(140,68){$\disk_t$}
\put(60,20){\oval(40,30)[t]}
\put(80,20){\oval(80,80)[t]}
\put(40,10){\hbox to0pt{\hss$s^{\mathrm w}$\hss}}
\put(80,10){\hbox to0pt{\hss$s^{\mathrm w}_i$\hss}}
\put(120,10){\hbox to0pt{\hss$s^{\mathrm w}_j$\hss}}
\put(75,35){$\alpha$}
\put(115,55){$\beta$}
\end{picture}
\hskip2cm
\begin{picture}(160,80)
\put(80,20){\circle3}
{\linethickness{1pt}\put(0,20){\line(1,0){78.5}}\put(160,20){\line(-1,0){78.5}}}
\put(80,9){\hbox to0pt{\hss$s^{\mathrm w}$\hss}}
\put(80,20){\oval(160,120)[t]}
\put(140,65){$M$}
\put(80,23){\oval(20,20)[t]}
\put(90,21){\vector(0,-1)0}
\put(30,45){\circle3}
\put(110,50){\circle3}
\put(81.5,21.5){\line(1,1){27}}
\put(31.5,44){\line(2,-1){47}}
\put(30,51){\hbox to0pt{\hss$s^{\mathrm w}_i$\hss}}
\put(113,55){\hbox to0pt{\hss$s^{\mathrm w}_j$\hss}}
\put(50,37){$\alpha$}
\put(93,40){$\beta$}
\end{picture}
}
\end{figure}

We also introduce the \emph{oriention} of the winding vertex,
which is either a locally $\theta$-increasing or a locally
$\theta$-decreasing linear
ordering $\prec$ of fibres approaching $s^{\mathrm w}$.
(This orientation is not the one
induced by the orientation of $S^1_\theta$, which is often
used in the work of Birman and Menasco.)
Take an arbitrary page $\disk_t$ such that at least two arcs
in $\disk_t\cap M$ have $s^{\mathrm w}$ as an endpoint.
Let $\alpha=\gamma_{t,s^{\mathrm w},s^{\mathrm w}_i}$,
$\beta=\gamma_{t,s^{\mathrm w},s^{\mathrm w}_j}$,
be such arcs. We orient $s^{\mathrm w}$ so as to have
$\alpha\prec\beta$ if and only if
$s^{\mathrm w}_i\in(s^{\mathrm w},s^{\mathrm w}_j)$,
see Fig.~\ref{windorient}.
This orientation does not depend on the choice of the
page $\disk_t$ and of the pair $\alpha,\beta$ in $\disk_t$,
and can be described informally as follows.
A small part of the surface $M$ near the winding vertex
forms a ``screw'' whose axis is $S^1_\varphi$. The orientation
of $s^{\mathrm w}$ indicates which way this screw
must rotate in order to screw up in the direction opposite
to the orientation of $S^1_\varphi$.

Now we put $\delta_i=+1$, if the saddle $x^{\mathrm w}_i$
is coherently oriented with the winding
vertex, and $\delta_i=-1$ otherwise.

\begin{lemma}\label{newlemma}
If two successive interior vertices $s^{\mathrm w}_i$ and $s^{\mathrm w}_{i+1}$
in the star of the winding vertex are both two-valent, then
$\delta_i=+1$.
\end{lemma}

\begin{proof}
Let $s$ be the vertex in the stars of
$s^{\mathrm w}_i$ and $s^{\mathrm w}_{i+1}$
different from $s^{\mathrm w}$. The union
of these stars is shown in Fig.~\ref{twotwo} on the left.
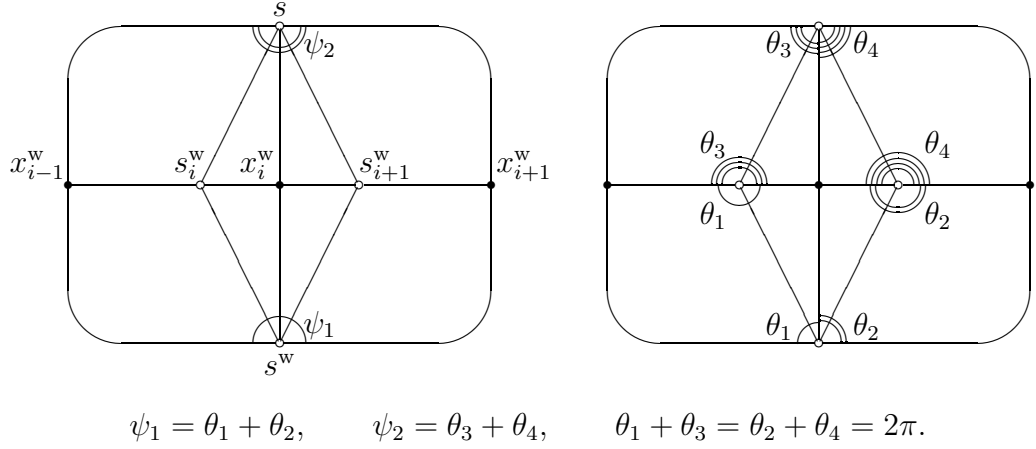
\begin{figure}[ht]
\caption{Two successive two-valent vertices in $\sigma^{\mathrm w}$}\label{twotwo}
\centerline{\hskip20pt
\begin{picture}(200,160)
\put(20,80){\circle*3}
\put(100,80){\circle*3}
\put(180,80){\circle*3}
\put(100,20){\circle3}
\put(100,140){\circle3}
\put(70,80){\circle3}
\put(130,80){\circle3}
\put(20,80){\line(1,0){48}}
\put(180,80){\line(-1,0){48}}
\put(71.5,80){\line(1,0){57}}
\put(100,21.5){\line(0,1){117}}
\put(100.5,21){\line(1,2){29}}
\put(99.5,21){\line(-1,2){29}}
\put(70.5,81){\line(1,2){29}}
\put(129.5,81){\line(-1,2){29}}
\put(101.5,80){\oval(157,120)[r]}
\put(98.5,80){\oval(157,120)[l]}
\put(100,8){\hbox to0pt{\hss$s^{\mathrm w}$\hss}}
\put(60,85){$s^{\mathrm w}_i$}
\put(130,85){$s^{\mathrm w}_{i+1}$}
\put(85,85){$x^{\mathrm w}_i$}
\put(182,85){$x^{\mathrm w}_{i+1}$}
\put(-2,85){$x^{\mathrm w}_{i-1}$}
\put(100,144){\hbox to0pt{\hss$s$\hss}}
\put(100,20){\oval(20,20)[t]}
\put(100,140){\oval(20,20)[b]}
\put(100,140){\oval(16,16)[b]}
\put(109,130){$\psi_2$}
\put(109,24){$\psi_1$}
\end{picture}
\begin{picture}(200,160)
\put(20,80){\circle*3}
\put(100,80){\circle*3}
\put(180,80){\circle*3}
\put(100,20){\circle3}
\put(100,140){\circle3}
\put(70,80){\circle3}
\put(130,80){\circle3}
\put(20,80){\line(1,0){48}}
\put(180,80){\line(-1,0){48}}
\put(71.5,80){\line(1,0){57}}
\put(100,21.5){\line(0,1){117}}
\put(100.5,21){\line(1,2){29}}
\put(99.5,21){\line(-1,2){29}}
\put(70.5,81){\line(1,2){29}}
\put(129.5,81){\line(-1,2){29}}
\put(101.5,80){\oval(157,120)[r]}
\put(98.5,80){\oval(157,120)[l]}
\put(100,20){\oval(21,21)[tr]}
\put(100,20){\oval(17,17)[tr]}
\put(100,20){\oval(16,16)[tl]}
\put(100,140){\oval(13,13)[bl]}
\put(100,140){\oval(12,12)[br]}
\put(100,140){\oval(17,17)[bl]}
\put(100,140){\oval(16,16)[br]}
\put(100,140){\oval(21,21)[bl]}
\put(100,140){\oval(20,20)[br]}
\put(100,140){\oval(24,24)[br]}
\put(70,80){\oval(17,17)[t]}
\put(70,80){\oval(16,16)[b]}
\put(70,80){\oval(13,13)[t]}
\put(70,80){\oval(21,21)[t]}
\put(130,80){\oval(16,16)[t]}
\put(130,80){\oval(20,20)[t]}
\put(130,80){\oval(17,17)[b]}
\put(130,80){\oval(21,21)[b]}
\put(130,80){\oval(12,12)[t]}
\put(130,80){\oval(24,24)[t]}
\put(80,23){$\theta_1$}
\put(113,23){$\theta_2$}
\put(80,130){$\theta_3$}
\put(113,130){$\theta_4$}
\put(55,65){$\theta_1$}
\put(140,65){$\theta_2$}
\put(55,92){$\theta_3$}
\put(140,92){$\theta_4$}
\end{picture}
}
\vskip-0.5cm
$$\psi_1=\theta_1+\theta_2,\qquad\psi_2=\theta_3+\theta_4,\qquad
\theta_1+\theta_3=\theta_2+\theta_4=2\pi.$$
\end{figure}
Let $\psi_1$ be the ``$\theta$-angle'' between separatrices $s^{\mathrm w}
x^{\mathrm w}_{i-1}$ and $s^{\mathrm w}x^{\mathrm w}_{i+1}$ at
$s^{\mathrm w}$, \emph{i.e.}, the absolute value of the integral of
$d\theta$ along an arc indicated in the picture. Similarly,
let $\psi_2$ be the $\theta$-angle between
separatrices $sx^{\mathrm w}_{i-1}$ and $sx^{\mathrm w}_{i+1}$
at $s$.

It is not hard to see that $\psi_1+\psi_2=4\pi$ (see Fig.~\ref{twotwo}).
Since $s$ is not a winding vertex, we have $\psi_2<2\pi$.
Therefore, $\psi_1>2\pi$, which implies that some page
$\disk_t$ contains both arcs $\gamma_{t,s^{\mathrm w},s^{\mathrm w}_i}$
and $\gamma_{t,s^{\mathrm w},s^{\mathrm w}_{i+1}}$.

If $s^{\mathrm w}_i\in(s^{\mathrm w},s^{\mathrm w}_{i+1})$,
then the orientations of
the saddle $x^{\mathrm w}_i$ and of the winding vertex
are both clockwise (according to Fig.~\ref{twotwo}).
If $s^{\mathrm w}_{i+1}\in(s^{\mathrm w},s^{\mathrm w}_i)$,
then both orientations are counterclockwise. So, in both cases
they are coherent.
\end{proof}

The following lemma is a generalization of Lemma~2 of~\cite{bm4}.
(Although closed braids are now replaced with arc-presentations,
the argument of~\cite{bm4} does not need to be changed to establish
the claim below in the case of a split or composite link.)

\begin{lemma}\label{cases}
If we have $c(M)>0$, then the foliation $\cal F$
has at least one of the following:
\rm\begin{enumerate}
\item \it a pole\/\rm;
\item \it a univalent boundary vertex\/\rm;
\item \it a good two-valent interior vertex\/\rm;
\item \it a good three-valent interior vertex\/\rm;
\item \it two pairs of two-valent vertices, in each the vertices
are connected by a fibre of $\cal F$ intersecting $L$\rm;
\item \it a winding vertex and a three-valent interior vertex $s^{\mathrm w}_i$
in $\sigma^{\mathrm w}$\rm;
\item \it a winding vertex and a two-valent boundary vertex
$s^{\mathrm w}_i$ in $\sigma^{\mathrm w}$
such that $\delta_{i-1}=\delta_i$\rm;
\item \it a winding vertex and a four-valent interior
vertex
$s^{\mathrm w}_i$ in $\sigma^{\mathrm w}$
such that $\delta_{i-1}=\delta_i$\rm.
\end{enumerate}
\end{lemma}

\begin{proof}
We denote by $\chi$ the Euler characteristic of $M$,
$\chi\in\{1,2\}$. We also set $\epsilon=1$ if
$\cal F$ has a boundary saddle, and $\epsilon=0$ otherwise.
Denote the number of interior vertices of valence $k$ by $V_k$,
the number of boundary vertices of valence $k$ by $V_k^{\mathrm b}$,
the number of interior saddles by $S$, and the number of poles by $P$.
By the definition of complexity, we have
\begin{equation}\label{c(m)}
S+P+\epsilon=c(M).
\end{equation}

Notice that we always have $V_1=0$. Indeed, the opposite
would mean that there is an event $E$ in the flow
$\disk_t\cap(L\cup M)$ corresponding to a saddle
such that some arc participating in $E$ and existing just before
the event survives after the event, which is impossible.
Notice also that either $V_0>0$ or $V_0^{\mathrm b}>0$
imply $c(M)=0$: if there
is no singularity in the star of a vertex, then this star
is the whole surface $M$, and there are exactly two
vertices, either interior or boundary ones.

Thus, we may assume that $V_0=V_0^{\mathrm b}=V_1=0$.
If $P>0$, then case~1) takes place, and we are done. So,
we assume $P=0$.

Counting the topological indices of singularities gives
\begin{equation}\label{chi}
\sum\limits_kV_k+\frac12\sum\limits_kV_k^{\mathrm b}
-S-\frac\epsilon2=\chi.
\end{equation}

Counting the number of separatrices in two different ways
gives
\begin{equation}\label{sep}
\sum\limits_kk(V_k+V_k^{\mathrm b})=4S+3\epsilon.
\end{equation}

From~(\ref{chi}) and (\ref{sep})
we get
\begin{equation}\label{main}
\sum\limits_k(4-k)V_k+\sum\limits_k(2-k)V_k^{\mathrm b}=
4\chi-\epsilon.
\end{equation}

(i)
First, consider the case of a split link, which
is the simplest one because,
in this case, there are no boundary vertices and no bad vertices.
Relation~(\ref{main}) can be rewritten as
\begin{equation}\label{main1}
2V_2+V_3=8+\sum\limits_{k\geqslant4}(k-4)V_k,
\end{equation}
which implies
\begin{equation}\label{main2}
2V_2+V_3\geqslant8.
\end{equation}
We have either $V_2>0$ or $V_3>0$, which are Cases 3) and 4) of
the lemma.

(ii)
The case of a composite link is similar, but now the
inequality~(\ref{main2}) does not imply that there exists a
\emph{good} three- or two-valent vertex, because some
vertices are bad. There may be at most four bad vertices.
So, either there exists a good three- or two-valent vertex or
we have $V_2+V_3\leqslant4$. In conjunction with~(\ref{main2}),
the latter implies $V_2=4$, $V_3=0$, {\it i.e.}, all
bad vertices are two-valent. This is case~5) of the lemma.

\begin{rem}
In this case, we also have $V_5=V_6=\ldots=0$, but not necessarily
$V_4=0$. One can show that $V_4$ must be even, but it \emph{can}
be arbitrarily large. The case $V_4>0$ is missing
in~\cite{bm4} and also in~\cite{cromwell}.
\end{rem}

(iii)
Now we consider the case of the unknot. Relation~(\ref{main})
reads:
\begin{equation}\label{main3}
2V_2+V_3+V_1^{\mathrm b}=4-\epsilon+\sum_{k\geqslant4}(k-4)V_k+
\sum_{k\geqslant2}(k-2)V_k^{\mathrm b}.
\end{equation}
If $V_1^{\mathrm b}>0$, then case~2) of the lemma takes place.

Suppose that $V_1^{\mathrm b}=0$. Then it follows from~(\ref{main3})
that $2V_2+V_3>0$. If all vertices are good (and hence there
is no winding vertex), we are done:
either case~3) or 4) of the lemma takes place.

(iii.a) Assume that some vertices of $\cal F$ are bad, but there is no winding
vertex.
In each family of parallel fibres of $\cal F$ connecting
a bad vertex with a boundary vertex, we select one fibre and
cut the surface $M$ along all the selected arcs
(see Fig.~\ref{cutting}).
\begin{figure}[ht]\caption{Cutting the disk $M$}\label{cutting}
\vskip20pt
\centerline{\begin{picture}(140,60)
\put(22,10){\oval(20,20)[r]}
\put(18,20){\oval(36,40)[l]}\put(20,0){\circle{3}}
\put(20,20){\circle{3}}\put(20,40){\circle{3}}
\put(20,2){\line(0,1){16}}\put(20,22){\line(0,1){16}}
\put(22,0){\line(1,0){36}}\put(60,0){\circle{3}}
\put(60,2){\line(0,1){16}}\put(60,20){\circle{3}}
\put(22,40){\line(1,0){16}}\put(40,40){\circle{3}}
\put(41,39){\line(1,-1){18}}\put(60,42){\oval(40,36)[t]}
\put(62,0){\line(1,0){38}}\put(100,18){\oval(40,36)[br]}
\put(120,20){\circle{3}}\put(120,22){\line(0,1){16}}
\put(120,40){\circle{3}}\put(61,21){\line(1,1){18}}
\put(80,40){\circle{3}}\put(82,40){\line(1,0){36}}
\put(122,40){\line(1,0){16}}\put(120,42){\line(0,1){16}}
\put(120,60){\circle{3}}\put(140,40){\circle{3}}
\put(122,38){\oval(36,36)[br]}\put(122,42){\oval(36,36)[tr]}
\put(118,42){\oval(76,36)[tl]}
\end{picture}}\end{figure}
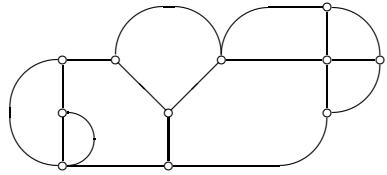
Let $d_1,\dots,d_l$ be the obtained disks. We say that
a disk $d_i$ is \emph{terminal} if there is
exactly one bad vertex of $\cal F$ in $\partial d_i$.
Clearly, there should be at least two terminal
disks ({\it e.g.}, there are six
in Fig.~\ref{cutting}). Therefore, we can find a terminal disk $d_i$
such that $\partial d_i$ does not contain
the boundary saddle of $\cal F$. Let $d_1$ be that disk.

Consider the restriction $\widetilde{\cal F}$ of the foliation $\cal F$
to the disk $d_1$. Let $\widetilde V_k$ and
$\widetilde V_k^{\mathrm b}$ be the number of
interior and boundary $k$-valent vertices of $\widetilde{\cal F}$,
respectively. By analogy with~(\ref{main3}), we have
$$
2\widetilde V_2+\widetilde V_3+\widetilde V_1^{\mathrm b}=
4+\sum_{k\geqslant4}(k-4)\widetilde V_k+
\sum_{k\geqslant2}(k-2)\widetilde V_k^{\mathrm b}\geqslant4.
$$
Since $V_1^{\mathrm b}=0$ and the disk $d_1$ is terminal, we have
$\widetilde V_1^{\mathrm b}\leqslant3$, which implies
$2\widetilde V_2+\widetilde V_3\geqslant1$. Thus, there is
at least one two- or three-valent interior vertex in $d_1$.
By construction, this vertex is a good vertex of $\cal F$,
so, one of the cases~3) or 4) takes place.

(iii.b)
Now assume that there is a winding vertex, but no
boundary saddle (\emph{i.e.} $\epsilon=0$). Assume also that all vertices
$s^{\mathrm w}_i$ in $\sigma^{\mathrm w}$ are pairwise distinct and
that Cases 6)--8) don't occur.
Let $V_k^{\mathrm w}$ (respectively, $V_k^{\mathrm{bw}}$)
be the number of $k$-valent interior (respectively, boundary) vertices
of $\mathcal F$ that lie
in $\sigma^{\mathrm w}$.
Consider the sequence $\Delta=(\delta_1,\delta_2,\dots,\delta_{q'})$
of $\pm1$s defined above. Let $m_+$ (respectively, $m_-$)
be the number of $+1$s (respectively, $-1$s) in $\Delta$,
and let $m_{--}$ be the number of subsequences $(-1,-1)$ in $\Delta$.
We have
\begin{equation}
m_{--}\geqslant m_--m_+-1.\label{m--}
\end{equation}
Whenever we have a subsequence of the form $(\delta_{i-1},\delta_i)=
(-1,-1)$ in $\Delta$, the vertex $s^{\mathrm w}_i$ cannot be
\begin{itemize}
\item
a two-valent interior vertex, since two saddles in the star of such
a vertex must be incoherently oriented (see~\cite{bm4});
\item a $k$-valent interior vertex with $2<k<5$
or a two-valent boundary vertex, by the assumption that cases 6)--8)
do not occur.
\end{itemize}
Thus the vertex $s^{\mathrm w}_i$ can be only a $k$-valent interior
vertex with $k\geqslant 5$ or a $k$-valent boundary vertex with
$k\geqslant3$. Taking into account that the winding vertex itself
has valence $q'=m_++m_-\geqslant2$, we obtain from (\ref{m--}):
\begin{equation}\label{2m-3}
\sum_{k\geqslant5}(k-4)V^{\mathrm w}_k+
\sum_{k\geqslant3}(k-2)V^{\mathrm{bw}}_k
\geqslant m_{--}+q'-2\geqslant2m_--3.
\end{equation}
It follows from Lemma~\ref{newlemma} that
\begin{equation}\label{v2m-}
V_2^{\mathrm w}\leqslant m_-.
\end{equation}
Combining (\ref{main}), (\ref{2m-3}), (\ref{v2m-}), and recalling
that $V_1^{\mathrm b}=0$, we obtain:
$$\begin{aligned}
2(V_2-V_2^{\mathrm w})+V_3&=
4+\sum_{k\geqslant5}(k-4)V_k+\sum_{k\geqslant3}(k-2)V^{\mathrm b}_k
-2V_2^{\mathrm w}\\
&\geqslant
4+\sum_{k\geqslant5}(k-4)V^{\mathrm w}_k+
\sum_{k\geqslant3}(k-2)V^{\mathrm{bw}}_k
-2V_2^{\mathrm w}
\geqslant1,
\end{aligned}
$$
which implies that there is
a two- or three-valent interior vertex in $M\setminus\sigma^{\mathrm w}$.
We are done unless there are bad vertices in $M\setminus\sigma^{\mathrm w}$.
In the latter case, we apply the same cutting as in case (iii.a)
and note that the presence of a bad vertex in $M\setminus\sigma^{\mathrm w}$
implies that one can find a terminal disk whose interior is disjoint from
$\sigma^{\mathrm w}$. Since the winding vertex
is away from such a terminal disk, we can apply the same argument as
in case (iii.a).

(iii.c)
Now let both a winding vertex and a boundary saddle
be present, but still assume that the vertices $s^{\mathrm w}_i$
are all distinct. Recall that, in this case,
the boundary saddle lies in $\sigma^{\mathrm w}$
by construction. Now we have $\epsilon=1$ and
$$
2(V_2-V_2^{\mathrm w})+V_3=
3+\sum_{k\geqslant5}(k-4)V_k+\sum_{k\geqslant3}(k-2)V^{\mathrm b}_k
-2V_2^{\mathrm w}=A+B+C,
$$
where
$$\begin{aligned}
A&=\sum_{k\geqslant5}(k-4)(V_k-V^{\mathrm w}_k)
+\sum_{k\geqslant3}(k-2)(V^{\mathrm b}_k-V^{\mathrm{bw}}_k),\\
B&=\sum_{k\geqslant5}(k-4)V^{\mathrm w}_k+
\sum_{k\geqslant3}(k-2)V^{\mathrm{bw}}_k-2m_-+3,\\
C&=2m_--2V^{\mathrm{bw}}_2.
\end{aligned}$$
By the same argument as in case (iii.b), we have $B,C\geqslant0$,
and by construction $A\geqslant0$.
So, we need to show that at least one of $A$, $B$, $C$ is non-zero.

Suppose the opposite, \emph{i.e.}, that we have
\begin{equation}
A=B=C=0.\label{abc}
\end{equation}
If $B=0$ holds, then~(\ref{m--}) must be sharp,
which implies $\delta_q=-1$. The equalities $C=0$,
$\delta_q=-1$, and Lemma~\ref{newlemma}
imply that $s^{\mathrm w}_q$ is two-valent. Let $s$ be the boundary vertex
different from $s^{\mathrm w}$
adjacent to the boundary saddle.
This vertex is not in $\sigma^{\mathrm w}$, so its valence must be equal
to two (otherwise we would have $A>0$).
Therefore, the star of $s$ does not contain any saddle except
$x^{\mathrm w}_q$ and $x^{\mathrm w}_{q-1}$, which implies that
$s^{\mathrm w}_{q-1}$ is a boundary vertex connected to $s$
by an arc of $L$, see Fig.~\ref{badcase}.
\begin{figure}[ht]
\caption{In case (iii.c), the equalities (\ref{abc})
imply $s^{\mathrm w}_{q-1}\in L$}\label{badcase}
\vskip2em
\centerline{\begin{picture}(120,120)
\put(20,20){\circle3}
\put(60,60){\circle3}
\put(20,100){\circle3}
\put(100,100){\circle3}
\put(21,21){\line(1,1){38}}
\put(61,61){\line(1,1){38}}
\multiput(21.5,98.5)(0.2,0)5{\oval(157,157)[br]}
\multiput(21.5,98.5)(0,-0.2)5{\oval(157,157)[br]}
\bezier400(20.5,21.5)(30,90)(98.5,99.5)
\put(20,21.5){\line(0,1){77}}
{\linethickness{1pt}\put(21.5,100){\line(1,0){77}}}
\put(21,99){\line(1,-1){38}}
\put(61,59){\line(1,-1){33}}
\put(94.5,25.5){\circle*3}
\put(45,75){\circle*3}
\put(18,8){\hbox to0pt{\hss$s^{\mathrm w}$}}
\put(96.5,13.5){$x^{\mathrm w}_q$}
\put(102,104){$s$}
\put(18,104){\hbox to0pt{\hss$s^{\mathrm w}_{q-1}$}}
\put(60,104){\hbox to0pt{\hss$L$\hss}}
\put(102,58){$L$}
\put(81,23){$-$}
\put(31,71.5){$+$}
\put(48,72){$x_{q-1}^{\mathrm w}$}
\put(63,57){$s_q^{\mathrm w}$}
\end{picture}}
\end{figure}
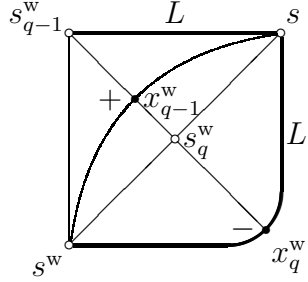
It should be at least two-valent, which implies $q\geqslant3$.

Let $r$ be the smallest integer with $2\leqslant r\leqslant q-1$ such
that $s^{\mathrm w}_r\in\partial M=L$. As we have just seen, the existence
of such an $r$ follows from~(\ref{abc}).
Let us cut the disk $M$ along a regular fibre connecting $s^{\mathrm w}$
and $s^{\mathrm w}_r$ into two disks. Let $d$ be the one of the two
that contains $s^{\mathrm w}_1$, see Fig.~\ref{diskd}.

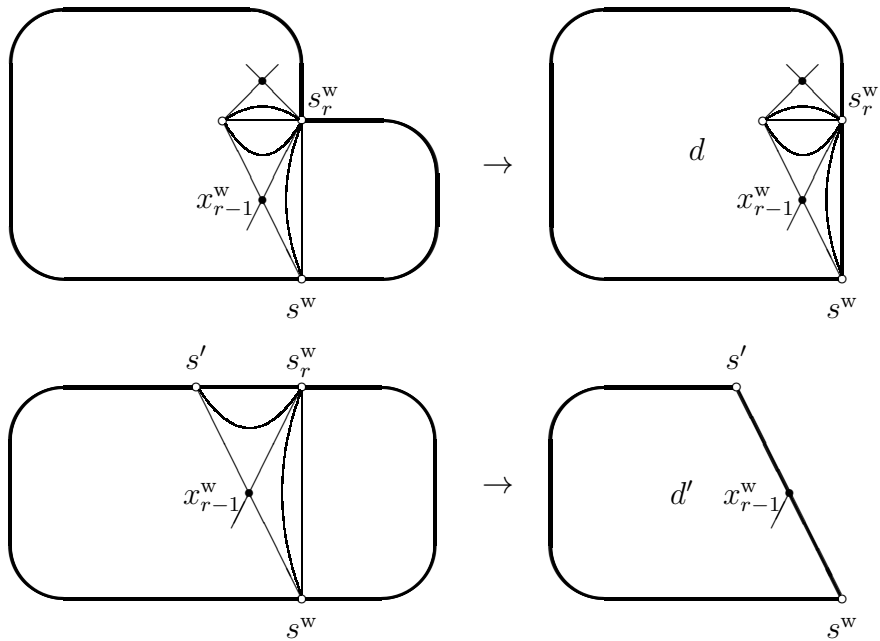
\begin{figure}
\caption{The smaller disks $d$ and $d'$}\label{diskd}
\vskip2em
\centerline{
\begin{picture}(160,120)
\put(110,20){\circle3}
\put(110,80){\circle3}
\bezier250(80.8,78.6)(95,55)(109,78.6)
\bezier200(81,80.5)(95,90)(108.5,80.5)
\bezier300(109.6,21.5)(98,50)(109.6,78.5)
\put(110,21.5){\line(0,1){57}}
\multiput(111.5,49.5)(0,0.2)5{\oval(98.5,60)[r]}
\multiput(111.5,50)(0.2,0)5{\oval(98.5,60)[r]}
\multiput(54.5,81.5)(0.2,0)5{\oval(110,80)[t]}
\multiput(55,81.5)(0,0.2)5{\oval(110,80)[t]}
\multiput(108.5,100)(-0.2,0)5{\oval(216.4,160)[bl]}
\multiput(108.5,99.5)(0,0.2)5{\oval(216.4,160)[bl]}
\put(80,80){\circle3}
\put(81,81){\line(1,1){20}}
\put(109,81){\line(-1,1){20}}
\put(95,95){\circle*3}
\put(80.7,78.4){\line(1,-2){28.6}}
\put(109.3,78.4){\line(-1,-2){20}}
\put(95,50){\circle*3}
\put(81.5,80){\line(1,0){27}}
\put(110,5){\hbox to0pt{\hss$s^{\mathrm w}$\hss}}
\put(70,47){$x^{\mathrm w}_{r-1}$}
\put(112,85){$s^{\mathrm w}_r$}
\end{picture}
\hskip0.5cm\raisebox{60pt}{$\rightarrow$}\hskip0.5cm
\begin{picture}(120,120)
\bezier250(80.8,78.6)(95,55)(109,78.6)
\bezier200(81,80.5)(95,90)(108.5,80.5)
\bezier300(109.6,21.5)(98,50)(109.6,78.5)
\put(110,20){\circle3}
\put(110,80){\circle3}
{\linethickness{1pt}\put(110,21.5){\line(0,1){57}}}
\multiput(54.5,81.5)(0.2,0)5{\oval(110,80)[t]}
\multiput(55,81.5)(0,0.2)5{\oval(110,80)[t]}
\multiput(108.5,100)(-0.2,0)5{\oval(216.4,160)[bl]}
\multiput(108.5,99.5)(0,0.2)5{\oval(216.4,160)[bl]}
\put(80,80){\circle3}
\put(81,81){\line(1,1){20}}
\put(109,81){\line(-1,1){20}}
\put(95,95){\circle*3}
\put(80.7,78.4){\line(1,-2){28.6}}
\put(109.3,78.4){\line(-1,-2){20}}
\put(95,50){\circle*3}
\put(81.5,80){\line(1,0){27}}
\put(110,5){\hbox to0pt{\hss$s^{\mathrm w}$\hss}}
\put(70,47){$x^{\mathrm w}_{r-1}$}
\put(112,85){$s^{\mathrm w}_r$}
\put(55,65){\hbox to0pt{\hss$d$\hss}}
\end{picture}
}
\centerline{
\begin{picture}(160,120)
\multiput(111.5,59.5)(0,0.2)5{\oval(97,80)[r]}
\multiput(111.5,60)(0.2,0)5{\oval(97,80)[r]}
\bezier300(71,99)(90,70)(109,99)
\bezier400(109.5,21.5)(95,60)(109.5,98.5)
{\linethickness{1.2pt}\put(71.5,100){\line(1,0){37}}}
\put(110,20){\circle3}
\put(110,100){\circle3}
\put(70,100){\circle3}
\put(110,21.5){\line(0,1){77}}
\put(70.7,98.6){\line(1,-2){38.6}}
\put(109.3,98.6){\line(-1,-2){26}}
\put(90,60){\circle*3}
\put(110,5){\hbox to0pt{\hss$s^{\mathrm w}$\hss}}
\put(110,107){\hbox to0pt{\hss$s^{\mathrm w}_r$\hss}}
\put(70,107){\hbox to0pt{\hss$s'$\hss}}
\multiput(68.5,60)(-0.2,0)5{\oval(137,80)[tl]}
\multiput(68.5,59.5)(0,0.2)5{\oval(137,80)[tl]}
\multiput(108.5,60)(-0.2,0)5{\oval(217,80)[bl]}
\multiput(108.5,59.5)(0,0.2)5{\oval(217,80)[bl]}
\put(65,56){$x^{\mathrm w}_{r-1}$}
\end{picture}
\hskip0.5cm\raisebox{60pt}{$\rightarrow$}\hskip0.5cm
\begin{picture}(120,120)
\put(110,20){\circle3}
\put(70,100){\circle3}
\multiput(70.7,98.6)(0.2,0.1)3{\line(1,-2){38.6}}
\multiput(70.7,98.6)(-0.2,-0.1)3{\line(1,-2){38.6}}
\put(90,60){\line(-1,-2){6.7}}
\put(90,60){\circle*3}
\put(110,5){\hbox to0pt{\hss$s^{\mathrm w}$\hss}}
\put(70,107){\hbox to0pt{\hss$s'$\hss}}
\multiput(68.5,60)(-0.2,0)5{\oval(137,80)[tl]}
\multiput(68.5,59.5)(0,0.2)5{\oval(137,80)[tl]}
\multiput(108.5,60)(-0.2,0)5{\oval(217,80)[bl]}
\multiput(108.5,59.5)(0,0.2)5{\oval(217,80)[bl]}
\put(65,56){$x^{\mathrm w}_{r-1}$}
\put(45,56){$d'$}
\end{picture}
}
\end{figure}

If $s^{\mathrm w}$ becomes a univalent boundary vertex of $d$
(\emph{i.e.}, $r=2$),
then its star $\sigma^{\mathrm w}\cap d$ in $d$ contains
no interior vertices, and $\mathcal F|_d$ has
at most two univalent boundary vertices, so we can apply the argument
of case (iii.a) to show that there is a good two- or three-valent
interior vertex in $d$. Let us assume that the valence of
$s^{\mathrm w}$ in $d$ is larger than one.

If $s^{\mathrm w}_r$ is not a univalent boundary vertex of $d$,
then we can apply the argument of case (iii.b) to $d$ without any change.

If $s^{\mathrm w}_r$ is a univalent boundary vertex of $d$,
let $s'$ be the vertex adjacent
to $x^{\mathrm w}_{r-1}$ but different from $s^{\mathrm w}$,
$s^{\mathrm w}_{r-1}$, $s^{\mathrm w}_r$.
Let us cut off a ``triangle'' whose ``hypotenuse'' consists
of separatrices connecting $x^{\mathrm w}_{r-1}$
to $s^{\mathrm w}$ and to $s'$ from the disk $d$.
By assumption, after the cutting,
$s^{\mathrm w}$ and $s'$ become boundary vertices
of valence at least two.
Thus, the restriction of $\mathcal F$ to the smaller disk
has a boundary saddle
and no univalent boundary vertices. So, we can repeat our
reasoning above for this disk and show that $s^{\mathrm w}_{r-1}$
is a boundary vertex, a contradiction.

(iii.d) Finally, suppose that some of the vertices
$s^{\mathrm w}_i$, $i=1,\dots,q$ coincide.
Choose $s^{\mathrm w}_i=s^{\mathrm w}_j$ with $j-i>0$
smallest possible. Let $d$ be a disk enclosed
by two non-parallel fibres connecting $s^{\mathrm w}$
with $s^{\mathrm w}_i=s^{\mathrm w}_j$ such that
$s^{\mathrm w}_k\in d$ if and only if $i<k<j$.
The disk $d$ has only two boundary vertices.
The vertex $s^{\mathrm w}_i$ cannot be a univalent
boundary vertex of $d$ because otherwise we would have
$s^{\mathrm w}_{i+1}=s^{\mathrm w}_{j-1}$, so the argument of case (iii.b)
applies to $d$.
\end{proof}

\subsection{Simplifying the characteristic surface}
This is the final part of the proof of
Propositions~\ref{triv}--\ref{comp}
and Theorem~\ref{th1}.

\begin{lemma}\label{l5}
Let $L$ be an arc-presentation satisfying the
assumptions of one of the Propositions\/~{\rm\ref{triv}--\ref{comp}},
$M$ an admissible characteristic surface for $L$ such that $c(M)>0$.
Then there exist an arc-presentation $L'$ and an admissible characteristic
surface $M'$ for $L'$ such that one of the following takes place
\rm\begin{enumerate}
\item\it
$L'$ is obtained from $L$ by finitely many exchange moves,
and we have $c(M')<c(M)$\rm;
\item\it $L'$ that is obtained from $L$ by finitely many exchange moves
and one destabilization move.
\end{enumerate}
\end{lemma}

\begin{proof}
It suffices to show how to obtain $L',M'$ in each
case listed in Lemma~\ref{cases}. In Cases~1, 3, 4 we just
follow~\cite{bm4,bm5,cromwell} and show additionally that the transformation
$L\mapsto L'$ can be decomposed into elementary moves.
Case~2 is specific in our situation.
In cases 5, 6, 7, and 8 our actions will be similar to those
in Cases 3 and 4. In Case~5,
we also consider a situation not covered by~\cite{bm4,bm5,cromwell}.
Two tricks play important r\^ole here:
changing the valence of a vertex (Cases~4 through 8) and reducing
a pair of vertices one of which is two-valent (Cases 3 and 5)
by using ``global exchange'' from Lemma~\ref{bigexchange}.
These tricks are the same in nature as those that were
used by Bennequin's in~\cite{ben} to prove an inequality
relating the knot genus, the number of strands, and the crossing
number of a braid representing the knot.

\begin{rem}
The Reeb foliation in $S^3$, which Bennequin considers,
is different from ours and from the one used in~\cite{bm4,bm5,cromwell}.
Our foliation $d\theta=0$ can be obtained from the Reeb foliation
by contracting one of the solid tori bounded by the compact leaf
into a circle, which becomes our circle $S^1_\varphi$. This difference
between foliations turns out to be unimportant for the arguments
that are used to show that the foliation $\mathcal F$
can be modified in a prescribed way once we know that it
contains a certain pattern.
\end{rem}

\medskip\noindent{\it Case\/}~1: $\cal F$ has a pole.\\
We remove all poles and closed fibres of the foliation $\cal F$.
The obtained surface has new boundary components,
each consisting of a saddle and a separatrix and lying
in a page $\disk_t$. For such boundary
component $\gamma$, we attach a disk bounded by $\gamma$
that lies in $\disk_t$ to the surface and then
deform the result slightly to obtain a smooth surface.
The saddle singularity at $\gamma$ disappears (see Fig.~\ref{loop}).
\begin{figure}[ht]\caption{Removing poles and closed fibres}\label{loop}
\centerline{\begin{picture}(70,70)
\put(50,0){\line(0,1){35}}\put(70,20){\line(-1,0){35}}
\put(50,20){\circle*{3}}\put(35,35){\oval(30,30)[t]}
\put(35,35){\oval(30,30)[bl]}\put(35,35){\oval(40,40)[t]}
\put(35,35){\oval(40,40)[bl]}\put(35,0){\oval(20,30)[tr]}
\put(70,35){\oval(30,20)[bl]}\put(70,0){\oval(30,30)[tl]}
\end{picture}\hskip40pt
\begin{picture}(70,70)
\put(50,0){\line(0,1){10}}\put(70,20){\line(-1,0){10}}
\put(35,35){\oval(40,40)[t]}\put(40,10){\oval(20,20)[tr]}
\put(60,30){\oval(20,20)[bl]}\put(40,30){\oval(20,20)[t]}
\put(40,30){\oval(20,20)[bl]}\put(35,35){\oval(40,40)[bl]}
\put(35,0){\oval(20,30)[tr]}\put(70,35){\oval(30,20)[bl]}
\put(70,0){\oval(30,30)[tl]}
\end{picture}}
\end{figure}
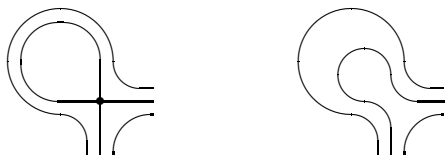
The complexity of the surface is decreased, the surface
remains admissible, the link $L$ is untouched.

In what follows we assume that there are no poles, and hence,
no closed regular fibres in $\cal F$.

\medskip\noindent{\it Case\/}~2: $\cal F$ has a univalent
boundary vertex.\\
Let $s_1$ be a univalent boundary vertex, $s_2$, $s_3$
the neighbouring vertices to $s_1$ at the knot $L$,
$\alpha$ and $\beta$ the two arcs of $L$
that connect $s_1$ with $s_2$ and $s_3$, respectively, and
$x$ the saddle connected with $s_1$ by a separatrix.
The saddle $x$ may be a boundary saddle as well as
an interior one.

Suppose $x$ is a boundary saddle. We may assume without loss
of generality that $x\in\beta$. The behaviour of the
foliation near $s_1$ is shown in Fig.~\ref{univalent} on the left.
\begin{figure}[ht]\caption{A univalent vertex of $\cal F$}\label{univalent}
\vskip20pt
\centerline{
\begin{picture}(80,120)
\put(20,20){\circle{3}}\put(60,20){\circle{3}}
\put(60,60){\circle*{3}}\put(60,100){\circle{3}}
{\linethickness{1pt}\put(0,20){\line(1,0){18}}
\put(22,20){\line(1,0){36}}\put(60,22){\line(0,1){76}}
\put(60,102){\line(0,1){18}}}\put(21,22){\line(1,2){38}}
\put(21,21){\line(1,1){38}}
\bezier300(22,21)(60,60)(59,21)
\bezier300(21,22)(60,60)(59,98)
\put(17,5){$s_2$}\put(57,5){$s_1$}\put(64,58){$x$}
\put(64,98){$s_3$}\put(37,10){$\alpha$}
\put(70,50){$\beta$}\end{picture}\hskip40pt
\begin{picture}(120,120)\put(20,20){\circle{3}}
\put(60,20){\circle{3}}\put(60,60){\circle*{3}}
\put(60,100){\circle{3}}{\linethickness{1pt}\put(0,20){\line(1,0){18}}
\put(62,20){\line(1,0){36}}\put(102,20){\line(1,0){18}}
\put(22,20){\line(1,0){36}}}\put(60,22){\line(0,1){76}}
\put(60,102){\line(0,1){18}}
\put(17,5){$s_2$}
\put(97,5){$s_3$}\put(57,5){$s_1$}\put(63,50){$x$}
\put(37,10){$\alpha$}\put(77,10){$\beta$}\put(100,20){\circle{3}}
\bezier300(21,22)(40,60)(60,60)
\bezier300(60,60)(80,60)(99,22)
\bezier300(20,22)(20,60)(59,99)
\bezier300(100,22)(100,60)(61,99)
\bezier300(21,21)(55,60)(59,21)
\bezier300(99,21)(65,60)(61,21)
\bezier300(21,22)(25,40)(45,60)
\bezier300(45,60)(55,70)(59,98)
\bezier300(99,22)(95,40)(75,60)
\bezier300(75,60)(65,70)(61,98)
\end{picture}}
\end{figure}
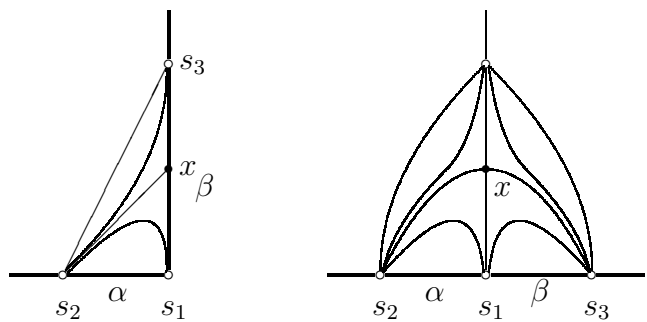
There is no topological obstruction to slide the arc $\alpha$ over $M$
toward $x$ in such a way that, at any moment, $\alpha$ coincides
with a regular fibre of $\cal F$ with endpoints $s_1,s_2$.

As a result of the sliding
finitely many arc exchange moves will occur.
We move $\alpha$
until $\alpha$ and $\beta$ become neighbouring arcs.
When $\alpha$ and $\beta$ become neighbours, we
apply a generalized destabilization move.

If $x$ is not a boundary saddle, then we do the
same thing, but with both arcs $\alpha$ and $\beta$.
The star of $s_1$ is shown in Fig.~\ref{univalent} on
the right. We slide $\alpha$ and $\beta$ toward
$x$ until it becomes possible to apply a generalized
destabilization move.

\begin{rem}
Note that, once we have applied a destabilization move
and the arc-presentation $L$ has become simpler,
we forget the surface $M$ and search
for a completely new admissible characteristic surface for the
obtained arc-presentation. However, the new characteristic
surface can be found by a modification of the old one,
which may yield in some optimization of the simplification
procedure. We do not discuss this in detail because the
optimization questions are out of the scope of this paper.
\end{rem}

\medskip\noindent{\it Case\/}~3: $\cal F$ has a good two-valent
interior vertex.\\
Let $s_1$ be a good two-valent interior vertex and let
$x_1$, $x_2$
be the two saddles from the star of $s_1$.
There must be also two vertices in the star of $s_1$.
Denote them by $s_2$, $s_3$. We may assume
without loss of generality that $s_2$ is also an interior vertex.
The star of $s_1$ is shown in Fig.~\ref{2-v}.
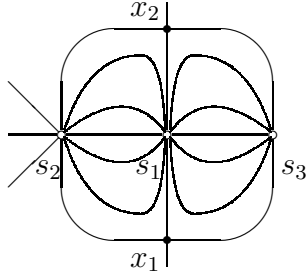
\begin{figure}[ht]\caption{A good two-valent vertex of $\cal F$}
\label{2-v}
\centerline{\begin{picture}(110,120)
\put(20,60){\circle{3}}\put(60,60){\circle{3}}\put(100,60){\circle{3}}
\put(60,20){\circle*{3}}\put(60,100){\circle*{3}}\put(0,60){\line(1,0){18}}
\put(22,60){\line(1,0){36}}\put(62,60){\line(1,0){36}}
\put(60,58){\oval(80,76)[b]}\put(60,62){\oval(80,76)[t]}
\put(60,10){\line(0,1){48}}\put(60,62){\line(0,1){48}}
\put(0,40){\line(1,1){19}}\put(0,80){\line(1,-1){19}}
\bezier150(21,62)(30,90)(50,90)\bezier150(50,90)(58,90)(59,62)
\bezier150(21,61)(45,80)(59,61)\bezier150(99,62)(90,90)(70,90)
\bezier150(70,90)(62,90)(61,62)\bezier150(99,61)(75,80)(61,61)
\bezier150(21,58)(30,30)(50,30)\bezier150(50,30)(58,30)(59,58)
\bezier150(21,59)(45,40)(59,59)\bezier150(99,58)(90,30)(70,30)
\bezier150(70,30)(62,30)(61,58)\bezier150(99,59)(75,40)(61,59)
\put(10,45){$s_2$}\put(48,45){$s_1$}\put(103,45){$s_3$}
\put(46,10){$x_1$}
\put(46,105){$x_2$}
\end{picture}}
\end{figure}
Let $t_1=\theta(x_1)$, $t_2=\theta(x_2)$. We may also assume
that the section $\disk_t\cap M$
contains the arc $\gamma_{t,s_1,s_2}$ when $t\in(t_1,t_2)$
and the arc $\gamma_{t,s_1,s_3}$ when $t\in(t_2,t_1)$, and
we have $s_1\in(s_2,s_3)$. The other cases are obtained
by flipping the orientation of $S^1_\varphi$ and/or $S^1_\theta$.

It was pointed out by Bill Menasco and Adam Sikora that
it is important to know here that $s_3$ is not a winding
vertex. That's why we consider all interior vertices in the star
of a winding vertex as bad by definition.

Pick a small $\varepsilon>0$. No arc
in $\disk_t\cap(L\cup M)$ will be interleaved
with the arc $(s_2+\varepsilon)(s_1-\varepsilon)$ if $t\in[t_1,t_2]$
and with $(s_1-\varepsilon)(s_3+\varepsilon)$ if $t\in(t_2,t_1)$.
Clearly, if no arc in $\disk_t\cap(L\cup M)$
is interleaved with an arc $ss'$, then this remains
true under a small variation of $t$. So, if $\varepsilon$ is
small enough, then no connected component
of $\disk_t\cap(L\cup M)$ is interleaved
with the arc $(s_2+\varepsilon)(s_1-\varepsilon)$
if $t\in[t_1-\varepsilon,t_2+\varepsilon]$
and with $(s_1-\varepsilon)(s_3+\varepsilon)$ if
$t\in[t_2+\varepsilon,t_1-\varepsilon]$.
Thus, the arc-presentation $L$ and the surface $M$
satisfy the assumptions of Lemma~\ref{bigexchange}, where
we should replace $s_1,s_2,s_3,t_1,t_2$ by
$s_2+\varepsilon,s_1-\varepsilon,s_3+\varepsilon,t_1-\varepsilon,
t_2+\varepsilon$, respectively, and we can
exchange the intervals $(s_2+\varepsilon,s_1-\varepsilon)$
and $(s_1-\varepsilon,s_3+\varepsilon)$.

The new arc-presentation $L'$ is obtained from $L$ by
a generalized exchange move, which, according to
Proposition~\ref{genexchange}, is the composition of
ordinary exchange moves. As a result of the exchange,
the vertex $s_1$ is moved to $s_1'=s_2+2\varepsilon$.
There is no other vertex of $L$ or $\cal F$
in $(s_2,s_1')$. The foliation of the surface is not
changed, therefore, there still exists a regular fibre $\alpha$ connecting
the vertices $s_2$ and $s_1'$. We can now isotop
the surface so that the obtained surface $M'$
is admissible and we have $|M'\cap S^1_\varphi|=
|M\cap S^1_\varphi|-2$. In order to see that
there is no obstruction, we can do the following.
First, we replace the interval $(s_2-\varepsilon,s_1'+\varepsilon)$
of the binding circle by an arc parallel to $\alpha$
(see Fig.~\ref{-2}). Then we apply an isotopy
that will restore the binding circle.

\begin{figure}[ht]\caption{Simplifying the surface $M$}\label{-2}
\centerline{
\begin{picture}(350,60)
\put(0,0){\begin{picture}(95,60)
\put(25,30){\circle{3}}
\put(55,30){\circle{3}}
\put(0,30){\line(1,0){23}}
\put(27,30){\line(1,0){26}}
\put(57,30){\line(1,0){23}}
\put(10,0){\line(1,2){14}}
\put(70,0){\line(-1,2){14}}
\bezier200(26,32)(40,60)(54,32)
\put(82,25){$S^1_\varphi$}
\put(15,0){$M$}
\put(47,45){$\alpha$}
\put(15,35){$s_2$}
\put(57,35){$s_1'$}
\end{picture}}
\put(115,28){\hbox to0pt{\hss$\rightarrow$\hss}}
\put(135,0){\begin{picture}(80,60)
\put(0,30){\line(1,0){18}}
\put(62,30){\line(1,0){18}}
\bezier200(18,30)(40,75)(62,30)
\put(10,0){\line(1,2){16}}
\put(70,0){\line(-1,2){16}}
\bezier200(26,32)(40,60)(54,32)
\end{picture}}
\put(235,28){\hbox to0pt{\hss$\rightarrow$\hss}}
\put(255,0){\begin{picture}(95,60)
\put(0,30){\line(1,0){80}}
\put(10,0){\line(1,2){7}}
\put(70,0){\line(-1,2){7}}
\bezier100(17,14)(22,24)(30,24)
\bezier100(63,14)(58,24)(50,24)
\put(30,24){\line(1,0){20}}
\put(82,25){$S^1_\varphi$}
\put(15,0){$M'$}
\end{picture}}\end{picture}}
\end{figure}
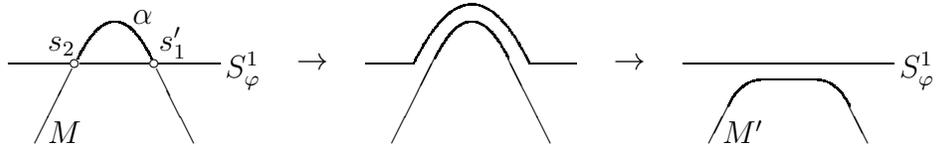

Under this isotopy, the boundary of the surface stays
unchanged. We may also assume that the obtained surface $M'$
is in general position (thus, it is still admissible) and there is no poles
of the foliation $\cal F'$ on $M'$ defined by $d\theta=0$.
Indeed, if there is one, we can apply the procedure described
in Case~1) and get rid of all poles and closed fibres of the foliation.

In the case $P=0$, from~(\ref{c(m)}) and (\ref{chi})
we have
$$c(M)=\sum\limits_kV_k+\frac12\sum\limits_kV_k^{\mathrm b}+
\frac\epsilon2-\chi.$$
Since the number of boundary vertices of $\cal F'$ is the same as
that of $\cal F$ and the number of interior vertices of $\cal F'$
is less than that of $\cal F$, we have $c(M')<c(M)$.

\begin{rem}
The operation of removing the vertices $s_2,s_1'$ has the following
effect on the combinatorial description of $M$:
in all pages $\disk_t$ with $t\in(t_1,t_2)$ the arc $s_2s_1'$
disappears, in all regular pages $\disk_t$ with $t\in(t_2,t_1)$
a pair of arcs one of which is $s_1's_3'$ and the other
has the form $s_2s$ is replaced by the arc $s_3's$,
the events at $t=t_1,t_2$ disappear, the events
in $t\in(t_2,t_1)$ are adjusted by replacing $s_2$
with $s_3'$.
\end{rem}

\medskip\noindent{\it Case\/}~4: $\cal F$ has a good three-valent
interior vertex.\\
Let $s_1$ be a good three-valent vertex and let $s_2,s_3,s_4$
be the other vertices from the star of $s_1$ enumerated
so that $s_2\in(s_1,s_3)$, $s_4\in(s_3,s_1)$.
Let $x_1$ (respectively, $x_2$) be the saddle (from the star of $s_1$)
connected by separatrices to $s_2,s_3$ (respectively, $s_3,s_4$).
The star of $s_1$ cannot contain a boundary saddle, since
otherwise there would be two boundary vertices
in the star, which would contradict to the goodness of $s_1$.
Thus, neither of the points $x_1,x_2$ is a boundary saddle. Let $s_5$ and
$s_6$ be the vertices other than $s_1,s_2,s_3,s_4$
connected by separatrices to $x_1$ and $x_2$, respectively
(see Fig.~\ref{3-v}).
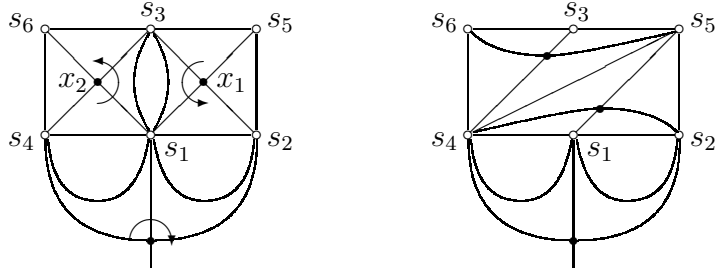
\begin{figure}[ht]\caption{Turning a three-valent vertex of $\cal F$
into a two-valent one}\label{3-v}
\centerline{
\begin{picture}(120,120)
\put(20,60){\circle{3}}
\put(60,60){\circle{3}}
\put(100,60){\circle{3}}
\put(60,20){\circle*{3}}
\put(60,20){\oval(16,16)[t]}
\put(68,18){\vector(0,-1)0}
\put(20,100){\circle{3}}
\put(60,100){\circle{3}}
\put(100,100){\circle{3}}
\put(60,10){\line(0,1){48}}
\put(20,62){\line(0,1){36}}
\put(100,62){\line(0,1){36}}
\put(22,100){\line(1,0){36}}
\put(62,100){\line(1,0){36}}
\put(22,60){\line(1,0){36}}
\put(62,60){\line(1,0){36}}
\put(21,61){\line(1,1){38}}
\put(21,99){\line(1,-1){38}}
\put(61,61){\line(1,1){38}}
\put(61,99){\line(1,-1){38}}
\put(40,80){\circle*{3}}
\put(80,80){\circle*{3}}
\put(40,80){\oval(16,16)[r]}
\put(38,88){\vector(-1,0)0}
\put(80,80){\oval(16,16)[l]}
\put(82,72){\vector(1,0)0}
\bezier200(20,58)(20,20)(60,20)
\bezier200(100,58)(100,20)(60,20)
\bezier200(21,58)(23,35)(40,35)
\bezier200(40,35)(55,35)(59,58)
\bezier200(99,58)(97,35)(80,35)
\bezier200(80,35)(65,35)(61,58)
\bezier100(59,62)(48,80)(59,98)
\bezier100(61,62)(72,80)(61,98)
\put(6,56){$s_4$}
\put(6,101){$s_6$}
\put(104,101){$s_5$}
\put(104,56){$s_2$}
\put(56,105){$s_3$}
\put(65,52){$s_1$}
\put(25,78){$x_2$}
\put(85,78){$x_1$}
\end{picture}\hskip40pt
\begin{picture}(120,120)
\put(20,60){\circle{3}}
\put(60,60){\circle{3}}
\put(100,60){\circle{3}}
\put(60,20){\circle*{3}}
\put(20,100){\circle{3}}
\put(60,100){\circle{3}}
\put(100,100){\circle{3}}
\put(60,10){\line(0,1){48}}
\put(20,62){\line(0,1){36}}
\put(100,62){\line(0,1){36}}
\put(22,100){\line(1,0){36}}
\put(62,100){\line(1,0){36}}
\put(22,60){\line(1,0){36}}
\put(62,60){\line(1,0){36}}
\put(21,61){\line(1,1){38}}
\put(61,61){\line(1,1){38}}
\bezier200(20,58)(20,20)(60,20)
\bezier200(100,58)(100,20)(60,20)
\bezier200(21,58)(23,35)(40,35)
\bezier200(40,35)(55,35)(59,58)
\bezier200(99,58)(97,35)(80,35)
\bezier200(80,35)(65,35)(61,58)
\put(22,61){\line(2,1){76}}
\put(70,70){\circle*{3}}
\put(50,90){\circle*{3}}
\bezier200(22,61)(60,70)(70,70)
\bezier200(70,70)(90,70)(99,61)
\bezier200(98,99)(60,90)(50,90)
\bezier200(50,90)(30,90)(21,99)
\put(6,56){$s_4$}
\put(6,101){$s_6$}
\put(104,101){$s_5$}
\put(104,56){$s_2$}
\put(56,105){$s_3$}
\put(65,52){$s_1$}
\end{picture}}
\end{figure}

Let $t_1=\theta(x_1)$, $t_2=\theta(x_2)$. We may assume without
loss of generality that $\theta$ increases when we go from
$x_1$ to $x_2$ transversely to the fibres connecting $s_1$
and $s_3$. In other words, we have $\gamma_{t,s_1,s_3}\subset M$
if $t\in(t_1,t_2)$. Let us now see what we know about the
flow of sections $\disk_t\cap(L, M)$ when $t$ runs
over $[t_1,t_2]$. The following events occur:

\begin{description}
\item[$t=t_1$:] the arcs $s_1s_2$ and $s_3s_5$ are replaced by
$s_1s_3$ and $s_2s_5$;
\item[$t\in(t_1,t_2)$:] some unknown events happen, but
none of them involves the arc $s_1s_3$;
\item[$t=t_2$:] the arcs $s_1s_3$ and $s_4s_6$ are replaced by
$s_1s_4$ and $s_3s_6$.
\end{description}

Together with our assumptions, this implies that
the following pairs of arcs
are non-interleaving: $(s_1s_2,s_3s_5)$, $(s_1s_3,s_2s_5)$,
$(s_1s_3,s_4s_6)$, $(s_1s_4,s_3s_6)$, which prescribes the following
circular order of the vertices in
the binding circle: $s_1,s_2,s_5,s_3,s_6,s_4$.

Since the unknown events that happen in $(t_1,t_2)$ do not
involve the arc $s_1s_3$, all the arcs participating in
them must have both endpoints either in $(s_1,s_3)$
or in $(s_3,s_1)$. In the first case, there is no topological
obstruction to moving these events ``to the future'' so that
they happen after $t_2$. In the latter case, we can
move them ``to the past'', {\it i.e.}, before~$t_1$.

As a result, finitely many arc exchange moves
are applied to $L$, and an isotopy preserving the foliation $\cal F$
to $M$. The main achievement is this:
\begin{enumerate}
\item there is no event between $t_1$ and $t_2$;
\item all three arcs $s_1s_2$, $s_3s_5$ and $s_4s_6$
are present at the moment preceding $t_1$;
\item all three arcs $s_2s_5$, $s_3s_6$ $s_1s_4$
are present right after the moment $t_2$.
\end{enumerate}

\begin{rem}
There is a minor mistake in~\cite{cromwell},
in the proof of Lemma 3: it is claimed (in different terms)
that the only obstruction to an isotopy making the events
corresponding to $x_1,x_2$ successive are arcs of $L$.
This is not necessarily so, the events in the time interval $(t_1,t_2)$
can be arbitrary, including those corresponding to saddles.
In particular, before the exchanges,
the arc $s_2s_5$, which must be present in
$\disk_{t_1+\varepsilon}\cap M$ for a small $\varepsilon$
may be already absent in $\disk_{t_2}\cap M$, being destroyed
at some moment $t\in(t_1,t_2)$.
\end{rem}

Now we replace the couple of events
$$(s_1s_2,s_3s_5)\mapsto(s_1s_3,s_2s_5)\mbox{ followed by }
(s_1s_3,s_4s_6)\mapsto(s_1s_4,s_3s_6)$$
by the following couple of events:
$$(s_3s_5,s_4s_6)\mapsto(s_3s_6,s_4s_5)\mbox{ followed by }
(s_1s_2,s_4s_5)\mapsto(s_1s_4,s_2s_5).$$
This results in an isotopy of the surface $M$
and changing the foliation as shown in Fig.~\ref{3-v}
on the right. For more explanations and figures
illustrating the change of the surface, see~\cite{bm4}.

After the change of the foliation, the vertex $s_1$ becomes
two-valent. It remains good because the star of $s_1$
becomes smaller. So, we can now proceed as in the case
of a good two-valent vertex (Case~3).

\medskip\noindent{\it Case\/}~5: $\cal F$ has two pairs two-valent
vertices, in each the vertices are connected by a fibre
intersecting $L$.
This case may occur only when $L$ is a composite link and $M$
is a factorizing sphere. This case was not considered
in full generality in~\cite{cromwell}, so, here we also fill
a gap in the proof given in~\cite{cromwell} of the
additivity of arc-index.

Let $s_1,s_2$ be a pair of two-valent
vertices connected by a fibre $\alpha$ of $\cal F$ that
intersects $L$. Since there are exactly two such pairs of vertices
and exactly two points in $L\cap M$, the union $U$ of the stars of $s_1$
and $s_2$ is pierced only once. The structure of $\cal F$
in the region $U$ is shown in Fig.~\ref{pair2-v}.
\begin{figure}[ht]\caption{A pierced region}\label{pair2-v}
\centerline{\begin{picture}(160,120)
\put(20,60){\circle{3}}
\put(60,60){\circle{3}}
\put(100,60){\circle{3}}
\put(140,60){\circle{3}}
\put(22,60){\line(1,0){36}}
\put(62,60){\line(1,0){36}}
\put(102,60){\line(1,0){36}}
\put(80,30){\circle*{3}}
\put(80,90){\circle*{3}}
\bezier400(21,61)(65,105)(80,90)
\bezier400(80,90)(90,80)(99,61)
\bezier400(21,59)(65,15)(80,30)
\bezier400(80,30)(90,40)(99,59)
\bezier400(139,61)(95,105)(80,90)
\bezier400(80,90)(70,80)(61,61)
\bezier400(139,59)(95,15)(80,30)
\bezier400(80,30)(70,40)(61,59)
\put(80,60){\vbox to0pt{\vss\hbox to0pt{\hss$*$\hss}\vss}}
\put(10,50){$s_3$}
\put(50,50){$s_1$}
\put(100,50){$s_2$}
\put(140,50){$s_4$}
\put(77,17){$x_1$}
\put(77,98){$x_2$}
\bezier350(21,61)(80,110)(60,62)
\bezier350(21,59)(80,10)(60,58)
\bezier350(139,61)(80,110)(100,62)
\bezier350(139,59)(80,10)(100,58)
\bezier150(61,61)(80,90)(99,61)
\bezier150(61,59)(80,30)(99,59)
\end{picture}}
\end{figure}
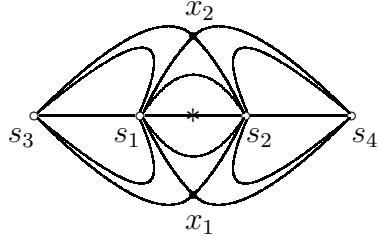
The `$*$' in the Figure stands for the point $L\cap U$. Let $\beta$
be the arc of $L$ that pierces $U$ and $x_1,x_2$ be
the saddles from~$U$. Let $t_1=\theta(x_1)$, $t_2=\theta(x_2)$,
$t_3=\theta(*)$.

Let $\partial\beta=\{s_5,s_6\}$. Since the arcs $s_1s_2$ and
$s_5s_6$ interleave, we may assume $s_5\in(s_1,s_2)$ and
$s_6\in(s_2,s_1)$. Denote by $\beta'$ the arc
of $L$ attached to $s_5$ but different from $\beta$. Let $s_3$
be the vertex from the star of $s_1$
other than $s_2$, and $s_4$ the vertex from the star of $s_2$
other than $s_1$.
We can always achieve the following by adjusting notation:
$s_1\in(s_3,s_2)$, $s_2\in(s_1,s_4)$, $t_3\in(t_1,t_2)$.

Suppose that the interval $(s_1,s_2)$ of the binding circle
does not intersect
the surface $M$. Let $A_1$ (respectively, $A_2$) be the set of arcs
in $L$ that lie in the sector $\theta\in(t_1,t_3)$ (respectively,
$\theta\in(t_3,t_2)$) and have the endpoints in $(s_1,s_2)$.
By using arc exchange moves, we can shift all the arcs from $A_1$
``to the past'' $(t_1-\varepsilon,t_1)$ and the arcs
from $A_2$ ``to the future'' $(t_2,t_2+\varepsilon)$,
keeping the relative order of arcs in each family $A_i$.

Now, for a sufficiently small $\varepsilon$, we can
exchange the intervals of the binding circle:
$(s_1+\varepsilon,s_5-\varepsilon)$ with $(s_3-\varepsilon,s_1+\varepsilon)$,
and $(s_5+\varepsilon,s_2-\varepsilon)$ with $(s_2-\varepsilon,s_4+
\varepsilon)$, since the assumptions of Lemma~\ref{bigexchange} are satisfied.
As a result of the exchanges, the vertices $s_1$, $s_2$ are moved
to $s_1'=s_5-2\varepsilon$ and $s_2'=s_5+2\varepsilon$, respectively,
becoming neighbouring to $s_5$. We can now reduce the number of
vertices of $\cal F$ in the same way as we did in Case~3.
The intersection point $*$ slides over $L$ from the arc $\beta$
to $\beta'$ (see Fig.~\ref{sliding}).
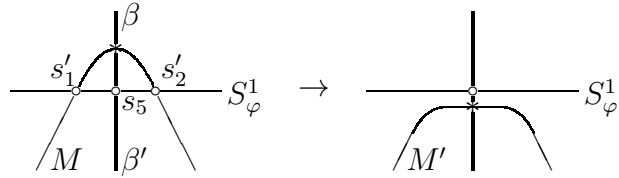
\begin{figure}[ht]\caption{Simplifying the surface $M$
near a piercing point}\label{sliding}\vskip15pt
\centerline{
\begin{picture}(230,60)
\put(0,0){\begin{picture}(95,60)
\put(25,30){\circle{3}}
\put(55,30){\circle{3}}
\put(0,30){\line(1,0){23}}
\put(27,30){\line(1,0){11}}
\put(42,30){\line(1,0){11}}
\put(40,30){\circle{3}}
{\linethickness{1pt}\put(40,0){\line(0,1){28}}
\put(40,32){\line(0,1){12}}
\put(40,48){\line(0,1){12}}}
\put(40,46){\vbox to 0pt{\vss\hbox to 0pt{\hss$*$\hss}\vss}}
\put(42,23){$s_5$}
\put(57,30){\line(1,0){23}}
\put(10,0){\line(1,2){14}}
\put(70,0){\line(-1,2){14}}
\bezier200(26,32)(40,60)(54,32)
\put(82,25){$S^1_\varphi$}
\put(15,0){$M$}
\put(42,0){$\beta'$}
\put(42,55){$\beta$}
\put(15,35){$s_1'$}
\put(57,35){$s_2'$}
\end{picture}}
\put(115,28){\hbox to0pt{\hss$\rightarrow$\hss}}
\put(135,0){\begin{picture}(95,60)
\put(0,30){\line(1,0){38}}
\put(40,30){\circle{3}}
\put(42,30){\line(1,0){38}}
\put(10,0){\line(1,2){7}}
\put(70,0){\line(-1,2){7}}
\bezier100(17,14)(22,24)(30,24)
\bezier100(63,14)(58,24)(50,24)
\put(30,24){\line(1,0){20}}
\put(40,24){\vbox to 0pt{\vss\hbox to 0pt{\hss$*$\hss}\vss}}
{\linethickness{1pt}\put(40,0){\line(0,1){22}}
\put(40,26){\line(0,1){2}}
\put(40,32){\line(0,1){28}}}
\put(82,25){$S^1_\varphi$}
\put(15,0){$M'$}
\end{picture}}\end{picture}}
\end{figure}
The rest of the argument is the same as in Case~3.

It remains to consider the situation in which $M$ intersects
the interval $(s_1,s_2)\subset S^1_\varphi$.
Actually, the procedure that has just been
considered may work in this case, too. One should include
in $A_1$, respectively $A_2$, all the events in $(t_1,t_3)$, respectively
$(t_3,t_2)$, involving arcs with endpoints in $(s_1,s_2)$.
But there is one really bad case, which really \emph{can} take place:
the arc $\beta$ can intersect the surface $M$ twice.
If the second intersection point lies in the segment of $\beta$
between $s_5$ and $*$, we cannot exchange the intervals.
We now show that, in the case when $M$
intersects $(s_1,s_2)\subset S^1_\varphi$,
the surface $M$ can be simplified by using a different method.

As we have seen at the proof of
Lemma~\ref{cases}, if none of the
Cases~1--4 occurs, then four of the vertices of
the foliation $\cal F$ are two-valent and all the others
are four-valent. It is not difficult to see that,
in this case, the surface $M$
can be cut along regular non-separatrix
fibres of $\cal F$ into two disks
and a number of annuli in which the foliation looks
as shown in Fig.~\ref{fol}.
\begin{figure}[ht]\caption{Patterns of $\cal F$ in Case~5}\label{fol}
\vskip20pt
\centerline{\begin{picture}(370,120)
\put(30,60){\circle{3}}
\put(160,60){\circle{3}}
\put(60,60){\circle{3}}
\put(130,60){\circle{3}}
\put(32,60){\line(1,0){26}}
\put(132,60){\line(1,0){26}}
\bezier600(31,59)(95,-40)(159,59)
\bezier600(31,61)(95,160)(159,61)
\bezier400(31,59)(80,-16)(129,59)
\bezier400(31,61)(80,136)(129,61)
\bezier400(61,59)(110,-16)(159,59)
\bezier400(61,61)(110,136)(159,61)
\put(95,95){\circle*{3}}
\put(95,25){\circle*{3}}
\bezier350(61,61)(95,110)(129,61)
\bezier350(61,59)(95,10)(129,59)
\bezier350(31,61)(80,110)(61,62)
\bezier350(31,59)(80,10)(61,58)
\bezier350(159,61)(110,110)(129,62)
\bezier350(159,59)(110,10)(129,58)
\bezier200(31,61)(55,80)(60,62)
\bezier200(31,59)(55,40)(60,58)
\bezier200(159,61)(135,80)(130,62)
\bezier200(159,59)(135,40)(130,58)
\bezier200(31,61)(70,115)(90,103)
\bezier50(90,103)(95,101)(100,103)
\bezier200(100,103)(120,115)(159,61)
\bezier200(31,59)(70,5)(90,17)
\bezier50(90,17)(95,19)(100,17)
\bezier200(100,17)(120,5)(159,59)
\put(1,57){$s_{2i+1}$}
\put(66,57){$s_{2i-1}$}
\put(111,57){$s_{2i}$}
\put(165,57){$s_{2i+2}$}
\put(215,60){\circle{3}}
\put(255,60){\circle{3}}
\put(235,60){\vbox to 0pt{\vss\hbox to 0pt{\hss$*\rlap{$\smash{_{_1}}$}$\hss}\vss}}
\put(310,60){\circle{3}}
\put(350,60){\circle{3}}
\put(330,60){\vbox to 0pt{\vss\hbox to 0pt{\hss$*\rlap{$\smash{_{_2}}$}$\hss}\vss}}
\bezier200(216,61)(235,95)(254,61)
\bezier200(216,61)(235,80)(254,61)
\bezier200(216,59)(235,25)(254,59)
\bezier200(216,59)(235,40)(254,59)
\put(217,60){\line(1,0){36}}
\put(202,57){$s_1$}
\put(259,57){$s_2$}
\bezier200(311,61)(330,95)(349,61)
\bezier200(311,61)(330,80)(349,61)
\bezier200(311,59)(330,25)(349,59)
\bezier200(311,59)(330,40)(349,59)
\put(312,60){\line(1,0){36}}
\put(280,57){$s_{2k-1}$}
\put(354,57){$s_{2k}$}
\end{picture}}
\end{figure}
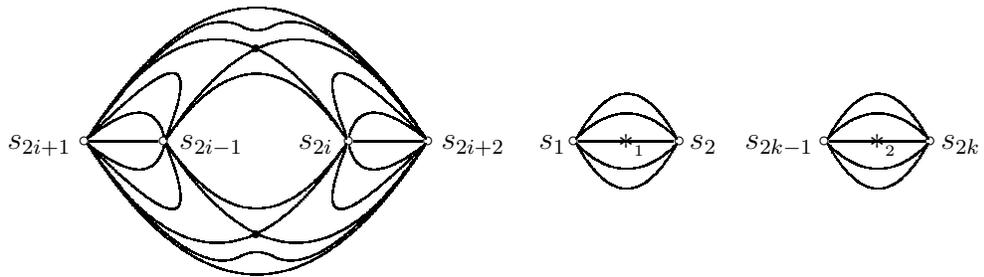
In the Figure, $i$ runs from $1$ through $k$, where
$2k$ is the total number of vertices of $\cal F$.
Certainly, we have $k\geqslant3$, since the
vertices $s_3,s_4$ do not lie in $(s_1,s_2)$
by construction. Notice that we have chosen the notation
for $s_1,s_2,s_3,s_4$ to be as before; $*_1$ and $*_2$
stand for the points $L\cap M$.

For a moment, we step back to Case~4 and Figure~\ref{3-v}.
What assumption did we use in the proof that
the foliation can be changed as shown in Fig.~\ref{3-v} on the right?
Examining the proof, one sees that we only needed
that the region filled by arcs connecting $s_1$ and $s_3$
was not pierced by $L$ and
the arcs $s_2s_5$ and $s_4s_6$ were separated from each
other in the binding circle $S^1_\varphi$
by the arc $s_1s_3$,
which is equivalent to saying that the saddles $x_1$ and $x_2$
in Figure~\ref{3-v} are coherently oriented.
In the case of a three-valent
vertex, we achieved this by appropriately numbering the
vertices $s_2,s_3$, and $s_4$. But, clearly,
we can apply the same procedure in the case
of a vertex $s_1$ such that the three vertices $s_2,s_3,s_4$
in the star of $s_1$ viewed from $s_1$ in this order
satisfy the condition $s_3\in(s_2,s_4)$.
The effect on the foliation is that the valence of $s_1$
is decreased by $1$ and the complexity of $M$
is unchanged.

Now we return to the case of the foliation shown in Fig.~\ref{fol}.
For any $i=1,\dots,k-2$ the foliation $\cal F$ have
two copies of the pattern shown in Fig.~\ref{fol1} on the left.
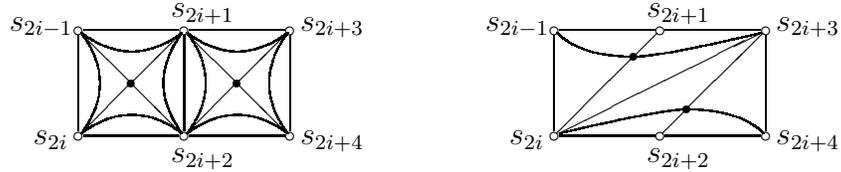
\begin{figure}[ht]\caption{Turning a four-valent vertex into
a two-valent one}\label{fol1}
\centerline{\begin{picture}(320,80)
\put(0,0){\begin{picture}(140,80)
\put(30,20){\circle{3}}
\put(70,20){\circle{3}}
\put(110,20){\circle{3}}
\put(30,60){\circle{3}}
\put(70,60){\circle{3}}
\put(110,60){\circle{3}}
\put(32,20){\line(1,0){36}}
\put(72,20){\line(1,0){36}}
\put(32,60){\line(1,0){36}}
\put(72,60){\line(1,0){36}}
\put(30,22){\line(0,1){36}}
\put(70,22){\line(0,1){36}}
\put(110,22){\line(0,1){36}}
\put(31,21){\line(1,1){38}}
\put(71,21){\line(1,1){38}}
\put(71,59){\line(1,-1){38}}
\put(31,59){\line(1,-1){38}}
\put(50,40){\circle*{3}}
\put(90,40){\circle*{3}}
\bezier150(31,21)(50,35)(69,21)
\bezier150(71,21)(90,35)(109,21)
\bezier150(31,59)(50,45)(69,59)
\bezier150(71,59)(90,45)(109,59)
\bezier150(31,21)(45,40)(31,59)
\bezier150(71,21)(85,40)(71,59)
\bezier150(109,21)(95,40)(109,59)
\bezier150(69,21)(55,40)(69,59)
\put(4,60){$s_{2i-1}$}
\put(114,60){$s_{2i+3}$}
\put(13,17){$s_{2i}$}
\put(114,17){$s_{2i+4}$}
\put(65,10){$s_{2i+2}$}
\put(65,65){$s_{2i+1}$}
\end{picture}}
\put(180,0){\begin{picture}(140,80)
\put(30,20){\circle{3}}
\put(70,20){\circle{3}}
\put(110,20){\circle{3}}
\put(30,60){\circle{3}}
\put(70,60){\circle{3}}
\put(110,60){\circle{3}}
\put(32,20){\line(1,0){36}}
\put(72,20){\line(1,0){36}}
\put(32,60){\line(1,0){36}}
\put(72,60){\line(1,0){36}}
\put(30,22){\line(0,1){36}}
\put(110,22){\line(0,1){36}}
\put(32,21){\line(2,1){76}}
\put(80,30){\circle*{3}}
\put(60,50){\circle*{3}}
\bezier200(32,21)(70,30)(80,30)
\bezier100(80,30)(100,30)(109,21)
\bezier200(108,59)(70,50)(60,50)
\bezier100(60,50)(40,50)(31,59)
\put(31,21){\line(1,1){38}}
\put(71,21){\line(1,1){38}}
\put(4,60){$s_{2i-1}$}
\put(114,60){$s_{2i+3}$}
\put(13,17){$s_{2i}$}
\put(114,17){$s_{2i+4}$}
\put(65,10){$s_{2i+2}$}
\put(65,65){$s_{2i+1}$}
\end{picture}}
\end{picture}}
\end{figure}
If the arcs $s_{2i-1}s_{2i}$ and $s_{2i+3}s_{2i+4}$
are separated from each other in $S^1_\varphi$
by the arc $s_{2i+1}s_{2i+2}$, then
we can change the foliation as shown in the Fig.~\ref{fol1}
on the right by applying the same procedure as in Case~4
to both copies. The vertices $s_{2i+1}$ and $s_{2i+2}$
will become two-valent, and we can simplify the
surface $M$ as in Case~3.

So, the only possibility left is that the
arcs $s_{2i-1}s_{2i}$ and $s_{2i+3}s_{2i+4}$
are not separated by $s_{2i+1}s_{2i+2}$ for all $i=1,\dots,k-2$.
We now show that, in this case, the interval $(s_1,s_2)$
is disjoint from $M$, which contradicts the assumption.

In order to help the imagination, we explain informally how the surface $M$
looks in the case under consideration.
Let $C$ be a simple non-closed curve in $S^3$ with endpoints in
$S^1_\varphi$ consisting of $k-1$ arcs, each lying
in a separate page $\disk_t$. Let $S$ be the boundary
of a small neighbourhood of $C$. If the neighbourhood
has been chosen appropriately, then the foliation $d\theta=0$
on $S$ looks exactly as that on $M$. Actually,
the surface $M$ can be isotoped to such an $S$
without changes in the foliation and combinatorial structure of $L$
for an appropriately chosen curve $C$.

Let $I_j$ be the following interval of the binding circle:
$$I_j=\left\{\begin{array}{ll}
(s_{2j-1},s_{2j}),&\mbox{if }j\mbox{ is odd},\\
(s_{2j},s_{2j-1}),&\mbox{if }j\mbox{ is even}.
\end{array}\right.$$
We want to show that $I_i\cap I_j=\varnothing$ if $i\ne j$,
which will imply $(s_1,s_2)\cap M=I_1\cap M=\varnothing$.
Recall that the vertices $s_1,s_2,s_3,s_4$ have been numbered
so that $I_1\cap I_2=\varnothing$.

Suppose that we have $I_i\cap I_{i+1}=\varnothing$, $i\leqslant k-2$.
There is
a saddle of $\cal F$ that is connected by separatrices to the vertices
$s_{2i+1},s_{2i+2},s_{2i+4},s_{2i+3}$ in this circular order.
These vertices must appear in $S^1_\varphi$
in the same or in the opposite circular order. Since
the arcs $s_{2i-1}s_{2i}$ and $s_{2i+3}s_{2i+4}$
are not separated by $s_{2i+1}s_{2i+2}$, we see that
the vertices $s_{2i+1},s_{2i+2},s_{2i},s_{2i-1}$
go in $S^1_\varphi$ in the same circular order
as the vertices $s_{2i+1},s_{2i+2},s_{2i+4},s_{2i+3}$ do.
Therefore, the assumption $I_i\cap I_{i+1}=\varnothing$
implies $I_{i+1}\cap I_{i+2}=\varnothing$.
Since $I_1\cap I_2=\varnothing$,
we have $I_i\cap I_{i+1}=\varnothing$ for all $i=1,\dots,k-1$.

Now, we show by induction that, for $i<j$, we have
either $I_i\cap I_j=\varnothing$ or $I_i\subset I_j$.
This is true for $i=1$, $j=2$.
Suppose, this is true for all $i<j\leqslant m$.
We have already shown that $I_m\cap I_{m+1}=\varnothing$.
Let $i\in\{1,\dots,m-1\}$. In any regular page
$\disk_t$ one of the following is present:
\begin{enumerate}
\item the arc $s_{2i-1}s_{2i}$, which separates $I_i$ from
the rest of the binding circle;
\item the arcs $s_{2i-1}s_{2i+1}$ and $s_{2i}s_{2i+2}$, which
separate $I_i\cup I_{i+1}$ from the rest of the circle;
\item the arcs $s_{2i-3}s_{2i-1}$ and $s_{2i-2}s_{2i}$, which
separate $I_i\cup I_{i-1}$ from the rest of the circle.
\end{enumerate}
In order to see this, refer to the stars of $s_{2i-1}$ and $s_{2i}$.
One of regular pages contains the arcs $s_{2m-1}s_{2m+1}$
and $s_{2m}s_{2m+2}$. By the induction hypothesis,
the vertices $s_{2m-1}$ and $s_{2m}$ are outside
of $I_{i-1}\cup I_i\cup I_{i+1}$. Therefore,
the vertices $s_{2m+1}$ and $s_{2m+2}$ are
outside of $I_i$. This implies $I_i\cap I_{m+1}=\varnothing$
or $I_i\subset I_{m+1}$.

From the symmetry argument, we conclude that,
for any $i>j$, we also have either $I_i\cap I_j=\varnothing$
or $I_i\subset I_j$. Therefore, $I_i\cap I_j=\varnothing$
for all $i\ne j$. This concludes the proof of the Lemma
in the case of a composite link.

\medskip\noindent{\it Case\/}~6: $\mathcal F$ has
a three-valent interior vertex $s^{\mathrm w}_i$ in the star
of a winding vertex $s^{\mathrm w}$. Let us assume for the moment that
there is no boundary saddle. We apply the same trick as
in Case~4, changing the valence of the vertex $s^{\mathrm w}_i$
to two. The result depends on whether we have $\delta_{i-1}=
\delta_i$ or not. Refer again to Fig.~\ref{3-v}, where $s_1$
is the vertex $s^{\mathrm w}_i$.

If $\delta_{i-1}=\delta_i$, then according to Fig.~\ref{3-v}
the winding vertex must be $s_3$. As a result of the change of
the foliation, the vertex $s_1$, which becomes two-valent,
escapes the star of $s^{\mathrm w}$. If it becomes good,
we proceed as in Case~3. If it is still bad, we can
apply the same argument as in Case (iii.b) of the proof of
Lemma~\ref{cases} to show that there is a good two- or three-valent interior
vertex or a one-valent boundary vertex somewhere else, and
proceed as in Case 2, 3, or 4 above.

If $\delta_{i-1}\ne\delta_i$, then according to Fig.~\ref{3-v}
the winding vertex must be $s_2$ or $s_4$. We may assume
that it is $s_4$. As a result of the change of the foliation
the valence of $s^{\mathrm w}$ increases, whereas the complexity
of the disk $M$ does not change. So, after finitely many such
operations either all three-valent interior vertices will
be removed from $\sigma^{\mathrm w}$ or we get a good
interior vertex of valence $\leqslant3$ or a univalent boundary vertex.

A slight complication occurs in the case of boundary saddle present,
when the three-valent vertex under consideration is $s^{\mathrm w}_q$,
and we have $\delta_{q-1}=\delta_q$.
Similarly to Case~4, one can show that now
the change of the foliation indicated in Fig.~\ref{1.5}
is possible.
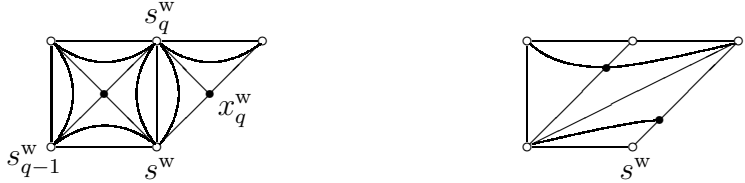
\begin{figure}[ht]\caption{The case when $s^{\mathrm w}_q$
is three-valent and $x^{\mathrm w}_{q}$ is a boundary saddle}
\label{1.5}
\centerline{\begin{picture}(320,80)
\put(0,0){\begin{picture}(140,80)
\put(30,20){\circle{3}}
\put(70,20){\circle{3}}
\put(30,60){\circle{3}}
\put(70,60){\circle{3}}
\put(110,60){\circle{3}}
\put(32,20){\line(1,0){36}}
\put(32,60){\line(1,0){36}}
\put(72,60){\line(1,0){36}}
\put(30,22){\line(0,1){36}}
\put(70,22){\line(0,1){36}}
\put(31,21){\line(1,1){38}}
\put(71,21){\line(1,1){38}}
\put(71,59){\line(1,-1){19}}
\put(31,59){\line(1,-1){38}}
\put(50,40){\circle*{3}}
\put(90,40){\circle*{3}}
\put(93,32){$x^{\mathrm w}_q$}
\bezier150(31,21)(50,35)(69,21)
\bezier150(31,59)(50,45)(69,59)
\bezier150(71,59)(90,45)(109,59)
\bezier150(31,21)(45,40)(31,59)
\bezier150(71,21)(85,40)(71,59)
\bezier150(69,21)(55,40)(69,59)
\put(13,13){$s^{\mathrm w}_{q-1}$}
\put(65,8){$s^{\mathrm w}$}
\put(65,67){$s^{\mathrm w}_q$}
\end{picture}}
\put(180,0){\begin{picture}(140,80)
\put(30,20){\circle{3}}
\put(70,20){\circle{3}}
\put(30,60){\circle{3}}
\put(70,60){\circle{3}}
\put(110,60){\circle{3}}
\put(32,20){\line(1,0){36}}
\put(32,60){\line(1,0){36}}
\put(72,60){\line(1,0){36}}
\put(30,22){\line(0,1){36}}
\put(32,21){\line(2,1){76}}
\put(80,30){\circle*{3}}
\put(60,50){\circle*{3}}
\bezier200(32,21)(70,30)(80,30)
\bezier200(108,59)(70,50)(60,50)
\bezier100(60,50)(40,50)(31,59)
\put(31,21){\line(1,1){38}}
\put(71,21){\line(1,1){38}}
\put(65,8){$s^{\mathrm w}$}
\end{picture}}
\end{picture}}
\end{figure}
The vertex $s^{\mathrm w}_q$ will escape the star of $s^{\mathrm w}$
and become two-valent, so, we can proceed as above.

\medskip\noindent{\it Case\/}~7:
$\mathcal F$ contains a two-valent vertex $s^{\mathrm w}_i$
in the star of a winding vertex, and we have $\delta_{i-1}=\delta_i$.
By the same trick as in Case~4 we can change the valence of
the vertex $s^{\mathrm w}_i$ two one, and then simplify the knot
as we did in Case~2.

\medskip\noindent{\it Case\/}~8:
$\mathcal F$ contains a four-valent interior vertex $s^{\mathrm w}_i$
in the star of a winding vertex, and we have $\delta_{i-1}=\delta_i$.
By changing the valence of the vertex $s^{\mathrm w}_i$ to three,
we also remove it from the star of the winding vertex. So, either
it becomes good and we can proceed as in Case~4, or it is still bad and
then there must be a good two- or three-valent interior vertex
or a univalent boundary vertex somewhere else.
\end{proof}

In order to complete the proof of Propositions~\ref{triv}--\ref{comp},
we need only to consider the case
$c(M)=0$. If $L$ is an arc-presentation of
the unknot, then (\ref{c(m)}) and (\ref{chi})
imply $c(L)=V_2^{\mathrm b}=2$.

If $L$ is a split or composite link, we may assume
without loss of generality that $L\cap S^1_\varphi=\{0,\pi\}$.
Then the line $x=\pi$ will be splitting or factorizing,
respectively, for the rectangular diagram of $L$.

A couple of remarks are in order.

\begin{rem}\label{strengthen}
The assertions of Propositions~\ref{split} and~\ref{comp}
can be strengthened a little by skipping destabilization moves
in the formulation. Indeed, it is trivial to show that,
for any sequence $L_0\mapsto L_1\mapsto L_2\mapsto\ldots\mapsto L_N$
of exchange and destabilization moves, there
exists a sequence $L_0\mapsto L_1'\mapsto L_2'\mapsto\ldots
\mapsto L_{N'-1}'\mapsto L_{N'}'=L_N$ in which the first $k$ moves
are exchanges and the last $N'-k$ ones destabilizations.
If $L_N$ is split or composite, then $L_k'$ is also split or
composite, respectively.
\end{rem}

\begin{rem}
One can easily show that, for any two sequences
$$D_0^{(1)}\mapsto D_1^{(1)}\mapsto\ldots\mapsto D_{N_1}^{(1)}\quad\mbox{and}
\quad D_0^{(2)}\mapsto D_1^{(2)}\mapsto\ldots\mapsto D_{N_2}^{(2)}$$
of cyclic permutation, exchange, and destabilization moves
of rectangular diagrams, there exists a sequence
$$D_0^{(1)}\#D_0^{(2)}\mapsto D_1\mapsto D_2\mapsto\ldots
\mapsto D_N=D_{N_1}^{(1)}\#D_{N_2}^{(2)}$$
of cyclic permutation,
exchange, and destabilization moves.
A similar statement is true for the distant union operation.
Thus, we can obtain a \emph{complete} decomposition of an arc-presentation
by using exchange and destabilization moves.
\end{rem}

\section{Applications}

\subsection{An algorithm for recognizing the unknot, recognizing split
links, and link factorization}
Propositions~\ref{triv}--\ref{comp} allow
to construct a simple algorithm for
recognizing the unknot,
decomposing a given link into
the distant union of non-split links, and
factorizing a non-split link to prime links.
The algorithm works as follows.

We search for all arc-presentations obtained by
finite sequences of exchange moves from
a given one. Since, for any $n$, the number of combinatorial classes
of arc-presentations of complexity $n$ is finite,
the process will terminate in finitely many steps.
Let $L_1,L_2,\ldots,L_N$ be the arc-presentations found.
If some $L_i$ admits a destabilization move $L_i\mapsto L'$,
we replace $L$ by $L'$ and proceed as before.
In finitely many steps, we get a list
of arc-presentations none of which can be simplified anymore.
Denote them again by $L_1,\dots,L_N$ and their rectangular
diagrams by $D_1,\dots,D_N$. For each $i=1,\dots,N$
we find a maximal decomposition
$$D_i=(D_{i1}\#\ldots\#D_{ip})\sqcup\ldots\sqcup
(D_{iq}\#\ldots\#D_{ir})\sqcup D_{i,r+1}\sqcup\ldots\sqcup
D_{i,N_i},$$
where all diagrams $D_{i1},\ldots,D_{ir}$ are non-trivial.
We choose $i$ so as to have maximal possible value of $N_i$.
Then all diagrams $D_{i1},\ldots,D_{ir}$ will represent non-trivial links,
and all $D_{i,r+1},\ldots,D_{iN_i}$ prime non-split links or trivial
knots.

\begin{rem}
It is natural to ask how fast the algorithm just described is.
Unfortunately, the answer is not optimistic in the sense that
the algorithm is too hard to be implemented in practice:
we cannot provide an estimation for the running time
much better than $n^{2n}$, where $n$ is the complexity
of the given diagram. However, we think that the algorithm
can be improved considerably. We describe the point very
briefly.

We say that an arc-presentation $L$ is \emph{simplifiable} if
there exists a sequence of exchange moves $L\mapsto L_1\mapsto\ldots
\mapsto L_k$ such that a destabilization move can be applied
to $L_k$. So, the problem is actually this: how to detect that
a given arc-presentation is simplifiable and find the corresponding
sequence of exchanges as fast as possible? Also, we need a fast
method to decide whether a given non-simplifiable arc-presentation
can be transformed into a split or composite one by using exchange moves.
There is a hope that a solution much better than an exhaustive search
might exist. So far, this problem has not been investigated.
\end{rem}

\subsection{An upper bound on the crossing number needed for untangling}
We shall use the notation $c_\times(D)$
for the crossing number of a planar diagram $D$.
We shall also denote by $c_{\mathrm{arc}}(D)$
the complexity of the arc-presentation corresponding to
a rectangular diagram $D$. We introduce two different
notations, because a rectangular diagram
can be viewed in both ways, as an arc-presentation and as
an ordinary planar diagram.

\begin{theorem}
Let $D_0$ be an ordinary planar diagram of either a trivial
knot, a split link, or a composite link.
Then there exists
a sequence of Reidemeister moves $D_0\mapsto D_1\mapsto\ldots
\mapsto D_N$ such that the diagram $D_N$ is trivial, split, or composite,
respectively, and we have $c_\times(D_i)\leqslant2(c_\times(D_0)+1)^2$
for all $i$.
\end{theorem}
\begin{proof}
It is shown in~\cite{cromwell2} that there is an algorithm
for converting a planar diagram $D_0$ to
an arc-presentation that has no more than $c_\times(D_0)+2$ arcs provided
that $D_0$ satisfies certain restrictions.
By using the same method, one can show that the estimation $2c_\times(D_0)+2$
works for \emph{any} connected diagram $D_0$.
Let $D^{(1)}$ be the corresponding rectangular diagram.
Following the lines of the conversion procedure,
one sees that there is a sequence
of Reidemeister moves $D_0\mapsto D_1\mapsto\ldots
\mapsto D_{i_1}=D^{(1)}$ such that $c_\times(D_i)\leqslant
c_\times(D^{(1)})$ for all $i=1,\ldots,i_1$.

Notice that a cyclic permutation of horizontal (vertical)
edges of a rectangular diagram can be decomposed
into one stabilization move, a few exchange moves, and one
destabilization move (see Fig.~\ref{decompcycl}).
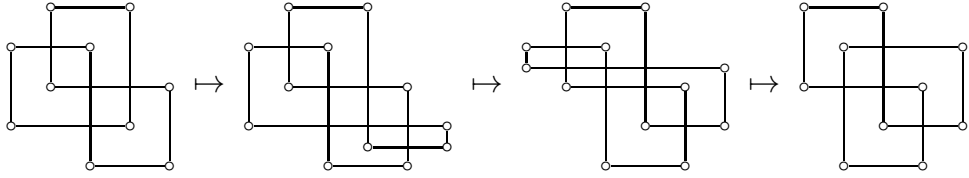
\begin{figure}[ht]\caption{Cyclic permutation via stabilization,
exchanges, and destabilization}\label{decompcycl}\vskip10pt
\centerline{\begin{picture}(360,60)
\put(0,0){\begin{picture}(60,60)
\put(0,15){\circle{3}}
\put(2,15){\line(1,0){41}}
\put(45,15){\circle{3}}
\put(45,17){\line(0,1){41}}
\put(45,60){\circle{3}}
\put(43,60){\line(-1,0){26}}
\put(15,60){\circle{3}}
\put(15,58){\line(0,-1){26}}
\put(15,30){\circle{3}}
\put(17,30){\line(1,0){41}}
\put(60,30){\circle{3}}
\put(60,28){\line(0,-1){26}}
\put(60,0){\circle{3}}
\put(58,0){\line(-1,0){26}}
\put(30,0){\circle{3}}
\put(30,2){\line(0,1){41}}
\put(30,45){\circle{3}}
\put(28,45){\line(-1,0){26}}
\put(0,45){\circle{3}}
\put(0,43){\line(0,-1){26}}
\end{picture}}
\put(90,0){\begin{picture}(75,60)
\put(0,15){\circle{3}}
\put(2,15){\line(1,0){71}}
\put(75,15){\circle{3}}
\put(75,13){\line(0,-1){4}}
\put(75,7){\circle{3}}
\put(73,7){\line(-1,0){26}}
\put(45,7){\circle{3}}
\put(45,9){\line(0,1){49}}
\put(45,60){\circle{3}}
\put(43,60){\line(-1,0){26}}
\put(15,60){\circle{3}}
\put(15,58){\line(0,-1){26}}
\put(15,30){\circle{3}}
\put(17,30){\line(1,0){41}}
\put(60,30){\circle{3}}
\put(60,28){\line(0,-1){26}}
\put(60,0){\circle{3}}
\put(58,0){\line(-1,0){26}}
\put(30,0){\circle{3}}
\put(30,2){\line(0,1){41}}
\put(30,45){\circle{3}}
\put(28,45){\line(-1,0){26}}
\put(0,45){\circle{3}}
\put(0,43){\line(0,-1){26}}
\end{picture}}
\put(195,0){\begin{picture}(75,60)
\put(0,37){\circle{3}}
\put(2,37){\line(1,0){71}}
\put(75,37){\circle{3}}
\put(75,35){\line(0,-1){18}}
\put(75,15){\circle{3}}
\put(73,15){\line(-1,0){26}}
\put(45,15){\circle{3}}
\put(45,17){\line(0,1){41}}
\put(45,60){\circle{3}}
\put(43,60){\line(-1,0){26}}
\put(15,60){\circle{3}}
\put(15,58){\line(0,-1){26}}
\put(15,30){\circle{3}}
\put(17,30){\line(1,0){41}}
\put(60,30){\circle{3}}
\put(60,28){\line(0,-1){26}}
\put(60,0){\circle{3}}
\put(58,0){\line(-1,0){26}}
\put(30,0){\circle{3}}
\put(30,2){\line(0,1){41}}
\put(30,45){\circle{3}}
\put(28,45){\line(-1,0){26}}
\put(0,45){\circle{3}}
\put(0,43){\line(0,-1){4}}
\end{picture}}
\put(285,0){\begin{picture}(75,60)
\put(32,45){\line(1,0){41}}
\put(75,45){\circle{3}}
\put(75,43){\line(0,-1){26}}
\put(75,15){\circle{3}}
\put(73,15){\line(-1,0){26}}
\put(45,15){\circle{3}}
\put(45,17){\line(0,1){41}}
\put(45,60){\circle{3}}
\put(43,60){\line(-1,0){26}}
\put(15,60){\circle{3}}
\put(15,58){\line(0,-1){26}}
\put(15,30){\circle{3}}
\put(17,30){\line(1,0){41}}
\put(60,30){\circle{3}}
\put(60,28){\line(0,-1){26}}
\put(60,0){\circle{3}}
\put(58,0){\line(-1,0){26}}
\put(30,0){\circle{3}}
\put(30,2){\line(0,1){41}}
\put(30,45){\circle{3}}
\end{picture}}
\put(75,30){\vbox to 0pt{\vss\hbox to 0pt{\hss$\mapsto$\hss}\vss}}
\put(180,30){\vbox to 0pt{\vss\hbox to 0pt{\hss$\mapsto$\hss}\vss}}
\put(285,30){\vbox to 0pt{\vss\hbox to 0pt{\hss$\mapsto$\hss}\vss}}
\end{picture}}
\end{figure}
Therefore, we can find a sequence
of (de)stabilization and exchange moves
$D^{(1)}\mapsto D^{(2)}\mapsto\ldots\mapsto D^{(m)}$
such that the final diagram $D^{(m)}$ is trivial, split,
or composite, respectively, and, for all $j=1,\dots,m$, we have
$c_{\mathrm{arc}}(D^{(j)})\leqslant 2c_\times(D_0)+3$.

It is easy to show that, for any rectangular diagram $D$,
we have
$$c_\times(D)\leqslant(c_{\mathrm{arc}}(D)-1)^2/2.$$
It remains to notice that, for any exchange
or (de)stabilization move $D^{(j)}\mapsto D^{(j+1)}$,
there exists a finite sequence of Reidemeister moves
$$D^{(j)}\mapsto D_{i_j+1}\mapsto D_{i_j+2}\mapsto\ldots
\mapsto D_{i_{j+1}}=D^{(j+1)}$$
such that
$c_\times(D_{i_j+l})\leqslant\max(c_\times(D^{(j)}),c_\times(D^{(j+1)}))$.
\end{proof}

\subsection{Results of Birman and Menasco on closed braids}
The results of~\cite{bm4} and \cite{bm5}
on closed braid representatives of links
can be deduced from the main result of this paper, which is not
a surprise because we used a modification of a
technique that was originally worked out by
Birman and Menasco for studying closed braids.

First of all, we briefly recall what the elementary moves
for braids are:
\begin{enumerate}
\item\emph{conjugation}
$$b\mapsto cbc^{-1},\quad b,c\in B_n;$$
\item\emph{stabilization}
$$b\mapsto b\sigma_n^{\pm1}\in B_{n+1},\quad b\in B_n\subset B_{n+1};$$
\item\emph{destabilization}
$$b\sigma_n^{\pm1}\mapsto b,\quad b\in B_n\subset B_{n+1};$$
\item\emph{exchange}
$$b_1\sigma_n^{\pm1}b_2\sigma_n^{\mp1}\mapsto
b_1\sigma_n^{\mp1}b_2\sigma_n^{\pm1},\quad
b_1,b_2\in B_n\subset B_{n+1}.$$
\end{enumerate}
The moves 1)--3) are Markov's moves, exchange moves were introduced by
J.\,Bir\-man and W.\,Menasco. In~\cite{bm4}, \cite{bm5}
a more general transform is called an exchange move,
so the one above is a particular case of that transform.
However, one can easily show that the more general, ``multistrand'',
exchange move can be decomposed into a finite number of
``elementary'' exchange moves and braid isotopies.

An exchange move can be decomposed
into finite sequence of Markov's moves, but not always within
the class of braids with the same or smaller number of strands.
The main result of~\cite{bm4} and \cite{bm5} is the following.

\begin{theorem}[Birman, Menasco]
If $b_0$ is a braid whose closure is a trivial knot,
then there exists a finite
sequence $b_0\mapsto b_1\mapsto\ldots\mapsto b_N$
of conjugations, destabilizations, and exchanges such that
$b_N$ is a trivial braid on one strand.

If $b_0$ is a braid whose closure is either a split
link or a composite link, there exists a finite
sequence $b_0\mapsto b_1\mapsto\ldots\mapsto b_N$
of conjugations and exchanges such that
$b_N$ is  a split braid or
a composite braid, respectively.
\end{theorem}

For the definition of split and composite braids, see~\cite{bm4}.

\begin{proof}
We give here a sketch only, inviting the reader to recover the details.

The construction of~\ref{arcbraid} applied to
a split or a composite rectangular diagram $D$ gives
a braid $b_D$ which is conjugate to a split or composite
braid, respectively.

By a straightforward check one establishes the following:
\begin{itemize}
\item if $D\mapsto D'$ is an exchange of vertical edges,
then $b_D=b_{D'}$;
\item if $D\mapsto D'$ is a cyclic permutation of edges,
then either $b_D\mapsto b_{D'}$ is a conjugation or we have
$b_D=b_{D'}$;
\item if $D\mapsto D'$ is an exchange of the upper two
horizontal edges, then $b_D\mapsto b_{D'}$ is either
a conjugation or an exchange, depending on the orientation
of the edges and their relative position (it is assumed
that we enumerate the braid strands from the bottom to the top);
\item if $D\mapsto D'$ is a destabilization move
involving the upper two horizontal edges of $D$, then
either $b_D\mapsto b_{D'}$ is a destabilization move
or we have $b_D=b_{D'}$, depending on the orientation
of the edges.
\end{itemize}

Let $b$ be a braid whose closure is a trivial knot. According
to Proposition~\ref{BD} there exists a rectangular diagram
$D$ such that $b=b_D$. From Proposition~\ref{triv} we know that
there exists a sequence
$D\mapsto D_1\mapsto\ldots\mapsto D_N$
of exchanges, cyclic permutations, and destabilizations which
takes the diagram $D$ to a trivial diagram $D_N$.

Each exchange move $D_i\mapsto D_{i+1}$ can be decomposed into
some number of cyclic permutations of horizontal edges and an exchange move
that involves the upper two edges. Similarly, each
destabilization $D_i\mapsto D_{i+1}$ can be decomposed
into some number of cyclic permutations and a destabilization
involving the upper two horizontal edges. Therefore, the braid
$b_{D_N}$, which is the unit of $B_1$, can be obtained from
the initial braid $b$ by finitely many exchanges, conjugations,
and destabilizations.

For split and composite braids, the proof is similar. One should
use Propositions~\ref{split} and \ref{comp}, respectively, and
also Remark~\ref{strengthen} (to show that destabilizations are
not needed in this case).
\end{proof}

\begin{rem}
Notice that, as a by-product, we received a new proof of Markov's theorem,
which is now a corollary to Proposition~\ref{arc-pres}.
\end{rem}

\subsection{Two tests for knottedness}

How one usually concludes that a given knot is truly knotted?
The standard way to prove that
a given knot is non-trivial is to compute some knot invariant like
the Alexander polynomial or the Jones polynomial. Such computations
seem to be hard for large knots. The algorithms known up to now
for computing polynomial invariants take at least
exponential time in the complexity of a link diagram.

However, there is (at least) one well known exception from the above
mentioned situation: if an alternating planar diagram of a knot
has no separating vertex, then it presents
a non-trivial knot or a non-split link (see~\cite{murasugi}).
To check that a given diagram is alternating and has no separating
vertex is very easy.

Here, we provide two sufficient conditions for knottedness
in terms of arc-pre\-sen\-ta\-tions (rectangular diagrams). Verification
of the conditions requires not more than quadratic time in the complexity
of the diagram provided that the diagram is encoded in a reasonable
way.

From Propositions~\ref{triv}--\ref{comp} one concludes the following

\begin{coro}\label{c1}
If no exchange and no destabilization move can be
applied to a non-trivial arc-presentation
$L$, then $L$ is a non-trivial non-split prime link.
\end{coro}

\begin{proof}
We need only to comment on why $L$ is prime, {\it i.e.}, non-composite.
The point is that, if $L$ is composite as a link,
then Proposition~\ref{comp} implies that $L$ is composite
as an arc-presentation. This contradicts the assumption,
since a composite arc-presentation always admits exchange moves.
\end{proof}

The assumption of Corollary~\ref{c1} means that
$L$ does not have ``trivial'' arcs, {\it i.e.}, arcs connecting
two neighbouring vertices, and
any two neighbouring arcs (vertices) of $L$ interleave.
In terms of the corresponding rectangular diagram,
this means that any two neighbouring vertical (horizontal)
edges interleave provided that we regard the rightmost and leftmost
(the top and the bottom) edges also as neighbouring.
Diagrams of this kind will be said to be \emph{rigid}.

There are numerous examples of rigid diagrams. Any torus knot
or link can be presented by a rigid diagram.
A rigid diagram of the Whitehead link is displayed
on the left of Fig.~\ref{wh}. Notice that, in the contrast
to the case of alternating planar diagrams, rigid and non-rigid
rectangular diagrams
of the same complexity can present equivalent prime links.
For example, the diagram of the Whitehead link
on the right of Fig.~\ref{wh} is not rigid, though it has the same
complexity as the one on the left.

In Fig.~\ref{alex0} we show a rigid diagram of
a two-component link with zero multivariable
Alexander polynomial and zero linking number. (The latter
two invariants are often used to establish non-splitness.
The example was constructed in response to a question by
Jozef Przytycki, who asked the author whether rigid diagrams
always have non-zero Alexander polynomial. So, the answer is `no'.
Recall that the non-triviality of alternating links was proved
in~\cite{murasugi} by showing that they have non-trivial Jones polynomial.)
\begin{figure}[ht]\caption{Two rectangular diagrams of the Whitehead link}
\label{wh}
\centerline{\begin{picture}(120,120)
\put(0,120){\circle{3}}\put(2,120){\line(1,0){36}}\put(40,120){\circle{3}}
\put(40,118){\line(0,-1){76}}\put(40,40){\circle{3}}
\put(42,40){\line(1,0){36}}\put(80,40){\circle{3}}\put(80,42){\line(0,1){36}}
\put(80,80){\circle{3}}\put(78,80){\line(-1,0){76}}\put(0,80){\circle{3}}
\put(0,82){\line(0,1){36}}\put(20,100){\circle{3}}\put(22,100){\line(1,0){96}}
\put(120,100){\circle{3}}\put(120,98){\line(0,-1){76}}\put(120,20){\circle{3}}
\put(118,20){\line(-1,0){56}}\put(60,20){\circle{3}}
\put(60,22){\line(0,1){36}}\put(60,60){\circle{3}}
\put(62,60){\line(1,0){36}}\put(100,60){\circle{3}}
\put(100,58){\line(0,-1){56}}\put(100,0){\circle{3}}
\put(98,0){\line(-1,0){76}}\put(20,0){\circle{3}}\put(20,2){\line(0,1){96}}
\put(60,-20){\hbox to 0pt{\hss A rigid diagram\hss}}
\end{picture}\hskip40pt\begin{picture}(120,140)\put(0,120){\circle{3}}
\put(60,-20){\hbox to 0pt{\hss A non-rigid diagram\hss}}
\put(2,120){\line(1,0){56}}\put(60,120){\circle{3}}
\put(60,118){\line(0,-1){56}}\put(60,60){\circle{3}}
\put(58,60){\line(-1,0){56}}\put(0,60){\circle{3}}\put(0,62){\line(0,1){56}}
\put(20,100){\circle{3}}\put(22,100){\line(1,0){96}}
\put(120,100){\circle{3}}\put(120,98){\line(0,-1){76}}
\put(120,20){\circle{3}}\put(118,20){\line(-1,0){36}}
\put(80,20){\circle{3}}\put(80,22){\line(0,1){56}}
\put(80,80){\circle{3}}\put(78,80){\line(-1,0){36}}
\put(40,80){\circle{3}}\put(40,78){\line(0,-1){36}}
\put(40,40){\circle{3}}\put(42,40){\line(1,0){56}}
\put(100,40){\circle{3}}\put(100,38){\line(0,-1){36}}
\put(100,0){\circle{3}}\put(98,0){\line(-1,0){76}}
\put(20,0){\circle{3}}\put(20,2){\line(0,1){96}}
\end{picture}}
\vskip20pt
\end{figure}
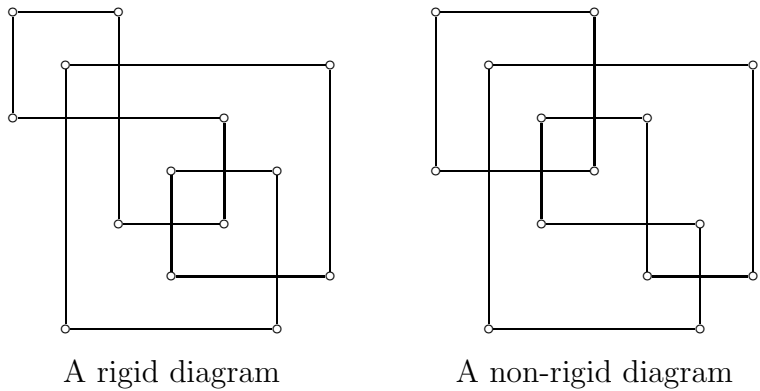
\begin{figure}\caption{A rigid diagram of a link with zero
Alexander polynomial}\label{alex0}
\vskip20pt
\centerline{\begin{picture}(90,90)
\put(10,90){\circle{3}}\put(12,90){\line(1,0){36}}
\put(50,90){\circle{3}}\put(50,88){\line(0,-1){46}}
\put(50,40){\circle{3}}\put(48,40){\line(-1,0){5}}
\put(37,40){\line(-1,0){4}}\put(27,40){\line(-1,0){4}}
\put(17,40){\line(-1,0){15}}\put(0,40){\circle{3}}
\put(0,42){\line(0,1){36}}\put(0,80){\circle{3}}
\put(2,80){\line(1,0){5}}\put(13,80){\line(1,0){5}}
\put(20,80){\circle{3}}\put(20,78){\line(0,-1){46}}
\put(20,30){\circle{3}}\put(22,30){\line(1,0){5}}
\put(33,30){\line(1,0){4}}\put(43,30){\line(1,0){15}}
\put(60,30){\circle{3}}\put(60,32){\line(0,1){36}}
\put(60,70){\circle{3}}\put(58,70){\line(-1,0){5}}
\put(47,70){\line(-1,0){24}}\put(17,70){\line(-1,0){5}}
\put(10,70){\circle{3}}\put(10,72){\line(0,1){16}}
\put(30,60){\circle{3}}\put(32,60){\line(1,0){15}}
\put(53,60){\line(1,0){4}}\put(63,60){\line(1,0){5}}
\put(70,60){\circle{3}}\put(70,58){\line(0,-1){46}}
\put(70,10){\circle{3}}\put(72,10){\line(1,0){5}}
\put(83,10){\line(1,0){5}}\put(90,10){\circle{3}}
\put(90,12){\line(0,1){36}}\put(90,50){\circle{3}}
\put(88,50){\line(-1,0){15}}\put(67,50){\line(-1,0){4}}
\put(57,50){\line(-1,0){4}}\put(47,50){\line(-1,0){5}}
\put(40,50){\circle{3}}\put(40,48){\line(0,-1){46}}
\put(40,0){\circle{3}}\put(42,0){\line(1,0){36}}
\put(80,0){\circle{3}}\put(80,2){\line(0,1){16}}
\put(80,20){\circle{3}}\put(78,20){\line(-1,0){5}}
\put(67,20){\line(-1,0){24}}\put(37,20){\line(-1,0){5}}
\put(30,20){\circle{3}}\put(30,22){\line(0,1){36}}
\end{picture}}
\end{figure}
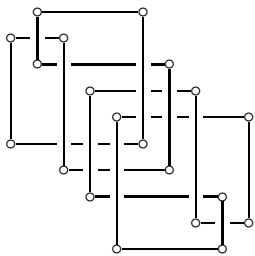

The second test for knottedness uses an analogue of the writhe.
Let us denote by $f_\varepsilon^+$ and $f_\varepsilon^-$
the self-homeomorphisms of $S^3$ defined by
$$f_\varepsilon^\pm(\varphi,\tau,\theta)=
(\varphi-\varepsilon,\tau,\theta\pm\varepsilon).$$
Let $L$ be an arc-presentation of an \emph{oriented} link.
Pick an $\varepsilon>0$ such that $\varepsilon$ is
less than the $\varphi$-distance between any two vertices of $L$
and $\theta$-distance between
any two arcs of $L$. Then the link $f_\varepsilon^\pm(L)$
will be disjoint from $L$.

We define the \emph{upper} (respectively, the \emph{lower})
\emph{writhe} of $L$ to be the linking number of $L$ with
$f_\varepsilon^+(L)$ (respectively, $f_\varepsilon^-(L)$).
We denote the upper and lower writhes of $L$
by $w_+(L)$ and $w_-(L)$, respectively. By the \emph{writhe}
of $L$ we shall mean the couple $(w_-(L),w_+(L))$.

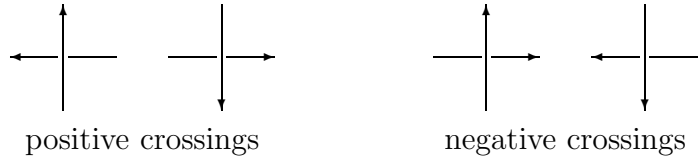
\begin{figure}[ht]\caption{Positive and negative crossings}
\label{crossings}
$$\begin{array}{ccc}
\begin{picture}(100,40)\put(20,0){\vector(0,1){40}}
\put(40,20){\line(-1,0){18}}\put(18,20){\vector(-1,0){18}}
\put(80,40){\vector(0,-1){40}}\put(60,20){\line(1,0){18}}
\put(82,20){\vector(1,0){18}}\end{picture}
&\hskip40pt&
\begin{picture}(100,40)\put(20,0){\vector(0,1){40}}
\put(0,20){\line(1,0){18}}\put(22,20){\vector(1,0){18}}
\put(80,40){\vector(0,-1){40}}\put(100,20){\line(-1,0){18}}
\put(78,20){\vector(-1,0){18}}\end{picture}\\
\mbox{positive crossings}&&\mbox{negative crossings}
\end{array}$$
\end{figure}

\begin{figure}[ht]\caption{Positive and negative turns}
\label{turns}
$$\begin{array}{ccc}
\begin{picture}(100,30)\put(20,10){\circle{3}}\put(80,20){\circle{3}}
\put(22,10){\line(1,0){18}}\put(78,20){\line(-1,0){18}}
\put(20,12){\line(0,1){18}}\put(80,18){\line(0,-1){18}}
\end{picture}
&\hskip40pt&
\begin{picture}(100,30)\put(20,20){\circle{3}}\put(80,10){\circle{3}}
\put(22,20){\line(1,0){18}}\put(78,10){\line(-1,0){18}}
\put(20,18){\line(0,-1){18}}\put(80,12){\line(0,1){18}}
\end{picture}\\
\mbox{positive turns}&&\mbox{negative turns}
\end{array}$$
\end{figure}
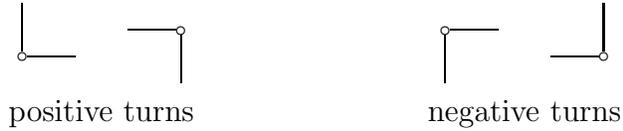

The numbers $w_\pm(L)$ can be easily computed from
the knowledge of a rectangular diagram $D$ of $L$.
We divide the set of crossings and the set of
turns of the diagram $D$ into two types,
which we call \emph{positive} and \emph{negative},
according to Figures~\ref{crossings} and \ref{turns}.

We define the \emph{writhe} $w(D)$ of a rectangular link
diagram $D$ in the standard way: $w(D)$ is the number of positive
crossings less the number of negative crossings in $D$.
Denote by $n_+(D)$ and $n_-(D)$, respectively, the number of positive
and negative turns in $D$.
The following three facts are established by a straightforward check.

\begin{prop}
Let $L$ be an arc-presentation of an oriented link, $D$ the corresponding
rectangular diagram. Then we have:
\begin{equation}\label{writhe}
w_+(L)=w(D)+\frac{n_+(D)}2,\qquad w_-(L)=w(D)-\frac{n_-(D)}2,
\end{equation}
which implies, in particular, the following:
\begin{equation}
w_+(L)-w_-(L)=c(L).
\end{equation}
\end{prop}

\begin{prop}\label{w}
The writhe of an arc-presentation does not
change under exchange moves.
\end{prop}

Equivalently,
the right-hand side of any of the formulae (\ref{writhe})
is not changed under cyclic permutation of
vertical or horizontal edges and interchange of
non-interleaved neighbouring vertical or horizontal edges of $D$.

\begin{prop}\label{w+-}
Let $L$ be an arc-presentation of an oriented link, $L\mapsto L'$
a stabilization move. Then we have either
$$w_+(L')=w_+(L)+1,\qquad w_-(L')=w_-(L')$$
or
$$w_-(L')=w_-(L)-1,\qquad w_+(L')=w_+(L).$$
\end{prop}

For the trivial arc-presentation $L$ of the unknot,
we have $w_+(L)=1$, $w_-(L)=-1$. For the distant union $L_1\sqcup L_2$
of two arc-presentations, we have $w_\pm(L_1\sqcup L_2)=w_\pm(L_1)+
w_\pm(L_2)$. In conjunction with
Propositions~\ref{triv}, \ref{split}, \ref{w},
and~\ref{w+-}, this implies the following.

\begin{coro}
If $L$ is an arc-presentation of the trivial link of $k$
components, then we have
\begin{equation}\label{test2}w_-(L)\leqslant-k,\quad k\leqslant w_+(L).
\end{equation}
\end{coro}

It happens very often that $w_-(L)\geqslant0$ or $w_+(L)\leqslant0$,
which immediately
implies that link $L$ is non-trivial.
For example, the torus link of type $(p,q)$
has an arc-presentation $L_{p,q}$ whose arcs are
$$\gamma_{-\frac{2\pi k}{p+q},\frac{2\pi k}{p+q},\frac{2\pi(k+p)}{p+q}},$$
where $k=1,\ldots,p+q$. A simple calculation gives:
$$w_-(L_{p,q})=pq-p-q,\qquad w_+(L_{p,q})=pq.$$
If both $p$ and $q$ are large, then the interval $(w_-,w_+)$
stays far away from $0$, and this remains true after a certain
number of modifications of the diagram of $L_{p,q}$.

\section*{Acknowledgement}
I am thankful to A.\,B.\,Sossinsky who read a preliminary
version of the paper and suggested many improvements, and
to J.\,Birman and S.\,V.\,Matveev for further useful remarks.

I am especially thankful to Bill Menasco and Adam Sikora for pointing
out a mistake in the earlier version of the paper and suggesting
an idea for correcting it.

The work is partially supported by Russian Foundation for
Fundamental Research, grant no.~02-01-00659,
by Leading Scientific School Support grant no.~00-15-96011,
and by LIFR MIIPT.

\end{document}